\documentclass[10pt]{amsart}
\usepackage{amsmath, amscd,amssymb}
\usepackage{amsfonts,color}
\usepackage{xypic}
\usepackage{mathrsfs}
\usepackage{cite, tocvsec2}
\usepackage{amstext, amsthm}
\usepackage{url}
\usepackage{shuffle}
\usepackage{amsopn,xspace,enumerate,multicol}
\usepackage{tikz-cd}
\usepackage[utf8]{inputenc}
\usepackage{textcomp}
\usepackage{newunicodechar}
\newunicodechar{−}{-}

\usepackage{array}

\usepackage{mathtools}
\usepackage{cellspace}

\usepackage{ulem}

\newcommand{\Lau}{\mathrm{Lau}}
\newcommand{\pol}{\mathrm{pol}}
\newcommand{\rat}{\mathrm{rat}}
\newcommand{\ser}{\mathrm{ser}}

		\newcommand{\incl}[1][r]
{\ar@<-0.2pc>@{^(-}[#1] \ar@<+0.2pc>@{-}[#1]}
		\newcommand{\eq}[1][r]
		{\ar@<-3pt>@{-}[#1]
		\ar@<-1pt>@{}[#1]|<{}="gauche"
		\ar@<+0pt>@{}[#1]|-{}="milieu"
		\ar@<+1pt>@{}[#1]|>{}="droite"
		\ar@/^2pt/@{-}"gauche";"milieu"
		\ar@/_2pt/@{-}"milieu";"droite"}
		\newcommand\lie[1]{\mathfrak{#1}}
		\newcommand\Q{{\mathbb Q}}
		
                \newtheorem{thm}{Theorem}[section]
				\newtheorem{cor}[thm]{Corollary}
				\newtheorem{prop}[thm]{Proposition} 
				\newtheorem{lem}[thm]{Lemma}

                \theoremstyle{definition}
				\newtheorem{defn}[thm]{Definition} 
				\newtheorem{rem}[thm]{Remark}

				\newtheorem{exm}[thm]{Example}

\title[]{On linearised and elliptic versions of the Kashiwara-Vergne Lie algebra} 

\author{Hidekazu Furusho}
\address{Graduate School of Mathematics, Nagoya University, Furo-cho, Chikusa-ku, Nagoya, 464-8602, Japan.}
\email{furusho@math.nagoya-u.ac.jp}

\author{Nao Komiyama}
\address{Department of Mathematics, Graduate School of Science, Osaka University Toyonaka, Osaka 560-0043, Japan}
\email{komiyama.nao.aww@osaka-u.ac.jp}
        
\author{Elise Raphael}
\address{Mathematics Department, University of Geneva, Switzerland}
\email{elise.raphael@unige.ch}

\author{Leila Schneps}
\address{Equipe Analyse Alg\'{e}brique, Institut de Math\'{e}matiques de Jussieu, Paris, France}
\email{Leila.Schneps@imj-prg.fr}

\date{January 11, 2026}

\begin{document}
                \begin{abstract}
                The goal of this article is to define a linearized or
				depth-graded version $\lie{lkrv}$, and a closely related elliptic version 
				$\lie{krv}_{ell}$, of the Kashiwara-Vergne Lie algebra $\lie{krv}$ 
				originally constructed by Alekseev and Torossian as the space of solutions 
				to the linearized Kashiwara-Vergne problem.  We show how the elliptic Lie 
				algebra $\lie{krv}_{ell}$ is related to earlier constructions 
				of elliptic versions $\lie{grt}_{ell}$ and $\lie{ds}_{ell}$ of the 
				Grothendieck-Teichm\"uller Lie algebra $\lie{grt}$ and the double shuffle 
				Lie algebra $\lie{ds}$ respectively. Based on the known relationships between 
				the three Lie algebra $\lie{grt}$, $\lie{ds}$ and $\lie{krv}$, we discuss
				the corresponding relationships between the linearized versions, and also 
				between the elliptic versions.
                \end{abstract}

\maketitle
\setcounter{tocdepth}{3}
\tableofcontents
\setcounter{section}{-1}
\section{Introduction}
				
This article studies two Lie algebras closely related to the
Kashiwara-Vergne Lie algebra $\lie{krv}$ defined in \cite{AT}: firstly, 
a depth-graded (or ``linearized'') version $\lie{lkrv}$, and secondly, 
an elliptic version $\lie{krv}_{ell}$.  
The results are
motivated by the comparison of $\lie{krv}$ with two other Lie algebras
familiar from the theory of multiple zeta values: the Grothendieck-Teichm\"uller Lie algebra $\lie{grt}$ and the double shuffle Lie algebra $\lie{ds}$.
Our definition of $\lie{lkrv}$ is an analog of the definition of the
depth-graded (or linearized) double shuffle Lie algebra $\lie{ls}$, whose structure has
given rise to many results and conjectures, in particular the famous
Broadhurst-Kreimer conjecture.
Our definition of $\lie{krv}_{ell}$ is an analog of the definition of the
elliptic double shuffle Lie algebra $\lie{ds}_{ell}$ (cf.~\cite{S3}), 
which itself is related on the one hand to $\lie{ls}$ and on
the other to the elliptic Grothendieck-Teichm\"uller Lie algebra 
$\lie{grt}_{ell}$. 
We explore all the relations between these different objects.
The main observation is that the spaces $\lie{lkrv}$ and 
$\lie{krv}_{ell}$ are defined by identical sets of properties, only applied to a different class of objects. 
This situation, which exactly
parallels the case of $\lie{ds}$ and $\lie{ds}_{ell}$, reveals a close 
and surprising relationship between the depth-graded and the elliptic versions 
of the Lie algebras $\lie{ds}$ and $\lie{krv}$, which remains invisible 
without the use of mould theory as a basic tool.

Like $\lie{grt}$ and $\lie{ds}$, the Lie algebra $\lie{krv}$ is equipped with a depth filtration; we write $gr$ for the associated graded.
We show that in analogy with the known injective map $gr\,\lie{ds}
\rightarrow \lie{ls}$, there is an injective map $gr\,\lie{krv}\hookrightarrow\lie{lkrv}$ (Proposition \ref{firstprop}).  We also show that there is an injective Lie algebra homomorphism $\lie{ls}\hookrightarrow \lie{lkrv}$ (Theorem \ref{firstthm}), and
that the parts of these spaces of depths $d=1,2,3$ are isomorphic
for all weights $n$ (Corollary \ref{firstcor}), which yields the dimensions of the bigraded parts of $\lie{lkrv}$ (and also $gr\,\lie{krv}$) of depths $1,2,3$ in all weights, since these dimensions are well-known for $\lie{ls}$. 

Passing to the elliptic situation, we define the elliptic version $\lie{krv}_{ell}$ as a subspace
of derivations of the free Lie algebra on two generators, and prove that
it is closed under the Lie bracket of derivations (Theorem \ref{krvellisLie}).
We also define an injective Lie algebra homomorphism 
                $\lie{krv}\hookrightarrow \lie{krv}_{ell} $
(in Theorem \ref{krvsection}) under the hypothetical isomorphism \eqref{Ecstatement}
in analogy with the section map $\lie{grt}\hookrightarrow\lie{grt}_{ell}$ ($\!\!$\cite{E1}) and the mould-theoretic double shuffle map $\lie{ds}\rightarrow\lie{ds}_{ell}$ ($\!\!$\cite{S3}). 
Finally, although we were not able to prove the existence of an injection $\lie{grt}_{ell}\hookrightarrow\lie{krv}_{ell}$,
we define a Lie subalgebra $\widetilde{\lie{grt}}_{ell}\subset \lie{grt}_{ell}$ such that 
the following diagram commutes (Theorem \ref{bigdiagram})
\begin{equation}\label{commdiag0}
\xymatrix{\lie{grt}\ar@{^{(}->}[r]\ar[d]&\lie{ds}\ar@{^{(}->}[r]\ar[d]&\lie{krv}\ar@{^{(}->}[d]\\
\widetilde{\lie{grt}}_{ell}\ar@{^{(}->}[r]&\lie{ds}_{ell}\ar@{^{(}->}[r]
&\lie{krv}_{ell}}
\end{equation}
Here the inclusion $\lie{grt}\hookrightarrow\lie{ds}$ is given in \cite{F2}
and $\lie{ds}\hookrightarrow\lie{krv}$ is given in \cite{S1}
\footnote{
The latter inclusion was initially contingent on an assumption, the 
validity for the senary relation \eqref{senaryrelation} for $ARI({\mathcal F}_\ser)_{\underline{al}/\underline{il}}$,
originally  asserted in \cite{Ec} without a detailed proof
(cf. \cite{FK2}).
However, recent work in \cite{S4,EF2, Ka} has unconditionally validated its validity.
\label{foot:S4 and EF2}
}
and see \S \ref{14} below for other details.
The main technique used for the constructions in this article is
the mould theory developed by J.~\'Ecalle, to which we review
in \S \ref{sec: mould theory}.

\vskip .2cm
\noindent {\bf Structure of the article.}
\S \ref{sec: Statements of main results} presents the principal results of this work.
In \S \ref{sec:2}, we reformulate the defining conditions of $\lie{krv}$, which lead to the first definition of $\lie{lkrv}$ and 
analyze the inclusion 
$gr\,\lie{krv}\hookrightarrow\lie{lkrv}$ (Proposition \ref{firstprop}).
The next section, \S \ref{sec: mould theory}, gives a brief introduction to mould theory and a translation of the defining conditions of $\lie{lkrv}$ into that language, and uses mould theory to prove 
$\lie{ls}\hookrightarrow \lie{lkrv}$ (Theorem \ref{firstthm}).
Finally, the Lie algebra structure on $\lie{grt}_{ell}$
(Theorems \ref{krvellisLie}),
the inclusion 
$\lie{krv}\hookrightarrow\lie{krv}_{ell}$
(Theorem \ref{krvsection}) and 
the commutativity of the above diagram \eqref{commdiag0}
(Theorem \ref{bigdiagram}) are given in the
three subsections of \S \ref{sec:The elliptic Kashiwara-Vergne Lie algebra}.
Appendix \ref{Appendix B} furnishes a long proof of Lemma \ref{Tnc},
a requisite step in establishing Theorem \ref{krvsection}.
                                                    
\vskip .2cm
\noindent {\bf Acknowledgements.} This article combines a fount of knowledge 
gathered mainly from work of J.~\'Ecalle, B.~Enriquez and A.~Alekseev.
We are grateful to all three for numerous helpful and interesting conversations.
H.F. has been supported by  JSPS KAKENHI Grants 24K00520 and 24K21510.
N.K. has been supported by grants JSPS KAKENHI JP23KJ1420.
													\vspace{.1cm}
                                                    
\section{Statements of main results}\label{sec: Statements of main results}
This section provides explicit statements of the main results concerning both the linearized and elliptic versions of the Kashiwara-Vergne Lie algebra.

\subsection{Special types of derivations of $\lie{lie}_2$} \label{subsec:Special types of derivations}
Let $\lie{lie}_2$ denote the degree completion of the free Lie algebra over 
$\mathbb{Q}$ on non-commutative variables $x$ and $y$. The Lie algebra 
$\lie{lie}_2$ has a weight grading by the degree (={\it weight}) of the polynomials, 
and a depth grading by the $y$-degree (={\it depth}) of the polynomials.  
We write $(\lie{lie}_2)_n$ for the graded part of weight $n$, 
$(\lie{lie}_2)^r$ for the graded part of depth $r$, and $(\lie{lie}_2)_n^r$ 
for the intersection, which is finite-dimensional.

All the Lie algebras we will study in this article (the well-known ones
$\lie{krv}$, $\lie{grt}$ and $\lie{ds}$ as well as the linearized $\lie{ls}$,
and the spaces $\lie{lkrv}$ and $\lie{krv}_{ell}$ that we introduce) can 
be viewed either as Lie subalgebras of particular subalgebras of the 
derivations of $\lie{lie}_2$, equipped with the bracket of derivations,
or as subspaces of $\lie{lie}_2$ equipped with particular Lie brackets coming 
from the Lie bracket of derivations.  
Both ways of considering our spaces 
are natural and useful, and we go back and forth between them as
convenient for our proofs.

Let $\lie{der_2}$ denote the algebra of derivations on $\lie{lie}_2$. 
It is a Lie algebra under the Lie bracket given by the commutator of 
derivations. For $a,b\in \lie{lie}_2$, we write $D_{b,a}$ for the derivation 
defined by $x\mapsto b$ and $y\mapsto a$.  The bracket is explicitly given by 
\begin{equation}\label{dabbracket}
[D_{b,a},D_{b',a'}]=D_{\tilde b,\tilde a}
\end{equation}
with 
\begin{equation}\label{bbracket}\tilde b=D_{b,a}(b')-D_{b',a'}(b), \ \ \ \ 
\tilde a=D_{b,a}(a')-D_{b',a'}(a).
\end{equation}

\vspace{.2cm}
\noindent $\bullet$ Let $\lie{oder}_2$ denote the Lie subalgebra 
of $\lie{der}_2$ of derivations $D=D_{b,a}$ that annihilate the bracket $[x,y]$
and such that neither $D(x)$ nor $D(y)$ have a linear term in $x$.
The map $\lie{oder}_2\rightarrow \lie{lie}_2$ given by $D\mapsto D(x)$
is injective (see Corollary \ref{liepush}).

\vspace{.1cm}
\noindent $\bullet$ Let $\lie{tder}_2$ denote the Lie subalgebra of 
$\lie{der}_2$ of {\it tangential derivations}, 
which are the derivations $E_{a,b}$ for elements $a,b\in\lie{lie}_2$ such
that $a$ has no linear term in $x$ and $b$ has no linear term in $y$, 
such that 
$$E_{a,b}(x)=[x,a],\ \ E_{a,b}(y)=[y,b].$$ 
The Lie bracket is explicitly given by
\begin{equation}\label{eabbracket}
[E_{a,b},E_{a',b'}]=E_{\tilde a,\tilde b}
\end{equation}
where 
\begin{equation}\label{tderbracket}
\tilde a=[a,a']+E_{a,b}(a')-E_{a',b'}(a),\ \ \ 
\tilde b=[b,b']+E_{a,b}(b')-E_{a',b'}(b).
\end{equation}

\vspace{.1cm}
{
Let $\lie{tder}_2^{(x)}$
be the Lie subalgebra of $\lie{tder}_2$
consisting of tangential derivations annihilating $x$, }
i.e. those
of the form $$d_b=E_{0,b}.$$  The derivation $d_b$ is defined by its values on $x$ and $y$
\begin{equation}\label{Iharader}
d_b(x)=0,\ \ \ \ d_b(y)=[y,b].
\end{equation}
The Lie bracket on {$\lie{tder}_2^{(x)}$} is given by $[d_b,d_{b'}]=d_{\{b,b'\}}$,  
where $\{b,b'\}$ is the {\it Poisson} (or Ihara) bracket given by
\begin{equation}\label{Poisson}
\{b,b'\}=[b,b']+d_b(b')-d_{b'}(b),
\end{equation}
i.e.~the second term of \eqref{tderbracket}.

\vspace{.1cm}
\noindent $\bullet$ {Let $\lie{sder}_{2}^{(x)}$ be the Lie subalgebra of 
$\lie{tder}_2^{(x)}$
consisting of derivations $D\in\lie{tder}_2$ such that 
$$
D(x)=0, \quad  D(y)=[y,b],\quad   D(z)=[z, c]
$$
with $z=-x-y$ 
for some $b,c\in\lie{lie}_2$.}

\vspace{.1cm}
\noindent $\bullet$ 
{
Similarly we define
$\lie{tder}_2^{(y)}$
to be the Lie subalgebra of $\lie{tder}_2$
consisting of tangential derivations annihilating $y$,
i.e. those
of the form $E_{b,0}$. 
The Lie bracket is given by 
\begin{equation}\label{eq: opposite Poisson}
[E_{b,0},E_{b',0}]=E_{\{b,b'\}^o,0}
\quad \text{ with }\quad
    \{b,b'\}^o:=\{b',b\}.
\end{equation}
This relation follows from the identities $E_{b,0}+d_{b}=ad(b)$ and \eqref{tderbracket},
since we have
$ [b,b']+E_{b,0}(b')-E_{b',0}(b)=[b',b]-d_b(b')+d_{b'}(b)=\{b',b\}$.
}

\vspace{.1cm}
\noindent $\bullet$ {Let $\lie{sder}_{2}^{(y)}$ be the Lie subalgebra of 
$\lie{tder}_2^{(y)}$
consisting of derivations $D\in\lie{tder}_2$ such that 
$D(x)=[x,a],  D(y)=0,D(z)=[z, c]$
for some $a,c\in\lie{lie}_2$.}

\vspace{.1cm}
\noindent $\bullet$ Let {$\lie{sder}_{2}^{(z)}$} denote the Lie subalgebra of 
$\lie{tder}_2$ 
consisting of derivations such that
$$E_{a,b}(z)=-[x,a]-[y,b]=0.$$
{We compute
\begin{align*}[E_{a,b},E_{a',b'}](x)&=E_{a,b}([x,a']-E_{a',b'}([x,a])\\
&=[[x,a],a']+[x,E_{a,b}(a')]-[[x,a'],a]-[x,E_{a',b'}(a)]\\
&=[x,E_{a,b}(a')-E_{a',b'}(a)+[a,a']]
\end{align*}
and similarly
$$[E_{a,b},E_{a',b'}](y)=[y,E_{a,b}(b')-E_{a',b'}(b)+[b,b']],$$
so the bracket of two
derivations $E_{a,b}$ and $E_{a',b'}$ is given by $E_{a'',b''}$ with
\begin{equation}\label{Eabbracket}a''=E_{a,b}(a')-E_{a',b'}(a)+[a,a'],\ \ \ b''=E_{a,b}(b')-E_{a',b'}(b)+[b,b'].
\end{equation}}



\vspace{.2cm}
We have the following diagram showing the connections between these subspaces:
\begin{equation}\label{derdiagram}
\xymatrix{
\lie{oder}_2\incl[r]{\ar[dd]_\simeq} &\lie{der}_2\\
&\lie{tder}_2\incl[u]\\
{\lie{sder}_{2}^{(z)}}\incl[ur]\ar[r]^\simeq
& {\lie{sder}_{2}^{(x)}}\incl[u].
}
\end{equation}

\begin{rem}
The correspondence  $D_{b,a}\mapsto E_{-a,b}$
induces a linear space isomorphism
\begin{equation}\label{eq: oder and sder}
i_{o,z}:    \lie{oder}_2\overset{\simeq}{\to}{\lie{sder}_2^{(z)}},
\end{equation}
but this isomorphism is not compatible with the Lie algebra structures.
\end{rem}

Let $\nu$ be the automorphism  defined by
\begin{equation}\label{nu}
	\nu(x)=z=-x-y,\ \ \nu(y)=y.
\end{equation}
{
\begin{lem}\label{isLie} 
Conjugation by $\nu$ induces an isomorphism of
Lie algebras
\begin{align}
\lie{sder}_2^{(z)}\buildrel\sim\over\rightarrow\lie{sder}_{2}^{(x)}\label{conjbynu}\\
E_{a,b}\mapsto d_{\nu(b)}.\notag
\end{align}
\end{lem}
}
\begin{proof}
Recall that $E_{a,b}\in \lie{sder}_2^{(z)}$ maps $x\mapsto [x,a]$ and $y\mapsto [y,b]$,
and $d_{\nu(b)}\in \lie{sder}_{2}^{(x)}$ is the Ihara derivation defined by
$x\mapsto 0$, $y\mapsto [y,\nu(b)]$. 

Let us first show that $d_{\nu(b)}$ is the conjugate of $E_{a,b}$ by 
$\nu$, i.e. $d_{\nu(b)}=\nu\circ E_{a,b}\circ \nu$ (since $\nu$ is
an involution).
It is enough to show they agree on $x$ and $y$, so we compute
$$\nu\circ E_{a,b}\circ\nu(x)=\nu\circ E_{a,b}(z)=0=
d_{\nu(b)}(x)$$
and
$$\nu\circ E_{a,b}\circ\nu(y)=\nu\circ E_{a,b}(y)=
\nu\bigl([y,b]\bigr)=[y,\nu(b)]=d_{\nu(b)}(y).$$
This shows that $\nu\circ E_{a,b}\circ\nu$ is indeed equal to $d_{\nu(b)}$.
To show that $d_{\nu(b)}$ lies in $\lie{sder}_{2}^{(x)}$, we check that
$d_{\nu(b)}(z)$ is a bracket of $z$ with another
element of $\lie{lie}_2$: 
$$d_{\nu(b)}(z) =\nu\circ E_{a,b}\circ\nu(z)=\nu\circ E_{a,b}(x)=
\nu([x,a])=[z,\nu(a)].$$
The same argument goes the other way to show that conjugation by 
$\nu$ maps an element of $\lie{sder}_{2}^{(x)}$ to an element of $\lie{sder}_2^{(z)}$,
which yields the isomorphism (\ref{conjbynu}) as linear spaces. To
see that it is also an isomorphism of Lie algebras, it suffices to note
that conjugation by $\nu$ preserves the Lie bracket of derivations
in $\lie{der}_2$, i.e.
$$\nu\circ [D_1,D_2]\circ \nu=[\nu\circ D_1\circ \nu,\nu\circ D_2\circ\nu],$$
since $\nu$ is an involution.
Since the Lie brackets on $\lie{sder}_2^{(z)}$ and $\lie{sder}_{2}^{(x)}$ are just
restrictions to those subspaces of the Lie bracket on the space of
all derivations, conjugation
by $\nu$ carries one to the other.
\end{proof}

\subsection{Definition of the Kashiwara-Vergne Lie algebra $\lie{krv}$}\label{11}
The free associative algebra ${\rm Ass}_2=\Q\langle\langle x,y\rangle\rangle$ 
on non-commutative generators $x,y$ 
(i.e.~the ring of formal power series in $x$ and $y$) 
is the completion with respect to the degree
of the universal enveloping algebra of the free Lie algebra $\lie{lie}_2$
on $x$ and $y$.

\begin{defn}
(1). The {\it trace linear space} $\lie{tr}_2$ (cf.~\cite{AT}) is defined to be the quotient of
${\rm Ass}_2$ by the equivalence relation given between words in $x$ and $y$ 
by $w\sim w'$ if $w'$ can be obtained from $w$ by a cyclic permutation of 
the letters of the word $w$, and extended linearly to polynomials. 
The natural projection is denoted 
$${tr}:{\rm Ass}_2 \rightarrow \lie{tr}_2.$$
For any polynomial $f\in {\rm Ass}_2$ with constant term $c$, 
we can decompose $f$ in two ways as
\begin{equation}\label{decomps}
f=c+f_xx+f_yy=c+xf^x+yf^y
\end{equation}
for uniquely determined polynomials $f_x,f_y,f^x,f^y$ in $Ass_2$.
        
(2).  The {\it divergence} map is given by
\[{\rm div} :  \begin{array} {rcl}
\lie{tder_2} &\longrightarrow & \lie{tr}_2 \\
u= E_{a,b} &\longmapsto &{\rm tr} (a_xx +b_yy).  \end{array} 
\]
\end{defn}

\begin{defn}\label{def:Kashiwara-Vergne Lie algebra}
{\rm The {\it Kashiwara-Vergne Lie algebra} $\lie{krv}_2$ is defined to be the 
	subspace of $\lie{sder}_2^{(z)}$ of derivations $E_{a,b}$ 
		such that there exists a one-variable power series $h(x)\in\Q[x]$ 
		of degree $\ge 2$ such that
		\begin{equation}\label{defkrv}
	{\rm div}(E_{a,b}) =tr\bigl(h(x+y)-h(x)-h(y)\bigr).
	\end{equation}}
	\end{defn}

	\vspace{.2cm}
	This definition comes from \cite{AT}, where it was shown that $\lie{krv}_2$
	is actually a Lie subalgebra of {$\lie{sder}_2^{(z)}$}. This Lie algebra 
	inherits a weight-grading from that of $\lie{lie}_2$, for which  
	$E_{a,b}$ is of weight $n$ if $b$ (and thus also $a$) is a Lie polynomial
	of homogeneous degree
	$n$.  In particular, the weight 1 part of $\lie{krv}_2$ is spanned by the
	single element $u=E_{y,x}$, and the weight 2 part is zero.  In this
	article, we do not consider the weight 1 part of $\lie{krv}_2$.  For
	convenience, we set $\lie{krv}=\oplus_{n\ge 3} (\lie{krv}_2)_n$,
	where $(\lie{krv}_2)_n$ denotes the weight graded part of $\lie{krv}_2$
	of weight $n$.  We have
	$$\lie{krv}_2=(\lie{krv}_2)_1\oplus \lie{krv}=\Q[E_{y,x}]\oplus \lie{krv}.$$
	Because the other Lie algebras in the literature that are most often 
	compared with the Kashiwara-Vergne Lie algebra have no weight 1 or weight 
	2 parts, it makes most sense to compare them with $\lie{krv}$.  Thus it is
	$\lie{krv}$ that we study for the remainder of this article.

	The Lie algebra $\lie{krv}$ also inherits a depth filtration from the 
	depth grading on $\lie{lie}_2$, for which $E_{a,b}$ is of depth $r$ 
	if $r$ is the smallest number of $y$'s occurring in any monomial of $b$.
	We write $gr\,\lie{krv}$ for the {\it completed} associated graded for this
	depth filtration, so that $gr\,\lie{krv}$ is a Lie algebra that is 
	bigraded for the weight and the depth; we write $gr_n^r\,\lie{krv}$ for the
	part of weight $n$ and depth $r$.  Essentially, an element of
	$gr\,\lie{krv}$ is a derivation {$E_{\bar a,\bar b}\in
	\lie{sder}_2^{(z)}$} where 
	$\bar a,\bar b$ are the lowest-depth parts (i.e.~the parts of
			lowest $y$-degree) of elements $a,b\in\lie{lie}_2$ such that 
	$E_{a,b}\in \lie{krv}$. If $\bar b$ is of homogeneous $y$-degree
	$r$, then $\bar a$ is of homogeneous $y$-degree $r+1$.

\begin{exm}
 The smallest element of $\lie{krv}$ is in weight $3$
	and is given by $E_{a,b}$ with
	$$a=[[x,y],y], \ \ b=[x,[x,y]].$$
	Since $\bar{a}=a$ and $\bar{b}=b$, this is also equal to $E_{\bar{a},\bar{b}}
	\in gr\,\lie{krv}$.
	The next smallest element of $\lie{krv}$ is in weight $5$, and the depth-graded
	part $E_{\bar{a},\bar{b}}$ is given by
	$$\bar{a}=[x,[x,[[x,y],y]]]-2[[x,[x,y]],[x,y]],\ \ \ \ \ 
	\bar{b}=[x,[x,[x,[x,y]]]].$$
    \end{exm}
    By application of mould theory, 
     $\lie{krv}$ is identified with the subspace 
     $ARI(\mathcal F_\ser)_{al+tsen/circconst}$ of moulds given in \eqref{eq: ARIal+tsen/circconst}
    via the isomorphism \eqref{krviso}, \eqref{VtoW} and \eqref{eq:ma Wkrv to ARIal+tsen*circconst}.

\subsection{The Grothendieck-Teichm\"uller Lie algebra $\lie{grt}$ and the double shuffle Lie algebra $\lie{ds}$}\label{12} 
Recall that the 
{\it Grothendieck-Teichm\"uller Lie algebra} $\lie{grt}$ is the space of polynomials $b\in \lie{lie}_2$ satisfying the famous {\it pentagon relation},
equipped with the Poisson bracket (\ref{Poisson}). 
This algebra was first introduced by Y.~Ihara in \cite{I}, with three defining
relations, as a particular derivation algebra of $\lie{lie}_2$ (via
the association $b\mapsto d_b$ as in (\ref{Iharader})); 
it was  subsequently shown that the pentagonal relation implies the
other two in \cite{F1}.  

Recall also that the {\it double shuffle Lie algebra}  $\lie{ds}$ is the 
space of polynomials $b\in \lie{lie}_2$ satisfying a particular set of
conditions on the coefficients called the {\it stuffle relations},
studied in the first place by Racinet (cf.~\cite{R}), who gave a quite
difficult proof that $\lie{ds}$ is also a Lie algebra under
the Poisson bracket (\ref{Poisson}).  This proof was later somewhat streamlined in the appendix of \cite{F2}.  
Another proof with a different approach, identifying the space as a stabilizer was given in \cite{EF}.
Putting together basic elements from \'Ecalle's mould theory,
particularly with the identification of  $\lie{ds}$ with $ARI({\mathcal F}_\ser)_{\underline{al}/\underline{il}}$ (cf. \cite{S2}),
also yields a completely different and very simple proof of this result 
($\!\!$\cite{SS}).

	There is a commutative triangle of injective Lie algebra homomorphisms:
	\begin{equation}\label{triangle}\xymatrix{
		\lie{grt}\ar@{^{(}->}[rr]\ar@{^{(}->}[dr]&& \lie{ds}\ar@{^{(}->}[dl]\\
						&\lie{krv}&.}
						\end{equation}
						The existence of the injection $\lie{grt}\rightarrow\lie{ds}$ was proven
						in \cite{F1}; it is given by $b(x,y)\mapsto b(x,-y)$. The existence of the
						injection $\lie{grt}\rightarrow \lie{krv}$ was proven in \cite{AT};
						it is given by $b(x,y)\mapsto b(z,y)$ where $z=-x-y$. Finally, 
						the existence of the injection from $\lie{ds}$ to $\lie{krv}$ which is given by $b(x,y)\mapsto b(z,-y)$
						was proven in \cite{S1}  under the assumption of 
the senary relation
and which was subsequently validated through independent proofs in \cite{S4} and \cite{EF2}
(see footnote \ref{foot:S4 and EF2}). 
                        We add that
						these morphisms respect the weight gradings and depth filtrations on
						all three spaces.

\subsection{Results on the linearized Kashiwara-Vergne Lie algebra $\lie{lkrv}$}
For $i\ge 1$, set $C_i=ad(x)^{i-1}(y)$ for $i\ge 1$, and let
$\lie{lie}_C$ denote the degree completion of the Lie algebra 
freely generated over $\Q$ by $C_1,C_2,\ldots$.  
By Lazard elimination, $\lie{lie}_C$ is free on the $C_i$ and
		\begin{equation}\label{lieC}
		\lie{lie}_2\simeq \Q\cdot x\oplus \lie{lie}_C.
		\end{equation}
		Thus, Lazard elimination shows that every polynomial $b\in \lie{lie}_2$ 
		having no linear term in $x$ can be written uniquely as a Lie polynomial 
		in the $C_i$.

		\begin{defn}\label{defn:push}
Let the {\it push-operator} be defined on monomials in $x,y$ by
	\begin{equation}\label{pushoperator}
	push(x^{a_0}yx^{a_1}y\cdots yx^{a_r}) = x^{a_r}yx^{a_0}y\cdots yx^{a_{r-1}}.
		\end{equation}
	The push is considered to act
		trivially on constants and powers of $x^n$, so we can extend it to all
		of $Ass_2$ by linearity.
        Any element $b\in Ass_2$ is said to be 

		\vspace{.1cm}
	\noindent $\bullet$ {\it push-invariant} if $push(b)=b$, and 

\vspace{.1cm}
\noindent $\bullet$ {\it push-neutral} if 
$b^r+push(b^r)+\cdots+push^r(b^r)=0$ for
all $r\ge 1$, where $b^r$ denotes the depth $r$ part of $b$.  Finally,
we say that $b$ is 

\vspace{.1cm}
\noindent $\bullet$ {\it circ-neutral} if $b^y$ is push-neutral in
depths $r>1$ when writing $b=c+xb^x+yb^y$ ($c\in\Q$).
\end{defn}

Our $\lie{lkrv}$ is constructed by the push operator.

\begin{defn}\label{pollkrv} 
{ The {\it linearized Kashiwara-Vergne Lie 
algebra} $\lie{lkrv}$ is the space of elements $b\in\lie{lie}_C$ of
degree $\ge 3$ such that

\vspace{.2cm}
\noindent (i) $b$ is push-invariant, and

\vspace{.1cm}
\noindent (ii) $b$ is circ-neutral.}
\end{defn}

\vspace{.3cm}
Our first result is that  $\lie{lkrv}$ forms a bigraded Lie algebra.

\begin{prop}\label{lkrvLie} 
The space $\lie{lkrv}$ is bigraded by weight and depth, and forms a Lie algebra under the Poisson bracket defined in
(\ref{Poisson}).
\end{prop}

It is proved in \S \ref{subsec: circ-neutrality and the second defining relation}.
We note that it is also  proved in \cite{FK} Remark 2.18 and Theorem 2.23 
that 
the depth $>1$-part of $\lie{lkrv}$  forms a Lie algebra.


\vspace{.2cm}
In \S \ref{sec:2}, we show how we derive the definition of $\lie{lkrv}$ 
via a reformulation of the defining properties of 
$\lie{krv}$, in the sense that the defining properties of $\lie{lkrv}$ are merely truncations of the two reformulated defining properties of $\lie{krv}$ 
to their lowest-depth parts. 
	The reformulation 
    suffices
	to prove the following result on $\lie{lkrv}$,
    whose proof is given 
    in the end of \S \ref{subsec:The linearized Kashiwara-Vergne Lie algebra}.  

\begin{prop}\label{firstprop} 
There is an injective Lie algebra morphism
$$gr\,\lie{krv}\hookrightarrow\lie{lkrv}.$$
\end{prop}

We speculate that these two spaces are in fact isomorphic.  

\vspace{.2cm}
In using this type of definition for $\lie{lkrv}$, 
we are following the analogous situation of the well-known double shuffle Lie algebra $\lie{ds}$
and the associated {linearized double shuffle space} $\lie{ls}$.
\footnote{It is shown in \cite{M} that  it is also isomorphic to the linear dual of Goncharov's dihedral Lie coalgebra \cite{G1}.
    } 
    \begin{defn}[$\!\!$\cite{Br}]\label{def: ls}
    The {\it linearized double shuffle Lie algebra} $\lie{ls}$ is the bigraded linear space,
    actually shown to be bigraded Lie algebra under the bracket \eqref{Poisson},
        defined as the set of Lie polynomials $f\in \lie{lie}_2$
	of weight $n\ge 3$ such that the polynomial $f_yy$, rewritten in the variables 
	$y_n=x^{n-1}y$ for $n\ge 1$, is an element of the free Lie algebra on the 
	$y_i$. One also adds the extra assumption that if $f$ is of depth 1, then it 
	is of odd weight. 
    \end{defn}
The above additional assumption  is not needed for $\lie{lkrv}$ as it
	follows from the push-invariance condition in the definition.  By its
	very construction, there is an injective Lie algebra homomorphism
	\begin{equation}\label{grdstols}gr\,\lie{ds}\hookrightarrow \lie{ls},
	\end{equation}
	and it is speculated that these two spaces are isomorphic, but like
	for $\lie{lkrv}$, this is still an open question.

	\vspace{.2cm}
	The injective Lie algebra morphism (\ref{triangle}) from $\lie{ds}$ to
	$\lie{krv}$ yields a corresponding bigraded injective map: 
	\begin{equation}\label{grarrow}
	\xymatrix{gr\,\lie{ds}\,\ar@{^{(}->}[r]&gr\,\lie{krv}}
			\end{equation}
			Our next result shows that there is a Lie algebra map on the generalized spaces
			spaces $\lie{ls}$ and $\lie{lkrv}$ (without any recourse to \'Ecalle's theorem).

\begin{thm}\label{firstthm} 
There is a bigraded Lie algebra injection on linearized spaces	\begin{equation}\label{firstthmeq}
\lie{ls}\hookrightarrow \lie{lkrv}.
\end{equation}
For all $n\ge 3$ and $r=1,2,3$, the map is an isomorphism of the bigraded parts
$$\lie{ls}_n^r\simeq \lie{lkrv}_n^r.$$
\end{thm}
Via the isomorphisms 
\begin{equation}\label{eq:ma ls ARIal/al}
    ma:\lie{ls}\buildrel\sim\over\rightarrow  ARI({\mathcal F}_\ser)_{\underline{al}/\underline{al}}
\end{equation}
given in Lemma \ref{lem:ma ls ARIal/al}
and
$$ma:\lie{lkrv}\buildrel\sim\over\rightarrow ARI({\mathcal F}_\ser)_{al+push/circneut}$$
constructed in \eqref{lkrvARI},
the proof of the theorem reduces to establishing an inclusion
$$ARI({\mathcal F}_\ser)_{\underline{al}/\underline{al}}
\hookrightarrow
ARI({\mathcal F}_\ser)_{al+push/circneut}$$
as formalized in  Theorem \ref{mouldfirstthm}.

\begin{rem}		
If the isomorphisms $\lie{ls}\simeq gr\,\lie{ds}$ and
$\lie{lkrv}\simeq gr\,\lie{krv}$ hold, then the map in Theorem \ref{firstthm}
is the same as the map \eqref{grarrow}. 
Without those speculations, we can only say that 
the Lie algebra injection \eqref{firstthmeq} should extend the map \eqref{grarrow},
fitting into a commutative diagram
\begin{equation}\xymatrix{
gr\,\lie{ds}\,\ar@{^{(}->}[r]\ar@{^{(}->}[d]&gr\,\lie{krv}\ar@{^{(}->}[d]\\
\lie{ls}\ar@{^{(}->}[r]&\lie{lkrv}.}
\end{equation}
\end{rem}

\begin{rem}
We note that in \cite{FK}, the Lie algebra $\lie{lkrv}(\Gamma)$
is constructed for any finite abelian group $\Gamma$,
thereby extending the definition of $\lie{lkrv}$
In particular, when $\Gamma$ is trivial 
this construction coincides with $\lie{lkrv}$  for depths greater than 
$1$ (cf. \cite[Remark 2.18]{FK}).
Moreover, Propositions \ref{lkrvLie} and \ref{firstprop} 
and Theorem \ref{firstthm}
for
$\lie{lkrv}(\Gamma)$ were established in \cite{FK}.
\end{rem}

Adding a variety of known results in the depth 2 and depth 3 situations to this result, we obtain the following corollary.

\begin{cor}\label{firstcor}
The following spaces are isomorphic for $n\ge 3$ and $r=1,2,3$:
$$gr_n^r\lie{grt}\simeq gr_n^r\lie{ds}\simeq gr_n^r\lie{krv}\simeq
\lie{ls}_n^r\simeq \lie{lkrv}_n^r.$$
In particular, all of these spaces are zero when $r=1$ or $3$ and $n$ is
even, or when $r=2$ and $n$ is odd.
\end {cor}

\begin{proof}
The dimensions of the spaces $gr_n^r\lie{grt}$, $gr_n^r\lie{ds}$ and 
$\lie{ls}_n^r$ in depths are known to be equal to each other
in depths $r\le 3$ ($\!\!$\cite{R}, \cite{G}). 
Indeed more is known than merely the dimensions: 
\vskip .2cm
\noindent $\bullet$ the spaces $gr_n^1\lie{grt}$, $gr_n^1\lie{ds}$ and
$\lie{ls}_n^1$ are all $0$ when $n$ is even and 1-dimensional generated by 
$ad(x)^{n-1}(y)$ when $n$ is odd;
\vskip .2cm
\noindent the spaces $gr_n^2\lie{grt}$, $gr_n^2\lie{ds}$ and $\lie{ls}_n^2$ 
are all $0$ when $n$ is odd and spanned by the double Poisson brackets 
$\{ad(x)^{p-1}(y),ad(x)^{q-1}(y)\}$ for odd $p,q\le 3$ with $p+q=n$ when $n$ is even; 
\vskip .2cm
\noindent $\bullet$ the spaces $gr_n^3\lie{grt}$, $gr_n^3\lie{ds}$ and $\lie{ls}_n^3$ 
are all $0$ when $n$ even and spanned by the triple brackets
$\{ad(x)^{p-1}(y),\{ad(x)^{q-1}(y),ad(x)^{s-1}(y)\}\}$ with odd $p,q,s\ge 3$ 
and $p+q+s=n$ when $n$ is odd.
\vskip .2cm
\noindent (Note that the proof for $r=3$ and odd $n$ is much more difficult than the proof 
for $r=2$, and was discovered by Goncharov \cite{G}; as for the case $r\ge 4$, 
the analogous result is known to be false.)  By Theorem \ref{firstthm}, we see that
$\lie{lkrv}_n^r\simeq \lie{ls}_n^r$ for $r=1,2,3$
Finally, since it is known that $\lie{grt}$ injects into $\lie{krv}$ (cf.~\cite{AT}),
we have a corresponding injection $gr_n^r\lie{grt}\hookrightarrow gr_n^r\lie{krv}$,
so by Proposition \ref{firstprop}, shows that $gr_n^r\lie{krv}$ is sandwiched between
$gr_n^r\lie{grt}$ and $\lie{lkrv}_n^r$, which are equal for $r=1,2,3$.  This concludes the proof.  
\end{proof}

We speculate that $\lie{lkrv}_n^r\simeq \lie{ls}_n^r$ for all $n,r$, 
and calculations up to about $n=15$ bear this speculation out, but we were 
not able to prove the isomorphism for any other cases, not even the special
case $n\not\equiv r$ mod 2, where it is well-known that 
$gr_n^r\lie{grt}=gr_n^r\lie{ds}=\lie{ls}_n^r=0$
(cf.~\cite{IKZ}, \cite{Br} for classical proofs, or \cite{S2} for the 
 exposition of \'Ecalle's mould-theoretic proof). 

\subsection{The elliptic Grothendieck-Teichm\"uller Lie algebra $\lie{grt}_{ell}$ and the elliptic double shuffle Lie algebra $\lie{ds}_{ell}$}     
In two independent articles, H.~Tsunogai \cite{Ts} and B.~Enriquez \cite{E1} defined a Lie algebra that Enriquez calls the 
{\it elliptic Grothendieck-Teichm\"uller Lie algebra} 
$\lie{grt}_{ell}$, based on the idea that just as Ihara had defined 
$\lie{grt}$ as the algebra of derivations on $\lie{lie}_2$ (identified with 
		the braid Lie algebra on four strands) that extend to a particular
type of derivation on the braid Lie algebra 
on five strands, $\lie{grt}_{ell}$ is the Lie algebra of derivations on 
$\lie{lie}_2$ (now identified with the genus one braid Lie algebra on two
		strands) that extend to a very particular type of derivation 
of the genus one braid Lie algebra on three strands.
The construction of $\lie{grt}_{ell}$ shows that it is a Lie subalgebra of
$\lie{oder}_2$, and that there is a canonical surjection 
\begin{equation}\label{surj}
s:\lie{grt}_{ell}\rightarrow \lie{grt}.
\end{equation}
Let $\lie{r}_{ell}$ denote the kernel.  Enriquez \cite{E1} showed that there also exists a Lie
algebra morphism 
\begin{equation}\label{sect}
\gamma:\lie{grt}\rightarrow \lie{grt}_{ell}
\end{equation}
that is a section of (\ref{surj}), i.e.~such that 
$\gamma\circ s=id$ on $\lie{grt}$. Thus, there is a
semi-direct product isomorphism
\begin{equation}\label{semidirect}
\lie{grt}_{ell}\simeq \lie{r}_{ell}{\mathbb o}\gamma(\lie{grt}).
\end{equation}

While an elliptic version $\lie{ds}_{ell}$ of the double shuffle 
Lie algebra $\lie{ds}$ was constructed in \cite{S3} using mould theory
as $\Delta\bigl(ARI^\Delta_{\underline{al}*\underline{al}}\bigr)$
 (cf. Definition \ref{defndsell}),
 which induces an inclusion
 \begin{equation}\label{eq: from ls to dsell}
 \lie{ls}\hookrightarrow\lie{ds}_{ell}
 \end{equation}
 by definition and \eqref{eq:ma ls ARIal/al}
    and it is shown there that like $\lie{grt}_{ell}$, 
    $\lie{ds}_{ell}$ is a Lie subalgebra of $\lie{oder}_2$, 
    and that there is an injective Lie algebra homomorphism
    \begin{equation}
        \tilde\gamma:\lie{ds}\rightarrow \lie{ds}_{ell}
    \end{equation}
    that makes the diagram
    \begin{equation}\label{diagram grt and ds}
    \xymatrix{
	    \lie{grt}\ar@{^{(}->}[rr]\ar[d]_\gamma&&\lie{ds}\ar[d]^{\tilde\gamma}\\
			    \lie{grt}_{ell}\ar[dr]&&\lie{ds}_{ell}\ar[dl]\\
			    &\lie{oder}_2&.}
			    \end{equation}
 commutative.
\subsection{Results on the elliptic Kashiwara-Vergne Lie algebra $\lie{krv}_{ell}$ }\label{14}
Our definition of {\it elliptic Kashiwara-Vergne Lie algebra}
is based on that of the linearized Lie algebra $\lie{lkrv}$,
differing only from Definition \ref{pollkrv} by
the denominator appearing in (\ref{Bdenom}), which makes it
impossible to express it directly in terms of Lie elements 
like Definition \ref{pollkrv}.

	\begin{defn}\label{mouldkrvell}
{\rm The {\it elliptic Kashiwara-Vergne linear 
space} $\lie{krv}_{ell}$ is spanned by the elements $b\in\lie{lie}_C$
such that when writing the depth $r$ part 
$B^r(u_1,\dots,u_r)$ of $ma(b)\in ARI(\mathcal F_\ser)$  (see \eqref{eq: ma}})
and setting
			\begin{equation}\label{Bdenom}
		B_*^r(u_1,\ldots,u_r):=\frac{1}{u_1\cdots u_r(u_1+\cdots+u_r)}B^r
			(u_1,\ldots,u_r),
			\end{equation}
		we have

            \vspace{.1cm}
            
	 (i) $B_*^r$ is {\it push-invariant} for $r\ge 1$;
i.e.
	\begin{equation}\label{pushinv}
M^r(u_0,u_1,\ldots,u_{r-1})=M^r(u_1,\ldots,u_r)
	\end{equation}
    holds for $M^r=B^r$
	where $u_0=-u_1-\cdots-u_r$, and
    
(ii) $swap (B_*)^r$ (cf. Definition \ref{swap}) 
        is {\it circ-neutral} up to addition of a constant  $C^r\in\Q$ for $r>1$;
            that is,
	\begin{equation}\label{circneut}
	M^r(v_1,\ldots,v_r)+M^r(v_2,\ldots,v_r,v_1)+\cdots+
	M^r(v_r,v_1,\ldots,v_{r-1})=0
	\end{equation}
        holds for $M^r=swap (B_*)^r+C^r$.
			\end{defn}

Its mould theoretical reformulation is given as 
{$\Delta\bigl(ARI(\mathcal F_\Lau)^{\Delta}_{al+push*circneut}\bigr)$}
in Definition \ref{defnkrvell}.
			\vspace{.2cm}
			The first main result on $\lie{krv}_{ell}$ is of course that it is a 
			bigraded Lie algebra, but this comes from an injective map from 
			$\lie{krv}_{ell}$ into $\lie{oder}_2$ rather than into $\lie{sder}_2^{(z)}$ 
			as for $\lie{lkrv}$.

			\begin{thm}\label{krvellisLie} 
			(i) The space $\lie{krv}_{ell}$ is bigraded for the weight and
			the depth.

			\vspace{.1cm}
			\noindent (ii) For each $b\in\lie{krv}_{ell}$, there exists a unique polynomial
			$a\in \lie{lie}_C$, called the {\rm partner} of $b$, 
			such that $D_{b,a}\in \lie{oder}_2$.

			\vspace{.1cm}
			\noindent (iii) The image of the injective linear map $b\mapsto D_{b,a}$ 
			is a Lie subalgebra of $\lie{oder}_2$; in other words
			$\lie{krv}_{ell}$ is a Lie algebra under the Lie bracket 
			\begin{equation}\label{ellbracket}
			\langle b,b'\rangle=D_{b,a}(b')-D_{b',a'}(b)
			\end{equation}
coming from the bracket of derivations as in (\ref{dabbracket})
	and (\ref{bbracket}). 
	\end{thm}

	\vspace{.1cm}
	This theorem is proven in \S \ref{subsec:Definition of the elliptic Kashiwara-Vergne Lie algebra}.

	\vspace{.3cm}
	The following  is an analogue of the map \eqref{eq: from ls to dsell}
which will be key to the comparison between $\lie{lkrv}$ and 
	$\lie{krv}_{ell}$, and to the proof that $\lie{lkrv}$ is a Lie algebra.

	\begin{prop}\label{liemapprop} 
    There is an injective linear map
	\begin{align*}
	\lie{lkrv}&\hookrightarrow \lie{krv}_{ell}\notag\\
		b(x,y)&\mapsto [x,b(x,[x,y])].
		\end{align*}
		 In fact, this linear map is actually a Lie algebra homomorphism.
		\end{prop}
This claim is proven in \S \ref{subsec:Definition of the elliptic Kashiwara-Vergne Lie algebra} as Corollary \ref{lastone}.

			    \vspace{.2cm}
Our second main result on $\lie{krv}_{ell}$ is an analog of the
existence of $\gamma$ and $\tilde\gamma$ in the diagram \eqref{diagram grt and ds}.

\begin{thm}\label{krvsection} 
Assume 
the isomorphism
in \eqref{Ecstatement}.	
Then we obtain an injective Lie algebra morphism
$$\lie{krv}\hookrightarrow\lie{krv}_{ell}$$
\end{thm}
This theorem will be proven in \S \ref{subsec:the map from krv to krvell}.

\vspace{.2cm}
 Based on the known injective Lie algebra homomorphisms
$\xymatrix{\lie{grt}\ar@{^{(}->}[r]&\lie{ds}\ar@{^{(}->}[r]&\lie{krv}}$ 
evoked in  \S \ref{12} above,
we believe that there are corresponding injective Lie algebra homomorphisms
between the elliptic versions of these Lie algebras.
However, we were not able to prove that $\lie{grt}_{ell}$ as defined
in \cite{E1} injects into $\lie{ds}_{ell}$ or $\lie{krv}_{ell}$. 
To circumvent this difficulty,
we define a Lie subalgebra $\widetilde{\lie{grt}}_{ell}\subset \lie{grt}_{ell}$,
conjecturally isomorphic to $\lie{grt}_{ell}$, as follows.

\begin{defn}
For $n\ge 0$, let $\delta_{2n}\in \lie{oder}_2$ 
denote the derivation of $\lie{lie}_2$ defined by 
 $$\delta_{2n}(x)=ad(x)^{2n}(y),\ \ \delta_{2n}([x,y])=0.$$
Let $\lie{b}$ be the Lie subalgebra of $\lie{oder}_2$ generated
by the $\delta_{2n}$.
\end{defn}

 Enriquez showed in \cite{E1} that $\delta_{2n}\in \lie{r}_{ell}$ for $n\ge 0$,
 so $\lie{b}$ is a Lie subalgebra of $\lie{r}_{ell}$.
 Let $\lie{B}$ denote the Lie ideal generated by $\lie{b}\subset
 \lie{r}_{ell}$ under the semi-direct action of $\gamma(\lie{grt})$ on 
 $\lie{r}_{ell}$ of (\ref{semidirect}).  
 We set 
 \begin{equation}\label{grtelltilde}
\widetilde{\lie{grt}}_{ell}=\lie{B}{\Bbb o}\gamma(\lie{grt}).
\end{equation}

\vspace{.2cm}
Our third main result on $\lie{krv}_{ell}$ relates all these maps
via a commutative diagram.

\vspace{.1cm}
 \begin{thm}\label{bigdiagram} We have the following commutative diagram of injective Lie algebra homomorphisms:
 \begin{equation}\xymatrix{
 \lie{grt}\ar@{^{(}->}[r]\ar@{^{(}->}[d]&\lie{ds}\ar@{^{(}->}[r]\ar@{^{(}->}[d]&\lie{krv}\ar@{^{(}->}[d]\\
\widetilde{\lie{grt}}_{ell}\ar@{^{(}->}[r]\ar@{^{(}->}[rd]&\lie{ds}_{ell}\ar@{^{(}->}[r]\ar@{^{(}->}[d]&
\lie{krv}_{ell}\ar@{^{(}->}[dl]\\
&\lie{oder}_2.&}
\end{equation}
\end{thm}
The theorem will be proven in \S \ref{doubleshuf}.

\section{Reformulation of the definition of $\lie{krv}$ and definition of the linearized Lie algebra $\lie{lkrv}$}\label{sec:2}
 In this section, we give a convenient reformulation of the defining conditions 
 of $\lie{krv}$ given in Definition \ref{def:Kashiwara-Vergne Lie algebra}, which leads to a simple definition of the linearized version 
$\lie{lkrv}$ that passes easily into the language of moulds which 
will be essential for our subsequent proofs in \S\ref{sec: mould theory} and \S \ref{sec:The elliptic Kashiwara-Vergne Lie algebra}. 

\subsection{The first defining condition of $\lie{krv}$: specialness}
 The first of the two defining conditions of $\lie{krv}$ is that 
$\lie{krv}$ lies in $\lie{sder}_2^{(z)}$, i.e.~elements of $\lie{krv}$ are
special tangential derivations having the form
$E_{a,b}$ with $E_{a,b}(x)=[x,a]$, $E_{a,b}(y)=[y,b]$ and
$[x,a] + [y,b]=0$.

 The following equivalent formulations of the property of specialness as properties of the polynomial $b$ were given in \cite{S1}.

\begin{prop}[{$\!\!\!$\cite[Theorem 2.1]{S1}}]
\label{specialprop} 
Let $b \in \lie{lie}_C$ be homogeneous with degree $n \geq 1$;
write 
$$b=b_xx+ b_yy=xb^x+yb^y.$$
Then the following are equivalent:

\vspace{.2cm}
 (i) There exists a unique element $a\in \lie{lie}_C$ such that $[x,a]+[y,b]=0$;

\vspace{.1cm}
(ii) $b$ is push-invariant; 

\vspace{.1cm}
(iii) $b_y=b^y$.
\end{prop}

\begin{proof}
The claim when $n\geq 3$ can be proven as follows:
    The equivalent between (i) and (ii) is given in \cite[Theorem 2.1]{S1},
    where it is also shown that (ii) is equivalent to the antipalindrome for $b_y$.
    By \cite[(2.1)]{S1}, we see that it is equivalent to (iii).
    The case for $n=1,2$ can be proved by a direct computation.
\end{proof}

Thanks to this proposition, we can now reformulate the first defining 
condition of $\lie{krv}$ as follows: the pair of polynomials $a,b\in\lie{lie}_C$
satisfies $[x,a]+[y,b]=0$ if and only if $b$ is push-invariant and 
$a$ is its partner.  
\vspace{.2cm}

\subsection{The second defining condition of $\lie{krv}$: divergence}
We now consider the second defining condition of $\lie{krv}$, the divergence
condition.  Because $\lie{krv}$ is weight-graded, we may restrict attention 
to derivations $E_{a,b}$ of homogeneous weight $n$, i.e. such that
$a$ and $b$ are Lie polynomials of homogeneous degree $n\ge 3$. The
second defining condition (\ref{defkrv}) then simplifies to the existence of
 a constant $c\in\Q$ such that
\[{\rm tr}(xa_x +yb_y)=c\ {\rm tr}
 \bigl((x+y)^n-x^n-y^n\bigr)\ {\rm in}\ \lie{tr}_2. \]

  Let us reformulate this as a condition only on $b$, just as we did for the first
 defining condition.  Since $a \in \lie{lie}_2$, its trace is zero and thus
${\rm tr}(xa_x)={\rm tr}(a_xx)=-{\rm tr}(a_yy)=-{\rm tr}(ya_y)$, so
 \[{\rm tr}(xa_x +yb_y) = {\rm tr}(yb_y- ya_y).  \]
 Since $E_{a,b} \in \lie{sder}^{(z)}$, we have $ [x,a] = [b,y]$. Expanding 
 this in terms of the decompositions of $a$ and $b$, we obtain 
 \[ xa_xx + xa_yy -xa^xx - ya^yx = xb^xy+ yb^yy- yb_xx -yb_yy,\]
 from which we deduce that $a_y = b^x$ and $a^y=b_x$.
Thus 
 \[{\rm tr}(yb_y-ya_y)={\rm tr}(yb_y- yb^x)= {\rm tr}\bigl(y(b_y- b^x)\bigr).\]
 From Proposition \ref{specialprop},  we have $b_y=b^y$, so now, using
 the circularity of the trace, 
 the divergence condition can be reformulated as
 \[{\rm tr}\bigl((b^y-b^x)y\bigr) =c\ {\rm tr}\bigl((x+y)^n-x^n+y^n\bigr).  \]
 We use this to express it as a condition directly on $b^y-b^x$ as follows, 
 using the push-operator defined in (\ref{pushoperator}).

 \vspace{.2cm}
\begin{defn}\label{polypushconst}
{\rm A polynomial $b\in Ass_2$ of homogeneous weight $n>1$ is 
	said to be {\it push-constant for the value $c\in\Q$} if $b$ does not contain 
		the monomial $y^n$ and for each 
        {$0<r<n$}, writing $b^r$ for the depth $r$
		part of $b$, we have
		\[\sum_{i=0}^r push^i(b^r)=c\ \sum_w w\]
		where the sum in the right-hand factor is over all monomials of weight $n$ and 
		depth $r$ (for $push$, see \eqref{pushoperator}).
		Equivalently, $b$ is push-constant  for the value $c$ if it does not contain $y^n$ and 
		for all monomials $w\ne x^n$ {with depth $r$ and weight $n$}, we have 
			\[ \sum_{v \in Push(w)} (b | v) =c \]
				where $(b|v)$ denotes the coefficient of the monomial $v$ in $b$, and
				$Push(w)$ is the list (with possible repetitions) 
				$[w,push(w),\ldots,push^r(w)]$. 
                {The monomial $x^n$ is said to be push-constant for its own coefficient. If $b$ is push-constant for the value $c=(b|x^n)$, then we say that $b$ is {\it strictly push-constant.
				If $c=0$, then $b$ is said to be {\it push-neutral} (cf. Definition \ref{defn:push}).  } }}
				\end{defn}

\begin{exm}\label{eg:poshconstant}
The simplest example of a push-constant polynomial is obtained by separating the
full set of words of given degree $n$ and
depth $d$ into push-orbits, taking one 
representative $w$ from each push-orbit, and
summing over the lists $Push(w)$. For example
when $n=5$ and $d=2$, the full list of words
is given by the 10 words
$$\{x^3y^2,x^2yxy, x^2y^2x, xyx^2y,xyxyx,
xy^2x^2,yx^3y,yx^2yx,yxyx^2,y^2x^3\},$$
the push-orbits are given by
$$\{x^3y^2,yx^3y,y^2x^3\},\ 
\{x^2yxy,yx^2yx,xy^2x^2\},\
\{x^2y^2x,yxyx^2,xyx^2y\},\
\{xyxyx\},$$
the $Push$-lists are the same as the orbits except in the case of $xyxyx$ where the 
$Push(xyxyx)=[xyxyx,xyxyx,xyxyx]$ three times,
and the polynomial $b$ is thus given by
$$b=x^3y^2+x^2yxy+x^2y^2x+xyx^2y+3xyxyx+xy^2x^2+yx^3y+yx^2yx+yxyx^2+y^2x^3,$$
which is push-constant for the value $c=3$. Note that since $(b|x^5)=0$, this polynomial is not strictly push-constant.

Interesting strictly push-constant polynomials can be obtained from elements 
$\psi\in\lie{grt}$ by writing $\psi=x\psi^x+y\psi^y$ and taking $b=\psi^y$.  In this way we obtain for example the degree 4 polynomial $b$ which is strictly push-constant for the value $c=(b|x^4)=1$: 
\begin{align*}
b=x^4+3x^3y&-\frac{9}{2}x^2yx-\frac{1}{2}x^2y^2+\frac{11}{2}xyx^2+2xyxy+\frac{9}{2}xy^2x+4xy^3-2yx^3\\ &-\frac{1}{2}yx^2y-\frac{11}{2}yxyx-6yxy^2+2y^2x^2+4y^2xy-y^3x.
\end{align*}
\end{exm}
The following proposition shows that the divergence condition
comes down to requiring that $b^y-b^x$ be push-constant.

\begin{prop}[$\!\!$\cite{S1}]\label{circprop} 
	Let $b\in Ass_2$ be a push-invariant Lie polynomial of homogeneous degree $n$.
	Then $b$ satisfies the divergence condition
	$${\rm tr}((b^y-b^x)y) =c\ {\rm tr}\bigl((x+y)^n-x^n-y^n\bigr)$$
	if and only if $b^y-b^x$  is push-constant for the value $nc$. Furthermore,
	if this is the case then 
	\begin{equation}\label{constval}c=\frac{1}{n}(b\, |\, x^{n-1}y).
	\end{equation}
	\end{prop}

The following proof originates from the arguments given in   \cite[Section 3]{S1}.  

\begin{proof}
 Let $w$ be a monomial of degree $n$ and depth $r\ge 1$, 
	and let $C_w$ denote 
	the list of words obtained from $w$ by cyclically permuting the letters,
	so that $C_w$ contains exactly $n$ words (with possible repetitions). 
	Let $C_w^y$ denote the list obtained from $C_w$ by removing all words ending
	in $x$, so that $C_w^y$ contains exactly $r$ words.
	Write $C_w^y = [u_1y,\ldots, u_ry]$. Then we have the equality of lists 
	\[ [u_1,...,u_r] = Push(u_1). \]

	Let $c_w={\rm tr}(w)$, i.e.~$c_w$ is the equivalence class of $w$, which is
	the set of the words in the list $C_w$, without repetitions: thus
	$C_w$ is nothing other than $n/|c_w|$ copies of $c_w$. The divergence condition
	\[{\rm tr}\bigl((b^y- b^x)y\bigr) = c\ {\rm tr}\bigl((x+y)^n-x^n-y^n\bigr)\]
	translates as the following family of conditions for 
	one word in each equivalence class $c_w$:
	\begin{equation}\label{temp}
	\sum_{v \in c_w} \bigl((b^y-b^x)y \,|\, v\bigr) = c|c_w|,
	\end{equation}
	where each side is the coefficient of the class $c_w$ in the trace,
	i.e.~the sum of the coefficients of the words in $c_w$ in the original
	polynomial.

	If $r>1$, we can choose a word $uy\in C_w$ that starts in $y$. 
	Then from (\ref{temp}), the divergence condition on $b$ implies that
	\begin{align*}c&=\frac{1}{|c_w|}\sum_{v\in c_w} \bigl((b^y-b^x)y\,|\,v\bigr)\\
		&=\frac{1}{n} \sum_{v \in C_w} \bigl((b^y-b^x)y \,|\, v\bigr)\\
		&=\frac{1}{n}\sum_{v \in C_w^y} \bigl((b^y-b^x)y \,|\, v\bigr)\\
		&=\frac{1}{n}\sum_{u' \in Push(u)} ((b^y-b^x) \,|\, u').
		\end{align*}
		This is exactly the definition of $b^y-b^x$ being push-constant for the 
		value $nc$. 

		If $r=1$, then $w$ is of depth $1$, $|c_w|=n$ and $x^{n-1}y$ is the only word 
		in $c_w$ ending in $y$. Thus (\ref{temp}) comes down to 
		$$\bigl((b^y-b^x)y\,|\,x^{n-1}y\bigr)=nc.$$
		But since $b$ is a Lie polynomial, we have $(b|x^n)=(b^x|x^{n-1})=0$,
		so using $b^y=b_y$ (by Proposition \ref{specialprop}), we also have
$$\bigl((b^y-b^x)y\,|\,x^{n-1}y\bigr)=(b^y-b^x\,|\,x^{n-1})
	=(b^y|x^{n-1})$$ $$=(b_y|x^{n-1})=(b_yy|x^{n-1}y)=(b|x^{n-1}y),$$
	which proves that $nc=(b|x^{n-1}y)$ as desired.  Note that this condition
	means that if $b$ has no depth 1 part, then $b^y-b^x$ is push-neutral.
\end{proof}
	\vspace{.3cm}
	We now have a new way of expressing $\lie{krv}$, which is much easier to translate 
	into the mould language. 

	\begin{defn}[{cf. \cite[Theorem 1.2]{S1}}]\label{Vkrv} 
Let  $V_{\lie{krv}}$ be the completion of
the linear space spanned 
by polynomials $b\in\lie{lie}_C$ of homogeneous degree $n\ge 3$ such that

\vspace{.2cm}
(i)\ $b$ is push-invariant, and

\vspace{.1cm}
(ii)\ $b^y-b^x$ is push-constant for the value $(b\,|\,x^{n-1}y)$,

\vspace{.2cm}
\noindent
equipped with the Lie bracket
$$\{b,b'\}=[b,b']+E_{a,b}(b')-E_{a',b'}(b)$$
where $a$ and $a'$ are the (unique) partners of $b$ and $b'$ respectively.
\end{defn}

			Indeed, since Propositions \ref{specialprop} and \ref{circprop} show that 
			\begin{align}\label{krviso}
{\lie{krv}}\ \ &\buildrel\sim\over \rightarrow \ V_{\lie{krv}}\notag\\
		E_{a,b}\ &\mapsto \ \ b
		\end{align}
is an isomorphism of linear spaces (shown in \cite[Theorem 1.2]{S1})
and $\lie{krv}$ is known to be a
		Lie subalgebra of $\lie{sder}_2^{(z)}$, the bracket on $V_{\lie{krv}}$ is
		inherited directly from this and makes $V_{\lie{krv}}$ into a Lie algebra (cf. \cite[Theorem 1.2]{S1}).

\subsection{The linearized Kashiwara-Vergne Lie algebra $\lie{lkrv}$}\label{subsec:The linearized Kashiwara-Vergne Lie algebra}
		Using the above isomorphism of $\lie{krv}$ with the linear space
		$V_{\lie{krv}}$ given by $E_{a,b}\mapsto b$, let us now consider the
		depth-graded versions of the defining conditions of $V_{\lie{krv}}$,
		i.e.~determine what these conditions say about the lowest-depth parts of
		elements $b\in V_{\lie{krv}}$.  The push-invariance is a depth-graded condition,
		so it restricts to the statement that the lowest depth part of $b$ is still
		push-invariant; in particular, by Proposition \ref{specialprop} it admits of a 
		unique partner $a\in \lie{lie}_C$ such that $[x,a]+[y,b]=0$, i.e. such that 
		the associated derivation $E_{a,b}$ lies in $\lie{sder}_2^{(z)}$.

		In the second condition, if $b$ is of degree $n$ and depth $r=1$ and 
		$b^1$ denotes the lowest-depth part of $b$, then $(b^1)^y=x^{n-1}$,
		so the push-constance condition on $b^1$ is empty since 
		$(b^1)^y=(b|x^{n-1}y)x^{n-1}$. If $r>1$, however, then $(b|x^{n-1}y)=0$ and so
		the push-constance condition on $b^y-b^x$ is actually push-neutrality,
		which implies the push-neutrality of $(b^r)^y$ alone, since $(b^r)^y$ is
		the only part of the expression $b^y-b^x$ of minimal depth $r-1$. 
        
These observations lead directly to the definition of the 
{\it linearized version} $\lie{lkrv}$ {\it of the Kashiwara-Vergne Lie algebra} 
given in Definition \ref{pollkrv} above, and that by definition it is bigraded
by weight and depth. The statement of Proposition \ref{lkrvLie}, that $\lie{lkrv}$ is a Lie 
algebra under the bracket coming from the bracket of derivations in
$\lie{sder}_2^{(z)}$, 
namely
\begin{equation}\label{sderzbracket}
\{b,b'\}=[b,b']+E_{a,b}(b')-E_{a',b'}(b), 
\end{equation}
{will be proved at the end of \S \ref{subsec:Definition of the elliptic Kashiwara-Vergne Lie algebra}.}  
Up to this fact, the proof
of Proposition \ref{firstprop} now follows trivially from the equivalences above.
\vskip .3cm
{\it Proof of Proposition \ref{firstprop} assuming Proposition \ref{lkrvLie}.\footnote{{Proposition \ref{lkrvLie} will be proven in 
Proposition \ref{prop: ARI al push circneut Lie algebra}.
}}} 
The defining
properties of the associated graded $gr\,\lie{krv}$ are properties
satisfied by the the lowest-depth parts of elements of $\lie{krv}$. 
We identify $\lie{krv}$ with $V_{\lie{krv}}$ and use the version of 
its defining properties expressed in Definition \ref{Vkrv}.
Since both the properties of being a Lie element and being push-invariant
respect the depth, the same properties are satisfied by elements 
of $gr\,\lie{krv}$. For the divergence, the argument in the paragraph
preceding this proof shows that it implies no condition on the lowest-depth
part if the depth is 1, and it implies the push-neutrality of the
lowest-depth part if the depth is $>1$. We do not know if this property
along with being Lie and push-invariant, which together define $\lie{lkrv}$,
are all that is implied on the lowest-depth part of an element of $\lie{krv}$ by
its defining properties, but we certainly know that they all hold for the lowest-depth part, and therefore we obtain the desired inclusion of linear spaces
\begin{equation*}
gr\,\lie{krv}\hookrightarrow \lie{lkrv}.
\end{equation*}
We now show that this inclusion is a Lie algebra homomorphism. The space $gr(\lie{krv})$ is equipped with the Lie bracket inherited from the Lie bracket on $\lie{krv}$, which is itself inherited from $\lie{sder}^{(z)}$ since $\lie{krv}$ is a Lie subalgebra of $\lie{sder}^{(z)}$; this is the bracket given explicitly in \eqref{sderzbracket}. Now suppose that $B,B'\in \lie{krv}$ (more precisely, the derivations $E_{A,B}$ and
$E_{A',B'}$ are in $\lie{krv}$, where $A$ and $A'$ denote the uniquely determined partners of $B$ and $B'$ such that $E_{A,B}(z)=E_{A',B'}(z)=0$). Let $b,b'$ be the lowest-depth part of $B$ and $B'$, of depths $r$ and $s$ respectively, and let $a$ and $a'$  denote the lowest-depth parts of $A$ and $A'$, which are also the partners of $b$ and $b'$, in the sense that $E_{a,b}(z)=E_{a',b'}(z)=0$. Since $E_{a,b}(x)=[x,a]$ and $E_{a,b}(y)=[y,b]$, this means that $[x,a]+[y,b]=0$, so $a$ must be of depth $r+1$, and $a'$ must be of depth $s+1$. 

Now, the lowest depth part of $[E_{A,B},E_{A',B'}]$ necessarily arises from restricting the bracket to the lowest-depth parts of $A,B,A',B'$, i.e.~from the expression
\begin{equation}\label{bracket} 
[b,b']+E_{a,b}(b')-E_{a',b'}(b).
\end{equation}
The term $[b,b']$ is of homogeneous depth $r+s$. However, the term $E_{a,b}(b')$, is obtained by replacing each of the $s$ $y$'s of each monomial $b'$ by $[y,b]$, giving terms of depth $s+r$, and each $x$ of $b'$ by $[x,a]$, giving terms of depth $s+r+1$ since $a$ is of depth $r+1$. Similarly, $E_{a',b'}(b)$ contains both terms of depth $s+r$, from replacing each $y$ in $b$ by $[y,b']$, and terms of depth $s+r+1$, from replacing each $x$ in $b$ by $[x,a']$. Thus, to get the lowest depth part of \eqref{bracket}, we must ignore the terms coming from replacing $x$ in $b'$ by $[x,a]$ and replacing $x$ in $b$ by $[x,a']$. This comes down to applying the derivation $d_b$ to $b'$ (the Ihara derivation mapping $x\mapsto 0$ and $y\mapsto [y,b]$), and similarly applying $d_{b'}$ to $b$. So the lowest-depth part of \eqref{bracket} is actually given by
\begin{equation}\label{newbracket}
[b,b']+d_b(b')-d_{b'}(b),
\end{equation}
which is precisely the bracket on $\lie{sder}^{(x)}$, and corresponds in mould terms to the ari bracket. Since $\lie{lkv}$ is a Lie algebra equipped with this bracket by Proposition \ref{lkrvLie}, this shows that the linear map $gr(\lie{krv})\rightarrow\lie{lkrv}$ is actually an injective Lie algebra morphism.
 \ \qed
\begin{rem}
 No examples of elements of $\lie{lkrv}$ that are
      not truncations to lowest depth of elements of $\lie{krv}$ are known. 
      It would be interesting to try to prove the equality of $\lie{lkrv}$ 
      with $gr\,\lie{krv}$ by starting with a polynomial $\lie{lkrv}$ of depth $r>1$ 
      and finding 
      a way to construct a depth by depth lifting to an element of $\lie{krv}$.
\end{rem}

\section{Mould theory}\label{sec: mould theory}
In this section, we recall the language of \'Ecalle's moulds  theory and reformulate
the defining conditions of $\lie{lkrv}$ in this language.  We end the
section with the proof of Theorem \ref{firstthm} 
and its corollary in terms of moulds.  We hope that this section and
the next one, which explores the elliptic version of $\lie{krv}$,
will illustrate the way in which moulds are powerful tools in this context.

\subsection{Moulds and various operators}
We recall the definitions of moulds and  various operators  which are employed in this article.

Let ${\mathcal F}=\cup_{m\geqslant 0}{\mathcal  F}_m$ denote a family of functions (see \cite{FHK} for a precise definition).
In this article, unless otherwise specified, all statements are formulated for a general family of functions $\mathcal{F}$. However, for practical purposes, it suffices to keep the following specific cases in mind:
\begin{itemize}
\item ${\mathcal F}_\pol=\cup_{m\geqslant 0}{\mathcal  F}_{\pol, m}$
with ${\mathcal F}_{\pol, m}=\Q[u_1,\dots,u_m]$,
\item ${\mathcal F}_\rat=\cup_{m\geqslant 0}{\mathcal  F}_{\rat, m}$
with ${\mathcal F}_{\rat, m}=\Q(u_1,\dots,u_m)$,
\item ${\mathcal F}_\ser=\cup_{m\geqslant 0}{\mathcal  F}_{\ser, m}$
with ${\mathcal F}_{\ser, m}=\Q[[u_1,\dots,u_m]]$,
\item ${\mathcal F}_\Lau=\cup_{m\geqslant 0}{\mathcal  F}_{\Lau, m}$
with ${\mathcal F}_{\Lau, m}=\Q((u_1,\dots,u_m))$ which is the quotient field of ${\mathcal F}_{\ser, m}$.
\end{itemize}

The notion of moulds was introduced by J. \'{E}calle in \cite{Ec81}.
However  this article  adopts the formulation in \cite{S1}.

\begin{defn}
A {\it mould} valued in a family of function 
${\mathcal F}=\cup_{r\geqslant 0}{\mathcal  F}_r$ is a collection 
$A=\bigl(A^r(u_1,\ldots,u_r)\bigr)_{r\ge 0}$ where each  component
$A^r(u_1, u_2, ..., u_r)$ belongs to ${\mathcal F}_r$.
The set $\mathcal M(\mathcal F)$ of all moulds valued in $\mathcal F$ 
forms a $\Q$-linear space under componentwise addition and scalar multiplication.
We denote by $ARI(\mathcal F)$ the linear subspace consisting of all moulds 
with   $A^0=0\in\mathcal F_0$.
We say that a mould $A$ is \textit{concentrated in depth} $r$ if $A^s=0$ 
for all $s\ne r$, and we let $ARI({\mathcal F})^r\subset ARI({\mathcal F})$ be the subspace of moulds concentrated in depth $r$. 
Thus $ARI({\mathcal F})=
\prod_{r\ge 1} ARI({\mathcal F})^r$.
\end{defn}



For later use, we list several operations on $\mathcal M(\mathcal F)$ below:

\begin{defn}[{$\!\!$\cite[(2.4)--(2.11)]{Ec}}]
\label{def: various mould operators}
For any mould $ M=( M^m(u_1,\dots, u_m))_m$ in $\mathcal{M}(\mathcal{F})$,
the following $\Q$-linear operators  are defined componentwise
\footnote{
In \cite{Ec}, the operation $circ$ is denoted by ${pus}$. 
However, to avoid confusion with $push$, we adopt the notation $circ$ throughout this paper.
}: 
\begin{align*}
& pari(M)(u_1,\ldots,u_m)=(-1)^m M(u_1,\ldots,u_m), \\
&mantar( M)^m(u_1,\dots,u_m)=(-1)^{m-1} M^m(u_m,\dots,u_1),\\
&push( M)^m(u_1,\dots,u_m)= M^m(-u_1-\cdots-u_m, u_1,\dots,u_{m-1}),\\
&circ( M)^m(u_1,\dots,u_m)= M^m(u_m,u_1,\dots,u_{m-1}),\\
&neg( M)^m(u_1,\dots,u_m)= M^m(-u_1,\dots,-u_m),\\
&teru( M)^m(u_1,\dots,u_m)= M^m(u_1,\dots,u_m) \\
&\qquad+\frac{1}{u_m}\left\{ M^{m-1}(u_1,\dots,u_{m-2},u_{m-1}+u_m)- M^{m-1}(u_1,\dots,u_{m-2},u_{m-1})\right\}.
\end{align*}
\end{defn}

For convenience, we also introduce $\overline{\mathcal M}(\mathcal F)$ (resp. $\overline{ARI}(\mathcal{F})$), which is a copy of $\mathcal M(\mathcal F)$ (resp. $ARI(\mathcal{F})$). 
However, to distinguish between these two, we use the parameters $v_1, v_2, \ldots$ instead of $u_1, u_2, \ldots$ in $\overline{\mathcal M}(\mathcal F)$ and $\overline{ARI}(\mathcal{F})$.
We now introduce \'Ecalle's important {\it swap} operator on moulds.

\begin{defn}\label{swap}
The operator $swap: \mathcal M(\mathcal F)\to \overline{\mathcal M}(\mathcal F)$ is defined by
$$swap(B)(v_1,\ldots,v_r)=B(v_r,v_{r-1}-v_r,\ldots,v_1-v_2)$$
for $B\in \mathcal M(\mathcal F)$. 
The inverse operator mapping $\overline{\mathcal M}(\mathcal F)$ to $\mathcal M(\mathcal F)$ 
(which we also denote by $swap$, as the context is clear according to whether 
$swap$ is acting on a mould in $\mathcal M(\mathcal F)$ or one in $\overline{\mathcal M}(\mathcal F)$) is given by
$$swap(C)(u_1,\ldots,u_r)=C(u_1+\cdots+u_r,u_1+\cdots +u_{r-1},\ldots,u_1)$$
for $C\in \overline{\mathcal M}(\mathcal F)$.
Thus it makes sense to write $swap\circ swap=id$.
\end{defn}

							\vspace{.1cm}
							We also need to consider an important symmetry on moulds, based on the
{\it shuffle} operator on tuples of commutative variables, which is defined by
$$Sh\bigl((u_1,\ldots,u_i)(u_{i+1},\ldots,u_r)\bigr)=
\Bigl\{(u_{\sigma^{-1}(1)},\ldots,u_{\sigma^{-1}(r)})\,|\,\sigma\in S_r^i
\Bigr\},$$
where $S_r^i$ is the subset of permutations $\sigma\in S_r$ such that
$\sigma(1)<\cdots<\sigma(i)$ and $\sigma(i+1)<\cdots<\sigma(r)$.

\vspace{.2cm}
\begin{defn}
A mould $A\in ARI({\mathcal F})$ is {\it alternal} if in each depth $r\ge 2$ we have
\begin{equation*}
	\sum_{w \in Sh((u_1, ..., u_i)(u_{i+1}, ..., u_r))} A^r(w) = 0
	\ \ \ {\rm for\ \ } 1\leq i \leq \left[\frac {r}{2} \right].
\end{equation*}
By convention, the alternality condition is void in depth 1, i.e. all depth
1 moulds are considered to be alternal.
\end{defn}

\begin{exm}
    In depth 4, there are two alternality conditions, given by
		\begin{flalign*}
		&A(u_1, u_2, u_3,u_4) +   A(u_2, u_1, u_3,u_4) +  A(u_2, u_3, u_1,u_4)+  A(u_2, u_3, u_4,u_1) = 0 \\
	& A(u_1, u_2, u_3,u_4) +   A(u_3, u_1, u_2,u_4) +  A(u_3, u_4, u_1,u_2)
	+  A(u_1, u_3, u_2,u_4)\\ &\ \  +  A(u_1, u_3, u_4,u_2) +  A(u_3, u_1, u_4,u_2)= 0 
\end{flalign*}
\end{exm}

\begin{defn} \label{def: ARI al/al}
(1).
We write $ARI({\mathcal F})_{al}$ for the linear subspace of $ARI({\mathcal F})$ consisting of alternal 
moulds.  

(2).
Let $ARI(\mathcal F)_{al/al}$ denote the linear space of moulds
that are alternal and have alternal swap and
following \'Ecalle's notation \cite{Ec},
let $ARI(\mathcal F)_{\underline{al}/\underline{al}}$ denote the linear subspace of 
$ARI(\mathcal F)_{al/al}$ of moulds that are even in depth 1. 

(3).
Let $ARI(\mathcal F)_{al\ast al}$ denote the linear space of moulds
that are alternal and whose swap are alternal up to addition of constant-valued moulds.
Following \cite{S2},
we denote by  $ARI(\mathcal F)_{\underline{al}\ast\underline{al}}$
the linear subspace of  $ARI(\mathcal F)_{{al}\ast{al}}$
moulds that are even functions in depth 1.
\footnote{We use the notation $ARI_{P_1/P_2}$ for the linear space of moulds 
having property $P_1$ and whose swaps have property $P_2$.
We also use a slightly more general notation $ARI_{P_1*P_2}$ to 
denote the space of moulds having property $P_1$ and 
whose swap has property $P_2$
up to adding on a constant-valued mould.\label{foot: conventionARIsymbol}}

\end{defn}

\begin{exm}
An example of a mould in $ARI_{al*al}$
is the mould 
where $A$ is the 
concentrated 
$\{u_1u_2u_3(u_1+u_2+u_3)\}^{-1}A(u_1,u_2,u_3)$
in depth 3 with
$$
A(u_1,u_2,u_3)=-{{1}\over{4}}u_1^3u_2+{{1}\over{4}}u_1^3u_3-
{{1}\over{4}}u_1^2u_2^2+{{1}\over{2}}u_1^2u_3^2+{{1}\over{4}}u_1u_3^3
-{{1}\over{4}}u_2^2u_3^2-{{1}\over{4}}u_2u_3^3$$
$$-{{1}\over{12}}u_1^2u_2u_3
+{{1}\over{6}}u_1u_2^2u_3-{{1}\over{12}}u_1u_2u_3^2.$$
It is easy to check that $\Delta^{-1}(A)$ is alternal, but its swap 
is not alternal unless one adds on the constant $\frac{1}{3}$. 
\end{exm}

\subsection{The map ma}

Let ${ Ass}_2$ be
the degree-completed free associative algebra on $x,y$.
It forms a (topological) noncommutative co-commutative Hopf algebra
whose coproduct is given by $\Delta(x)=x\otimes 1+1\otimes x$ and
$\Delta(y)=y\otimes 1+1\otimes y$.
Let $e:{ Ass}_2\to \Q$ be the map taking the coefficient of $x$.
Define
\begin{equation}\label{eq: Ass dag 2}
{Ass}^\dag_2:=\{w\in { Ass}_2\bigm| (e\otimes id)\circ \Delta(w)=0 \}.
\end{equation}

\begin{lem}
(i).
The subspace $Ass_2^\dag$ constitutes  a Hopf subalgebra of $Ass_2$.

(ii).
Let  $Ass_C=\Q\langle\langle C\rangle\rangle=\Q\langle\langle C_1,C_2,\ldots\rangle\rangle$ 
be the  degree-completion of the free
non-commutative polynomial algebra on variables $C_i = {\rm ad}_x^{i-1}(y)$ 
for $i\ge 1$, for the degree given by deg$\,C_i=i$.
Then we have a natural identification of Hopf algebras:
\begin{equation}\label{eq: Ass dag 2 simeq Ass C}
    {Ass}^\dag_2\simeq Ass_C.
\end{equation}

(iii).
The algebra $Ass_2$ admits the following decomposition as completed tensor products:
$$Ass_2 = \mathbb{Q}[[x]] \hat\otimes 
Ass_2^\dag
$$
\end{lem}

\begin{proof}
The above claims follow from \cite[Lemma 46]{FHK} and \cite[Proposition 47]{FHK}.
\end{proof}

\begin{rem}\label{rem:other basis of Q[C]}
We consider $C_{a_1,\dots,a_{r-1},a_r}\in Ass_2$ defined by
$$
C_{a_1,\dots,a_{r-1},a_r}:=\left\{
\begin{array}{ll}
{\rm ad}_x(C_{a_1,\dots,a_{r-1},a_r-1})& (a_r\ge2), \\
C_{a_1,\dots,a_{r-1}}\cdot y& (a_r=1),
\end{array}
\right.
$$
for $r\ge1$ and $a_1,\dots,a_r\ge1$.
By definition, 
the set $\{C_{a_1,\dots,a_r}\ |\ r\ge1, a_1,\dots,a_r\ge1\}\cup \{1\}$  
forms a linear basis of $Ass_C$.
\end{rem}

From now on, we regard $Ass_C$ as a Hopf subalgebra of $Ass_2$.
Let 
\begin{equation}\label{eq: lie C}
    \lie{lie}_C=\lie{lie}_2\cap Ass_C
\end{equation}
where we regard $Ass_C$ as a subspace of $Ass_2$ under the identification
\eqref{eq: Ass dag 2 simeq Ass C}.
By \cite[Proposition 47]{FHK}, it is identified with
the depth $\geq 1$-part of $\lie{lie}_2$.
This identification induces the decomposition  \eqref{lieC} and 
through which we establish that
this is  free  generated by  $C_i$ with $i\geq 1$
by Lazard elimination (cf.\cite{Bbk}).

\begin{lem}\label{Liestuff} 

(i). For $r,n\ge 1$, let $Ass_C^{ (r,n)}$ denote the (finite-dimensional)
	subspace of $Ass_C$ spanned by elements corresponding to monomials $C_{a_1}\cdots C_{a_r}$ with
	$a_1+\cdots+a_r=n$ under the identification \eqref{eq: Ass dag 2 simeq Ass C}, and let
   ${\mathcal M}(\mathcal F_\ser)^{(r,n)}$
    denote the subspace of ${\mathcal M}(\mathcal F_\ser)$
    \footnote{We may choose instead ${\mathcal M}(\mathcal F_\pol)$.}
	consisting of polynomial moulds of degree $n-r$ concentrated in depth $r$.
	The map
	\begin{equation}\label{eq: ma}
    ma:\begin{array}{rl} 
	Ass_C^{(r,n)}& \longrightarrow  {\mathcal M}(\mathcal F_\ser)^{(r,n)}\\
		C_{a_1}\cdots C_{a_r}& \longmapsto (-1)^{n-r} u_1^{a_1-1}\cdots u_r^{a_r-1}	
		\end{array}
        \end{equation}
		is a linear space isomorphism.

(ii). Let 
$\lie{lie}_C^{(r,n)}=\lie{lie}_2\cap Ass_C^{ (r,n)}.$
For each $r\ge 1$, the map $ma$ restricts to a (finite-dimensional)
	linear space isomorphism
	$$ma:\lie{lie}_C^{(r,n)}\rightarrow ARI(\mathcal F_\ser)^{(r,n)}_{al}.$$
    with $ARI(\mathcal F_\ser)^{(r,n)}_{al}:=
    ARI(\mathcal F_\ser)_{al}\cap {\mathcal M}(\mathcal F_\ser)^{(r,n)}$.
	\end{lem}

The above maps yield the following isomorphisms of linear spaces
\begin{align*}
    ma:& Ass_C\overset{\sim}{\to}{\mathcal M}(\mathcal F_\ser), \\
    ma:& \lie{lie}_C\overset{\sim}{\to}ARI(\mathcal F_\ser)_{al}. 
\end{align*}

\begin{exm}
The mould $ma(C_3)=ma([x,[x,y]])$ is the mould concentrated in 
	depth $1$ given by $u_1^2$.  Similarly, $ma(C_2C_1-C_1C_2)=ma([[x,y],y])$
	is the mould concentrated in depth 2 given by $(-1)^{3-2}(u_1^1u_2^0 - u_1^0u_2^1)
	=u_2-u_1$.
\end{exm}

\begin{lem}\label{lem:bijection mi}
Denote $mi:=swap\circ ma$.
Then the map $mi:Ass_C \rightarrow \overline{\mathcal M}(\mathcal F_\ser)$ is bijective, and its restriction ($a_1,\dots,a_r\ge1$ with $a_1+\cdots+a_r=n$)
\begin{equation}\label{eq: mi}
mi:\begin{array}{rl} 
    Ass_C^{(r,n)}& \longrightarrow  \overline{\mathcal M}(\mathcal F_\ser)^{(r,n)}\\
    C_{a_1,a_2,\dots,a_r}& \longmapsto (-1)^{n-r} v_r^{a_1-1}v_{r-1}^{a_2-1}\cdots v_1^{a_r-1}
\end{array}
\end{equation}
is a linear isomorphism.
\end{lem}

\begin{proof}
We prove \eqref{eq: mi} by induction on $r\ge1$.
When $r=1$, the claim is obvious.
Assume that we have \eqref{eq: mi} for $r\le m$, that is, we have
$$
ma(C_{a_1,a_2,\dots,a_r})(u_1,\dots,u_s)
=\delta_{r,s}(-1)^{n-r}u_1^{a_1-1}(u_1+u_2)^{a_2-1}\cdots (u_1+\cdots+ u_r)^{a_r-1}
$$
for $r\le m$ and $s\ge0$.
When $r=m+1$, by Leibnitz rule, we have
\begin{align*}
C_{a_1,\dots,a_m,a_{m+1}}
&=\sum_{i=0}^{a_{m+1}-1} \binom{a_{m+1}-1}{i}
{\rm ad}_x^{i}(C_{a_1,\dots,a_m}){\rm ad}_x^{a_{m+1}-1-i}(y) \\
&=\sum_{i=0}^{a_{m+1}-1} \binom{a_{m+1}-1}{i}
C_{a_1,\dots,a_{m-1},a_{m+i}}C_{a_{m+1}-i}.
\end{align*}
Applying the map $ma$ to both sides, we get
\begin{align*}
&ma(C_{a_1,\dots,a_m,a_{m+1}})(u_1,\dots,u_s) \\
&=\sum_{i=0}^{a_{m+1}-1} \binom{a_{m+1}-1}{i}
ma(C_{a_1,\dots,a_{m-1},a_{m+i}})\times ma(C_{a_{m+1}-i})(u_1,\dots,u_s) \\
&=\delta_{m+1,s}\sum_{i=0}^{a_{m+1}-1} \binom{a_{m+1}-1}{i}
ma(C_{a_1,\dots,a_{m-1},a_{m+i}})(u_1,\dots,u_m) \cdot ma(C_{a_{m+1}-i})(u_{m+1}) \\
&=\delta_{m+1,s}\sum_{i=0}^{a_{m+1}-1} \binom{a_{m+1}-1}{i}
(-1)^{a_1+\cdots+a_m+i-m}u_1^{a_1-1}\cdots (u_1+\cdots+ u_{m-1})^{a_{m-1}-1} \\
&\hspace{3cm}\cdot (u_1+\cdots+ u_{m})^{a_{m}+i-1} \cdot (-1)^{a_{m+1}-i-1} u_{m+1}^{a_{m+1}-i-1} \\
&=\delta_{m+1,s}(-1)^{a_1+\cdots+a_{m+1}-(m+1)}u_1^{a_1-1}(u_1+u_2)^{a_{2}-1}\cdots (u_1+\cdots+ u_{m})^{a_{m}-1} \\
&\hspace{3cm} \cdot \sum_{i=0}^{a_{m+1}-1} \binom{a_{m+1}-1}{i}
(u_1+\cdots+ u_{m})^{i} u_{m+1}^{a_{m+1}-i-1} \\
&=\delta_{m+1,s}(-1)^{a_1+\cdots+a_{m+1}-(m+1)}u_1^{a_1-1}(u_1+u_2)^{a_{2}-1}\cdots (u_1+\cdots+ u_{m+1})^{a_{m+1}-1}.
\end{align*}
Applying the map $swap$ to both sides, we have
\begin{align*}
&mi(C_{a_1,\dots,a_m,a_{m+1}})(v_1,\dots,v_s) \\
&=swap\circ ma(C_{a_1,\dots,a_m,a_{m+1}})(v_1,\dots,v_s) \\
&=\delta_{m+1,s}(-1)^{a_1+\cdots+a_{m+1}-(m+1)}
v_{m+1}^{a_1-1}v_{m}^{a_{2}-1}\cdots v_1^{a_{m+1}-1}.
\end{align*}
So we obtain \eqref{eq: mi} for $r=m+1$.
Hence, we finish the proof.
\end{proof}
		%

		\begin{defn}\label{backwards} {\rm Let $\beta$ denote the {\it backwards writing operator}
			on words in $x,y$, meaning that $\beta(m)$ is obtained from a word
				$m$ by writing it from right to left.  The operator $\beta$ extends to
				polynomials by linearity.}
				\end{defn}

\begin{rem}\label{rem:explicit formulation of ma}
By \eqref{eq: Ass dag 2 simeq Ass C}, the map $ma$ in \eqref{eq: ma} can be regarded as 
$$
ma(b)(u_1,\dots,u_r)
=\sum_{\underline{a}=(0,a_1,\dots,a_r)}k_{\underline{a}}
u_1^{a_1}(u_1+u_2)^{a_2}\cdots (u_1+\cdots+u_r)^{a_r}
$$
for $r\ge0$ and $b\in Ass_2^\dag$ with
$$
b=\sum_{r\ge0}\sum_{\underline{a}=(a_0,\dots,a_r)}k_{\underline{a}}x^{a_0}yx^{a_1}\cdots yx^{a_r}.
$$
So by Definition \ref{swap}, we have
$$
swap(ma(b))(v_1,\dots,v_r)
=\sum_{\underline{a}=(0,a_1,\dots,a_r)}k_{\underline{a}}
v_r^{a_1}v_{r-1}^{a_2}\cdots v_1^{a_r}.
$$
\end{rem}


\subsection{Ari bracket}
We explain the Lie bracket $ari$ on $ARI(\mathcal F)$ and explain how it is related with  the bracket $\{.,.\}$ of \eqref{Poisson}.

We begin by introducing the standard
mould multiplication that \'Ecalle denotes $mu$.

\begin{defn}
 Let $\mathcal F$ be a family of functions.
 For moulds $A,B$ in $\mathcal M(\mathcal F)$, we define the multiplication
$$mu(A,B)(u_1,\ldots,u_r)=\sum_{i=0}^r A(u_1,\ldots,u_i)B(u_{i+1},\ldots,
u_r).$$
By $mu$, $\mathcal M(\mathcal F)$ forms a $\Q$-algebra. 
The associated Lie bracket $lu$ in $ARI(\mathcal F)$
is defined by 
$$lu(A,B)=mu(A,B)-mu(B,A).$$
We write $ARI(\mathcal F)_{lu}$ for $ARI(\mathcal F)$ viewed as a Lie algebra for the $lu$-bracket.
\end{defn}

The identical formulas yield a multiplication and Lie algebra (also called
$mu$ and $lu$) on $\overline{ARI}(\mathcal F)$.

\begin{rem}
If $f$ and $g$ are power series in $Ass_C$ and $A=ma(f)$, $B=ma(g)$
in $\mathcal M(\mathcal F_\ser)$, then
$mu$ is a mould translation of the usual non-commutative multiplication,
and $lu$ the usual Lie bracket:
$$mu(A,B)=ma(fg),\ \ \ \ lu(A,B)=ma\bigl([f,g]\bigr).$$
\end{rem}

In order to define \'Ecalle's $ari$-bracket, 
we first introduce three
derivations of $ARI(\mathcal F)_{lu}$ associated to a given mould $A\in ARI(\mathcal F)$. 
It
is non-trivial to prove that these operators are actually derivations
(cf.~\cite[Prop. 2.2.1]{S2}).

\begin{defn}[{$\!\!$\cite[\S 2.2]{Ec}}]
\label{def: amit and anit}
Let $\mathcal F$ be a family of functions.
Let $B\in ARI(\mathcal F)$.  Then the derivation $amit(B)$ of
$ARI(\mathcal F)_{lu}$ is given by
$$\bigl(amit(B)\cdot A\bigr)(u_1,\ldots,u_r)=
\sum_{0\le i<j<r} A(u_1,\ldots,u_i,u_{i+1}+\cdots+u_{j+1},u_{j+2},\ldots,u_r)
B(u_{i+1},\ldots,u_j),$$
and the derivation $anit(B)$ is given by
$$\bigl(anit(B)\cdot A\bigr)(u_1,\ldots,u_r)=
\sum_{0<i<j\le r} A(u_1,\ldots,u_{i-1},u_{i}+\cdots+u_j,u_{j+1},\ldots,u_r)
B(u_{i+1},\ldots,u_j).$$

\noindent We also have corresponding derivations $\overline{amit}(B)$ and
$\overline{anit}(B)$ of $\overline{ARI}(\mathcal F)_{lu}$ 
for $B\in\overline{ARI}(\mathcal F)$, given by the formulas 
$$\bigl(\overline{amit}(B)\cdot A\bigr)(v_1,\ldots,v_r)=
\sum_{0\le i<j<r} A(v_1,\ldots,v_i,v_{j+1},\ldots,v_r)B(v_{i+1}-v_{j+1},
\ldots,v_j-v_{j+1}),$$
$$\bigl(\overline{anit}(B)\cdot A\bigr)(v_1,\ldots,v_r)=
\sum_{0<i<j\le r} A(v_1,\ldots,v_i,v_{j+1},\ldots,v_r)B(v_{i+1}-v_i,
\ldots,v_j-v_i).$$
Finally, \'Ecalle defines the derivation $arit(B)$ on $ARI(\mathcal F)_{lu}$ by
$$arit(B)=amit(B)-anit(B),$$
and the $ari$-bracket on $ARI(\mathcal F)$ by
\begin{equation}\label{aribrackdef}
ari(A,B)=arit(B)\cdot A-arit(A)\cdot B+lu(A,B),
\end{equation}
as well as the derivation $\overline{arit}$ on $\overline{ARI}(\mathcal F)_{lu}$
and the bracket $\overline{ari}$ on $\overline{ARI}(\mathcal F)$ by the same
formulas with overlines.
\end{defn}

\begin{prop}\label{prop: al closed under ari}
    Let $\mathcal F$ be a family of functions. Then we have the following:
    
    (1). The linear space $ARI(\mathcal F)$ forms a Lie algebra under the ari-bracket \eqref{aribrackdef}.
    
    (2). The linear subspace $ARI(\mathcal F)_{al}$ forms a Lie subalgebra of  $ARI(\mathcal F)$.
    
    (3). The linear subspace $ARI(\mathcal F)_{al/al}$ forms a Lie subalgebra of $ARI(\mathcal F)$.
    
    (4). The linear subspace $ARI(\mathcal F)_{\underline{al}\ast\underline{al}}$ forms a Lie subalgebra of  $ARI(\mathcal F)$.

\end{prop}

\begin{proof}
    The proof of claim (1) is given in \cite{S2}, although the proof of the key formula (2.2.10) looks unprovided therein. A complete proof, together with an extended version of the claim, is given in \cite[Proposition 1.14]{FK}.
    The claim (2)  (resp. (3)) is first proven in \cite[Appendix A]{SS}
    (resp. \cite[Theorem 2.5.6]{S2}), while its generalization is systematically developed in \cite[Proposition 1.15 (resp. Proposition 1.24)]{FK}.
    The claim (4) can be deduced from the claim (3).
\end{proof}

\begin{rem}
The definitions of $amit$, $anit$, $arit$ and $ari$ are generalizations
to all moulds of familiar derivations of $Ass_C$.  Indeed,
if $f,g\in Ass_C$ and $A=ma(f)$, $B=ma(g)$ in ${\mathcal M}(\mathcal F_ \ser)$, then
$$amit(B)\cdot A=ma\bigl(D^l_g(f)\bigr)$$
where $D^l_g$ is defined by $x\mapsto 0$, $y\mapsto gy$,
$$anit(B)\cdot A=ma\bigl(D^r_{g}(f)\bigr)$$
where $D^r_{g}$ is defined by $x\mapsto 0$, $y\mapsto yg$, and thus
$$arit(B)\cdot A=ma\bigl(-d_{g}(f)\bigr)$$
where $d_g$ is the Ihara derivation $x\mapsto 0$, $y\mapsto [y,g]$
(see (\ref{Iharader})), and
\begin{equation}\label{ariPoisson}
ari(A,B)=ma\bigl([f,g]+d_f(g)-d_{g}(f)\bigr)=ma\bigl(\{f,g\}\bigr).
\end{equation}
corresponds to the Ihara or Poisson Lie bracket \eqref{Poisson}
on $\lie{lie}_C$.  (See \cite{S2}, Corollary 3.3.4). 
\end{rem}

\begin{lem}\label{lem: iy}
The map $b \mapsto -E_{b,0}$ induces a Lie algebra isomorphism
    \begin{equation}\label{eq iy}
        i_y:(\lie{lie}_C,\{,\})\overset{\sim}{\to} \lie{tder}_2^{(y)}.
    \end{equation}
\end{lem}

\begin{proof}
Bijectivity is clear. By \eqref{eq: opposite Poisson}, we have 
\[
[-E_{b,0}, -E_{b',0}] = -E_{\{b, b'\}, 0},
\]
which shows that $i_y$ is a Lie algebra homomorphism.
\end{proof}

Denote $Ass_{C,\,\ge1}$ to be the subspace of $Ass_C$ which is the depth $\ge1$-part of $Ass_C$, that is,
$$
Ass_C=\Q1 \oplus Ass_{C,\,\ge1}.
$$
Note that the space $Ass_{C,\,\ge1}$ is generated by 
$\{ C_{a_1}\cdots C_{a_r} \ |\ r\ge1, a_1,\dots,a_r\ge1 \}$.
\begin{lem}\label{lem: ma Lie alg isom}
The map $ma$ induces the Lie algebra isomorphism:
 \begin{align*}
      ma: &(Ass_{C,\,\ge1},\{,\})\to (ARI(\mathcal F_\ser), ari), \\
     ma: & (\lie{lie}_C,\{,\})\to (ARI(\mathcal F_\ser)_{al}, ari).
  \end{align*}  
\end{lem}

\begin{proof}
This follows from the above arguments and Lemma \ref{Liestuff}.    
\end{proof}




\begin{defn}
Let $\mathcal F$ be a family of functions.
We denote by $GARI(\mathcal F)$ the set of all moulds $A$ in $\mathcal M(\mathcal F)$
with constant term 1, that is, $A^0=1\in\mathcal F_0$.
Similarly, we introduce $\overline{GARI}(\mathcal F)$
which is a copy of $GARI(\mathcal F)$.
\end{defn}     

To introduce an operator $\overline{ganit}$, 
we prepare the following notations:
Set ${\bf v}=(v_1,\ldots,v_r)$, and let $W_{\bf v}$ denote the set of 
 decompositions $d_{\bf v}$ of ${\bf v}$ into chunks
\begin{equation}\label{chunks}
	 d_{\bf v}={\bf a}_1{\bf b}_1\cdots {\bf a}_s{\bf b}_s
\end{equation}
for $s\ge 1$, where with the possible exception of ${\bf b}_s$, the
${\bf a}_i$ and ${\bf b}_i$ are non-empty.  
Thus for instance, when $r=2$ there are two decompositions in
$W_{\bf v}$, namely ${\bf a}_1=(v_1,v_2)$ and ${\bf a}_1{\bf b}_1=(v_1)(v_2)$, and when $r=3$ there are four decompositions,
three for $s=1$: ${\bf a}_1=(v_1,v_2,v_3)$, ${\bf a}_1{\bf b}_1=
(v_1,v_2)(v_3)$, ${\bf a}_1{\bf b}_1=(v_1)(v_2,v_3)$, and one for 
$s=2$: ${\bf a}_1{\bf b}_1{\bf a}_2=(v_1)(v_2)(v_3)$.

\begin{defn}[{$\!\!$\cite[\S 2.2]{Ec}}]
For any mould $Q\in \overline{GARI}(\mathcal F)$, 
we define an operation $\overline{ganit}(Q)$ 
given by
	\begin{equation}\label{ganitQ}
	\bigl(\overline{ganit}(Q)\cdot P\bigr)({\bf v})=\sum_{{\bf a}_1{\bf b}_1\cdots
		{\bf a}_s{\bf b}_s\in W_{\bf v}} Q(\lfloor {\bf b}_1)\cdots Q(\lfloor 
				{\bf b}_s) \ P({\bf a}_1\cdots {\bf a}_s),
		\end{equation}
		where if ${\bf b}_i$ is the chunk $(v_k,v_{k+1},\ldots,v_{k+l})$, then we use
		the notation
		\begin{equation}\label{flex}
		\lfloor {\bf b}_i=(v_k-v_{k-1},v_{k+1}-v_{k-1},\ldots,v_{k+l}-v_{k-1}).
		\end{equation}
\end{defn}

We note that in \cite[Theorem 3.7]{K} it is shown that
$\overline{ganit}(Q)$  is an automorphism 
of the Lie algebra $\overline{ARI}(\mathcal F)_{lu}$.

\subsection{The special mould $pal$ and \'Ecalle's fundamental identity}\label{adariinvpal}
We are now ready to introduce the fundamental identity of 
\'Ecalle in Remark \ref{rem: Ecalle fundamental identity}, 
which is the key to 
the proof of Theorem \ref{krvsection} given in \S \ref{subsec:the map from krv to krvell}.


\begin{defn} \label{defn:pal and pil}
{\rm Let constants $c_r\in \Q$, $r\ge 1$, be defined by
setting $f(x)=1-e^{-x}$ and expanding $f_*(x)=\sum_{r\ge 1} c_rx^{r+1}$,
where $f_*(x)$ is the {\it infinitesimal generator} of $f(x)$, defined by 
$$f(x)=\Bigl(exp\bigl(f_*(x)\frac{d}{dx}\bigr)\Bigr)\cdot x.$$
Let 
$$lopil\in\overline{ARI}({\mathcal F}_\Lau)_{\overline{ari}}$$
be the mould defined by 
$lopil(v_1)=-{{1}\over{2v_1}}$ and for $r\ge 2$ by the 
simple expression
\begin{equation}\label{deflopil}
lopil(v_1,\ldots,v_r)=c_r\,\frac{v_1+\cdots+v_r}{v_1(v_1-v_2)\cdots (v_{r-1}-v_r)v_r}
\end{equation}
Set $$pil=exp_{\overline{ari}}(lopil)\in GARI({\mathcal F}_\Lau)$$ where $exp_{\overline{ari}}$ denotes the
exponential map associated to the pre-Lie law $\overline{preari}$ given by 
$$\overline{preari}(A,B)=\overline{arit}(B)\cdot A+mu(A,B)\ \ \ \ \ {\rm on}\quad \overline{ARI}(\mathcal F),$$
and set $$pal=swap(pil)\in GARI({\mathcal F}_\Lau).$$}
\end{defn}

\vspace{.1cm}
The mould $lopil$ is easily seen to be both alternal and circ-neutral
(see Definition \ref{circconstance} and Example \ref{example: lopil is circ-neutral}). 
It is also known (although surprisingly difficult to show) that the mould 
$$lopal=log_{ari}(pal)\in ARI(\mathcal F_\Lau)$$
is alternal (cf.~\cite[\S 4]{Ec2}, or \cite[Chap. 4]{S2}).  
Thus the moulds $pil$ and $pal$ are both exponentials of alternal
moulds; this is called being {\it symmetral}. 
The inverses of $pal$ 
in $GARI(\mathcal F_\Lau)$ and $pil$ in $\overline{GARI}(\mathcal F_\Lau)$
are given by
$$invpal=exp_{ari}(-lopal),\ \ \ \ invpil=exp_{\overline{ari}}(-lopil).$$

\begin{rem}
The key maps we will be using in our proof are the adjoint operators associated to $pal$ and $pil$, given by
\begin{equation}\label{adaripil}
Ad_{ari}(pal) = exp\bigl(ad_{ari}(lopal)\bigr),\ \ \ \ 
Ad_{\overline{ari}}(pil) = exp\bigl(ad_{\overline{ari}}(lopil)\bigr),
\end{equation}
where $ad_{ari}(P)\cdot Q=ari(P,Q)$.
The inverses of these adjoint actions are given by
\begin{equation}\label{adariinvpil}
Ad_{ari}(invpal) = exp\bigl(ad_{ari}(-lopal)\bigr),\ \
Ad_{\overline{ari}}(invpil) = exp\bigl(ad_{\overline{ari}}(-lopil)\bigr).
\end{equation}
These adjoint actions produce remarkable transformations of certain
mould properties into others, and form the heart of much of \'Ecalle's theory of multizeta values. 
\end{rem}

To discuss the fundamental identity, we prepare the following two moulds.

\begin{defn}
Let $pic\in \overline{GARI}(\mathcal F_\Lau)$ be the mould 
defined by 
$$pic(v_1,\ldots,v_r)=\frac{1}{v_1\cdots v_r}.$$
Similarly
let $poc\in\overline{GARI}(\mathcal F_\Lau)$ be the mould
defined by 
\begin{equation}\label{poc}
poc(v_1,\ldots,v_r)=\frac{-1}{v_1(v_1-v_2)\cdots (v_{r-1}-v_r)}.
\end{equation}
\end{defn}

\begin{rem}\label{rem: Ecalle fundamental identity}
\'Ecalle's {\it fundamental identity}
\footnote{
This identity, indicated in a personal communication by \'{E}calle, follows from the first fundamental identity given in (2.62) of \cite{Ec} (see also (2.8.4) of \cite{S2}). 
The complete statement and proof is given in \cite{S2} Theorem 4.5.2;
the proof relies among other things on a basic fact of mould theory stated by 
\'Ecalle and used constantly in the mould literature, namely that the operator 
$\overline{ganit}(pic)$ transforms alternal moulds in $\overline{ARI}(\mathcal F)$
to alternil moulds. 
A full proof of this fact was not written down 
until the recent article \cite{K} by N.~Komiyama, see Corollary 3.25.}
relates the two adjoint actions of (\ref{adaripil}). Valid for all 
push-invariant moulds $R$, it is given by
\begin{equation}\label{Ecallefund}
swap\cdot Ad_{ari}(pal)\cdot R=\overline{ganit}(pic)\cdot Ad_{\overline{ari}}(pil)\cdot swap(R).
\end{equation}
Since it is shown 
\begin{equation}\label{eq: ganit pic and poc}
    \overline{ganit}(pic)^{-1}=\overline{ganit}(poc)
\end{equation}
as automorphisms of $\overline{ARI}(\mathcal F_\Lau)_{lu}$
in the proof of \cite[Lemma 4.37]{B},
we can rewrite the above identity (\ref{Ecallefund}) as
\begin{equation*}
\overline{ganit}(poc)\cdot swap\cdot Ad_{ari}(pal)\cdot R=Ad_{\overline{ari}}(pil)\cdot swap(R).
\end{equation*}
Letting $N=Ad_{ari}(pal)\cdot R$, i.e. $R=Ad_{ari}(invpal)\cdot N$, we further rewrite it in terms of $N$ as
\begin{equation}\label{goodfund}
Ad_{\overline{ari}}(invpil)\cdot \overline{ganit}(poc)\cdot swap(N)=swap\cdot Ad_{ari}(invpal)\cdot N,
\end{equation}
which is valid whenever $R=Ad_{ari}(invpal)\cdot N$ is 
push-invariant.
This identity will be used later in the proof of Proposition \ref{keypropb}.
\end{rem}

\subsection{Push-invariance and the first defining relation of $\lie{lkrv}$}
We  show how to reformulate the first defining property of elements
of $\lie{lkrv}$ in terms of moulds.

\begin{defn}
A mould $B\in ARI(\mathcal F)$ is said to be {\it push-invariant} if 
$$push\,(B) = B.$$ 
\end{defn}

We denote  the subset of push-invariant moulds in $ARI(\mathcal F)$
by $ARI(\mathcal F)_{push}$.

\begin{prop}
\label{prop: pushclosed}
    Let $\mathcal F$ be a family of functions.
    The subspace  $ARI(\mathcal F)_{push}$
    forms a Lie subalgebra  of $ARI(\mathcal F)$
     under the ari-bracket. 
\end{prop}

\begin{proof}
 A complete demonstration of this result, along with its generalization to broader cases, is presented in \cite[Proposition 1.28]{FK}.
\end{proof}



Then the following proposition shows that this notion is precisely the 
translation into mould terms of the property of push-invariance for 
a Lie polynomial given in Definition \ref{defn:push} above.

\begin{prop}\label{pushinvprop} 
Define
$$
Ass_{C,\,\ge1}^{push}:=\{b\in Ass_{C,\,\ge1}\bigm| push(b)=b\}.
$$
Then the map $ma$ induces an isomorphism of linear spaces
\begin{equation*}
    ma:Ass_{C,\,\ge1}^{push}\overset{\sim}{\to}ARI(\mathcal F_\ser)_{push}.
\end{equation*}
Whence $Ass_C^{push}$ forms a Lie algebra under the bracket $\{,\}$ and
we have a Lie algebra  isomorphism 
\begin{equation*}
    ma:(Ass_{C,\,\ge1}^{push},\{,\})\overset{\sim}{\to}
    (ARI(\mathcal F_\ser)_{push},ari).
\end{equation*}
\end{prop}

\begin{proof}
By Lemma \ref{lem: ma Lie alg isom}, it is enough  to show that 
$b\in Ass_{C,\,\ge1}$ is  push-invariant  if and only if $ma(b)$ is a push-invariant mould, i.e. $push(ma(b))=ma(b)$.

Let $b\in Ass_{C,\,\ge1}$ with depth $1$.
We only prove the case $b=C_{a+1}$ for $a\ge0$.
Then we have
\begin{align*}
push(b)
&=\sum_{k=0}^{a}(-1)^k\binom{a}{k}push(x^{a-k}yx^k)
=\sum_{k=0}^a(-1)^k\binom{a}{k}x^kyx^{a-k} \\
&=(-1)^a\sum_{k=0}^a(-1)^k\binom{a}{k}x^{a-k}yx^k
=(-1)^ab.
\end{align*}
So $b$ is push-invariant if and only if $a$ is even.
When $b=C_{a+1}$, by \eqref{eq: ma}, the associated mould $ma(b)$ is given as
$$
ma(b)(u_1,\dots,u_r)
=\left\{\begin{array}{ll}
(-1)^au_1^a& (r=1), \\
0& (r\neq1),
\end{array}\right.
$$
so $a$ is even if and only if $ma(b)$ is push-invariant, that is, we obtain the claim for depth $1$.

				Now, let $b \in Ass_{C,\,\ge1}$ with depth $r\ge2$, 
				and put $f=yb \in Ass_{C,\,\ge1}$. 
                The associated moulds $ma(f)$ and $ma(b)$
				are related by the formula 
				\begin{align}\label{fandb}
				ma(f)(u_0, u_1,...,u_r) &= ma(C_1b)(u_0,\dots,u_r)= ma(C_1)(u_0) ma(b)(u_1,...,u_r) \\
                &= ma(b)(u_1,...,u_r). \nonumber
				\end{align}
                We write 
				$$b= \sum_{\underline{a}=(a_0,\ldots,a_r)} k_{\underline{a}}\, 
				x^{a_0}yx^{a_1}\cdots yx^{a_r},$$ 
                then we have
                $$
                f=x^0yb= \sum_{\underline{a}=(a_0,\ldots,a_r)} k_{\underline{a}}\, 
				x^0yx^{a_0}yx^{a_1}\cdots yx^{a_r}.
                $$
                By Remark \ref{rem:explicit formulation of ma}, we get
                \begin{equation}\label{eqn:explicit formula of swap(ma(f))}
                swap(ma(f))(v_0,\dots,v_r)
                = \sum_{\underline{a}=(a_0,\ldots,a_r)} k_{\underline{a}}
                v_r^{a_0}v_{r-1}^{a_1}\cdots v_1^{a_{r-1}}v_0^{a_r}.    
                \end{equation}
                Because we have
                $$
                push(b)
                =\sum_{\underline{a}=(a_0,\ldots,a_r)} k_{\underline{a}}\, 
				x^{a_r}yx^{a_0}\cdots yx^{a_{r-1}}
                =\sum_{\underline{a}=(a_0,\ldots,a_r)} k_{a_1,\dots,a_r,a_0}\, 
				x^{a_0}yx^{a_1}\cdots yx^{a_{r}},
                $$
                $b=push(b)$ if and only if
				$k_{\underline{a}}=k_{(a_1,\ldots,a_{r},a_0)}$ for each $\underline{a}$,
				this is equivalent to 
				\begin{align}\label{fcirc}
&swap\bigl(ma(f)\bigr)(v_0,\dots,v_r) 
=\sum_{\underline{a}=(a_0,\ldots,a_r)} k_{\underline{a}}
                v_r^{a_0}v_{r-1}^{a_1}\cdots v_1^{a_{r-1}}v_0^{a_r} \\
&= \sum_{\underline{a}=(a_0,\ldots,a_r)} k_{(a_1,\ldots,a_{r},a_0)}
                v_r^{a_0}v_{r-1}^{a_1}\cdots v_1^{a_{r-1}}v_0^{a_r} \notag \\
&= \sum_{\underline{a}=(a_0,\ldots,a_r)} k_{\underline{a}}
                v_r^{a_r}v_{r-1}^{a_0}\cdots v_1^{a_{r-2}}v_0^{a_{r-1}}
= swap\bigl(ma(f)\bigr)
	(v_r,v_0,\ldots,v_{r-1}) \notag
	\end{align}
    by \eqref{eqn:explicit formula of swap(ma(f))}.
	Using the definition of the $swap$, we rewrite (\ref{fcirc}) in terms of 
	$ma(f)$ as
	\begin{equation} 
	ma(f)(v_r,v_{r-1}-v_r,\ldots,v_0-v_1)=ma(f)(v_{r-1},v_{r-2}-v_{r-1},\ldots,v_0-v_1,v_r-v_0).
    \end{equation}
	By making the change of variables 
    $$
    v_k=\left\{\begin{array}{ll}
    u_0+\cdots+u_{r-k-1} & (0\le k\le r-1), \\
    u_0+\cdots+u_r & (k=r), \\
    \end{array}
    \right.
    $$
    in this equation, we obtain
	\begin{equation}
	ma(f)(u_0+\cdots+u_r,-u_1-\cdots-u_r,u_1,\ldots,u_{r-1}) = ma(f)(u_0,u_1,\ldots,u_r).\end{equation}
	Finally, using relation (\ref{fandb}), we write this in terms of $ma(b)$ as
	\begin{equation}
	ma(b)(-u_1-\cdots-u_r,u_1,\ldots,u_{r-1})=ma(b)(u_1,\ldots,u_r).
	\end{equation}
	which is just the condition of mould push-invariance $ma(b)$ 
	in depth $r$.
The last claim follows from Lemma \ref{lem: ma Lie alg isom} and 
Proposition \ref{prop: pushclosed}.
\end{proof}

For a family  $\mathcal F$ of functions,
we define
$$
 ARI(\mathcal F)_{al+push}:= ARI(\mathcal F)_{al}\cap ARI(\mathcal F)_{push},
$$
which  forms a Lie algebra under the $ari$-bracket
by Propositions \ref{prop: al closed under ari} and \ref{prop: pushclosed}.

\begin{cor}\label{cor: Lie C push closed under the bracket}
    The linear subspace 
    $$\lie{lie}_C^{push}:=Ass_{C,\,\ge1}^{push}\cap \lie{lie}_C$$
    forms Lie subalgebra of $\lie{lie}_C$ (cf. \eqref{eq: lie C}) under  the bracket $\{,\}$ and we have a Lie algebra isomorphism
    \begin{equation}\label{eq: ma lie C push}
    ma:(\lie{lie}_C^{push},\{\ ,\ \})\overset{\sim}{\to}
    (ARI(\mathcal F_\ser)_{al+push}, ari).
\end{equation}
\end{cor}

\begin{proof}
Since the bracket $\{,\}$ is compatible with $ari$ under the map $ma$
by Lemma \ref{lem: ma Lie alg isom}, the result follows from  Proposition \ref{pushinvprop}.
\end{proof}

We also obtain the following corollary.

\begin{cor}
    The map $i_y$  in Lemma \ref{lem: iy} induces Lie algebra inclusion
    \begin{equation}
            i_y: (\lie{lie}_C^{push},\{\ ,\ \}) \hookrightarrow
    \lie{tder}_2^{(y)}.
        \end{equation}
\end{cor}

\begin{proof}
    This follows from Lemma \ref{lem: ma Lie alg isom} and Corollary \ref{cor: Lie C push closed under the bracket}.
\end{proof}

\subsection{Circ-neutrality and the second defining relation of $\lie{lkrv}$}
\label{subsec: circ-neutrality and the second defining relation}
We show how to reformulate the second defining property of elements
of $\lie{lkrv}$ in terms of moulds and prove Proposition \ref{lkrvLie}.

For $B\in \overline{ARI}(\mathcal F)$ with ${\mathcal F}={\mathcal F}_\pol$ or
${\mathcal F}_\ser$, there is a unique decomposition $\{B_n\}_{n\ge1}\subset \overline{ARI}(\mathcal F_\pol)$ such that ${\rm deg} B_n(v_1,\dots,v_r)=n-r$ for $1\le r \le n$ and
$$
B=\sum_{n=1}^\infty B_n.
$$
Actually, if $B$ is represented by
$$
B(v_1,\dots,v_r)
=\sum_{k_1,\dots,k_r=0}^\infty \langle B\ |\ k_1,\dots,k_r \rangle v_1^{k_1}\cdots v_r^{k_r}
$$
for $r\ge1$, each mould $B_n$ ($n\ge1$) can be represented by
$$
B_n(v_1,\dots,v_r)
=\sum_{\substack{k_1+\cdots+k_r=n-r \\ k_i\ge0}} \langle B\ |\ k_1,\dots,k_r \rangle v_1^{k_1}\cdots v_r^{k_r}
$$
for $1\le r\le n$, and $B_n(v_1,\dots,v_r)=0$ for $r>n$.

\begin{defn}
\label{circconstance}
(1).
A mould $B\in \overline{ARI}(\mathcal F)$ is said to be {\it circ-neutral}
\footnote{This property is referred to as "pus-neutral" in \cite[Definition 1.26]{FK} (respectively, \cite[(2.73)]{Ec}) when it holds for all $r \geq 1$ (respectively, for all $r > 1$). The distinction between these two definitions disappears when the depth-one component vanishes; for example, this is the case under the push-invariant condition.
}
(cf. \cite[(2.73)]{Ec}) if  for $r>1$ we have
$$\sum_{i=0}^{r-1} circ^i(B)(v_1,\ldots,v_r)=0.$$ 
For $circ$, see Definition \ref{def: various mould operators}.

(2).
{
Let $B\in \overline{ARI}(\mathcal F)$ with ${\mathcal F}={\mathcal F}_\pol$ or
${\mathcal F}_\ser$, and consider the decomposition $\{B_n\}_{n\ge1}$ of $B$ as above.
The mould $B$ is called {\it circ-constant for the sequence $\{c_n\}_{n>1}\subset \Q$} if we have\footnote{When $r=n$, we note that }
\begin{equation}\label{eq: circ-neutral}
\sum_{i=0}^{r-1} circ^i(B_n)(v_1,\ldots,v_r)= 
c_n\left(\sum_{\substack{a_1+\cdots+a_r=n-r \\ a_i\ge 0}} v_1^{a_1}\cdots v_r^{a_r}\right)
\end{equation}
for $n\ge 2$ and for {$1< r\le n$}.
If $c_n=0$ for any $n\ge 2$, we say that $B$ is {\it circ-neutral}.
}
\end{defn}

\begin{exm}\label{example: lopil is circ-neutral}
We put $v_k:=v_{k-r}$ for $k>r$.
By \eqref{deflopil}, for $r>1$, we have
\begin{align*}
\sum_{i=0}^{r-1} circ^i(lopil)(v_1,\ldots,v_r)
&=\sum_{i=0}^{r-1} lopil(v_{1+i},\ldots,v_{r+i}) \\
&=c_r \sum_{i=0}^{r-1}
\frac{v_{1+i}+\cdots+v_{r+i}}{v_{1+i}(v_{1+i}-v_{2+i})\cdots (v_{r-1+i}-v_{r+i})v_{r+i}} \\
&=c_r \frac{v_{1}+\cdots+v_{r}}{(v_{1}-v_{2})\cdots (v_{r-1}-v_{r})(v_{r}-v_{1})}\sum_{i=0}^{r-1}
\frac{v_{r+i}-v_{r+1+i}}{v_{1+i}v_{r+i}}.
\end{align*}
Because $v_{r+1+i}=v_{1+i}$ for $0\le i\le r-1$, we get
$$
\sum_{i=0}^{r-1}
\frac{v_{r+i}-v_{r+1+i}}{v_{1+i}v_{r+i}}
=\sum_{i=0}^{r-1}
\left( \frac{1}{v_{1+i}}-\frac{1}{v_{r+i}} \right)
=0,
$$
so we obtain
$$
\sum_{i=0}^{r-1} circ^i(lopil)(v_1,\ldots,v_r)=0
$$
for $r>1$, that is, the mould $lopil$ is circ-neutral.
\end{exm}

The following is the corresponding definition in $Ass_2$-side.

\begin{defn}\label{def: circ-const}
(1). We define the $\Q$-linear map $circ:Ass_2 \rightarrow Ass_2$ by
\begin{equation}\label{eqn:definition of circ on Ass2}
circ(x^{a_0}y x^{a_1}y x^{a_2} \cdots yx^{a_r})
:=x^{a_0}yx^{a_r}yx^{a_1} \cdots yx^{a_{r-1}}
\end{equation}
for $r\ge0$ and $a_0,a_1,\dots,a_r\ge0$. 
{In particular, if $a_0=0$ so that we have a monomial of the form $yx^{a_1}\cdots yx^{a_r}$, then the $circ$-operator cyclically permutes the ``chunks'' $yx^{a_i}$.}
It is obvious that this map is a bijection.
By direct calculation, the above map induces a linear isomorphism on $Ass_{C,\,\ge1}$.

(2).
{For $b\in Ass_{C,\ge1}$, $n\ge1$ and for {$1< r\le n$}, let $b_n^r$ be the depth $r$ and the weight $n$ part of $b$, and put $b_{n,0}:=\sum_{r=2}^{n-1}b_n^r$.
 We say that the element $b\in Ass_{C,\ge1}$ is {\it circ-constant for the sequence} $\{c_n\}_{n>1}\subset \Q$, if the following three conditions hold for all $n>1$ :
 \begin{itemize}
 \item[(I).] the polynomial $b_{n,0}^y$ is push-constant for $c_n$ (cf. Definition \ref{polypushconst}),
 \item[(II).] $\langle b\,|\,y^n \rangle=\frac{c_n}{n}$.
 \end{itemize}
 Here we denote $b_{n,0}=xb_{n,0}^x + yb_{n,0}^y$.
If $c_n=0$ for $n>1$, then we say that $b$ is {\it circ-neutral}.
}
\end{defn}

Observe that if $b\in Ass_{C,\ge1}$ and $b^r$ denotes the depth $r$-part of $b$ for $r\ge 1$, then $b$ is circ-neutral if and only if 
$$
\sum_{i=0}^{r-1}circ^i(b^r)=0
$$
for $r\ge1$.

We note that the above definition of circ-neutrality agrees with the one
appearing in Definition \ref{defn:push}.


\begin{exm}
Let $\psi\in \lie{grt}$ be homogeneous of degree $n$
with $(\psi|x^{n-1}y)=1$.
{We saw in Example \ref{eg:poshconstant} that writing $b=\psi=x\psi^x+y\psi^y$, the polynomial 
$\psi^y$ is push-constant for $1$. 
We have $(\psi^y|x^4)=1$, so the polynomial $b'=y\psi^y+\frac{1}{n}y^n$ is strictly push-constant for $1$. 
This polynomial $b'$ is given by $b'=yb'_0+\frac{1}{n}y^n$ where $b'_0$ is the degree 4 push-constant polynomial $\psi^y$ given explicitly in Example \ref{eg:poshconstant}.}

For an example of a  circ-constant mould for the sequence $\{\delta_{n,5}\}_{n>1}$  in homogeneous degree $5$, 
we take $B=swap\bigl(ma(\psi)\bigr)$, 
which has the same coefficients as the above polynomial $y\psi^y$: it is given by 
\begin{align*} 
&B(v_1)=v_1^4\\ 
&B(v_1,v_2)=-2v_1^3+\frac{11}{2}v_1^2v_2-\frac{9}{2}v_1v_2^2+3v_2^3\\ 
&B(v_1,v_2,v_3)=2v_1^2-\frac{11}{2}v_1v_2-\frac{1}{2}v_2^2+\frac{9}{2}v_1v_3 
+2v_2v_3-\frac{1}{2}v_3^2\\ 
&B(v_1,v_2,v_3,v_4)=-v_1+4v_2-6v_3+4v_4\\
&B(v_1,v_2,v_3,v_4,v_5)=\frac{1}{5}.
\end{align*}
\end{exm}

The following result proves that the circ-constance of a homogeneous
Lie polynomial $b$ and that of the associated mould $ma(b)$ are 
related as above.


\begin{lem}\label{lem:swap.ma.circ=circ.swap.ma}
 The following  diagram is commutative: \footnote{See Lemma \ref{lem:bijection mi} for the definition of the map $mi$.}
\begin{equation}\xymatrix{
   Ass_{C,\,\ge1}\ar@{->}[r]^{circ}\ar@{->}[d]_{mi}
   &Ass_{C,\,\ge1}\ar@{->}[d]^{mi}\\
   \overline{ARI}(\mathcal F_\ser)\ar@{->}[r]^{circ}
   &\overline{ARI}(\mathcal F_\ser). }
\end{equation}
\end{lem}

\begin{proof}
Below we will show prove $circ \circ mi(f) = mi \circ circ(f)$ for $f\in Ass_{C,\,\ge1}$ with
$$
f=\sum_{r\ge1}\sum_{\underline{a}=(a_0,a_1,\dots,a_r)}\langle \,f\, |\,\underline{a}\,\rangle \ x^{a_0}y x^{a_1}y x^{a_2} \cdots yx^{a_r}.
$$
In other word, we will prove
\begin{equation}\label{eqn:circswapma=swapmacirc}
circ \circ mi(f)(v_1,\dots,v_r)
=mi \circ circ(f)(v_1,\dots,v_r)
\end{equation}
for $r\ge1$.

When $r=1$, three maps $circ$, $swap$ and $ma$ form identity maps, so we obtain \eqref{eqn:circswapma=swapmacirc}.
Let $r\ge2$.
By Remark \ref{rem:explicit formulation of ma}, we have
\begin{align}\label{eqn:exlicit formula of circswapma(f)}
circ \circ mi(f)(v_1,\dots,v_r)
&=mi(f)(v_r,v_1,\dots,v_{r-1}) \\
&=\sum_{\underline{a}=(0,a_1,\dots,a_r)}\langle \,f\, |\,\underline{a}\,\rangle \ v_{r-1}^{a_1} v_{r-2}^{a_2} \cdots v_1^{a_{r-1}} v_r^{a_r}. \notag
\end{align}
On the other hand, we have
\begin{align*}
circ(f)
&=\sum_{r\ge1}\sum_{\underline{a}=(a_0,a_1,\dots,a_r)}\langle \,f\, |\,\underline{a}\,\rangle \ circ(x^{a_0}y x^{a_1}y x^{a_2} \cdots yx^{a_r}) \\
&=\sum_{r\ge1}\sum_{\underline{a}=(a_0,a_1,\dots,a_r)}\langle \,f\, |\,\underline{a}\,\rangle \ x^{a_0}yx^{a_r}yx^{a_1} \cdots yx^{a_{r-1}} \notag \\
&=\sum_{r\ge1}\sum_{\underline{a}=(a_0,a_1,\dots,a_r)}\langle \,f\, |\,a_0, a_2,\dots,a_r,a_1\,\rangle \ x^{a_0}y x^{a_1}y x^{a_2} \cdots yx^{a_r}, 
\end{align*}
so we get
$$
\langle \,circ(f)\, |\,\underline{a}\,\rangle
=\langle \,f\, |\,a_0, a_2,\dots,a_r,a_1\,\rangle
$$
for $a_0,a_1,\dots,a_r\ge0$.
By Remark \ref{rem:explicit formulation of ma}, we have
\begin{align*}
mi \circ circ(f)(v_1,\dots,v_r)
&=\sum_{\underline{a}=(0,a_1,\dots,a_r)}\langle \,circ(f)\, |\,\underline{a}\,\rangle \ v_{r}^{a_1} v_{r-1}^{a_2} \cdots v_2^{a_{r-1}} v_1^{a_r} \\
&=\sum_{\underline{a}=(0,a_1,\dots,a_r)}\langle \,f\, |\,0, a_2,\dots,a_r,a_1\,\rangle \ v_{r}^{a_1} v_{r-1}^{a_2} \cdots v_2^{a_{r-1}} v_1^{a_r} \\
&=\sum_{\underline{a}=(0,a_1,\dots,a_r)}\langle \,f\, |\,\underline{a}\,\rangle \ v_{r}^{a_r} v_{r-1}^{a_1} \cdots v_2^{a_{r-2}} v_1^{a_{r-1}}.
\end{align*}
By comparing this and \eqref{eqn:exlicit formula of circswapma(f)}, we obtain \eqref{eqn:circswapma=swapmacirc} for $r\ge2$.
\end{proof}

\begin{prop}\label{circneutprop}
{Let $b\in Ass_{C,\ge1}$.
The following three statements are equivalent: \\
(i). The mould $B=mi(b)\in \overline{ARI}(\mathcal F_\ser)$ is circ-constant for $\{c_n\}_{n>1}\subset \Q$. \\
(ii). The element $b\in Ass_{C,\ge1}$ satisfies
\begin{equation}\label{eqn:middle condition of circ-constant for b^r_n}
\sum_{i=0}^{r-1} circ^i(b^r_n)
=(-1)^{n-r}c_n\left( \sum_{\substack{a_1+\cdots+a_r=n \\ a_i\ge 1}} C_{a_r,a_{r-1},\dots,a_1} \right)
\end{equation}
for all $n\ge2$ and {$1< r\le n$}.
Here the element $b^r_n\in Ass_{C,\ge1}$ means the depth $r$ and the weight $n$ part of $b$.\\
(iii). The element $b\in Ass_{C,\ge1}$ is circ-constant for $\{c_n\}_{n>1}\subset \Q$.}
\end{prop}
\begin{proof}
First, we prove the equivalence between $(i)$ and $(ii)$.
Note that, by direct calculation, we have
$$
B_n(v_1,\dots,v_r)=mi\left( b^r_n \right)(v_1,\dots,v_r)
$$
for $n\ge1$. 
We apply the map $mi$ to both sides of \eqref{eqn:middle condition of circ-constant for b^r_n}.
Then the left-hand side of \eqref{eqn:middle condition of circ-constant for b^r_n} is calculated as follows:
\begin{align}\label{eqn:cal. of left-hand side of middle condition of circ-constant}
\sum_{i=0}^{r-1} mi\circ circ^i(b^r_n)(v_1,\dots,v_r)
&=\sum_{i=0}^{r-1} circ^i(mi(b^r_n))(v_1,\dots,v_r) \\
&=\sum_{i=0}^{r-1} circ^i(B_n)(v_1,\dots,v_r). \nonumber
\end{align}
Here, we used Lemma \ref{lem:swap.ma.circ=circ.swap.ma} in the first equality.
On the other hand, the right-hand side of \eqref{eqn:middle condition of circ-constant for b^r_n} is calculated as follows:
\begin{align}\label{eqn:cal. of right-hand side of middle condition of circ-constant}
&(-1)^{n-r}c_n\left( \sum_{\substack{a_1+\cdots+a_r=n \\ a_i\ge 1}} mi(C_{a_r,a_{r-1},\dots,a_1}) \right) \\
&=(-1)^{n-r}c_n\left( \sum_{\substack{a_1+\cdots+a_r=n \\ a_i\ge 1}} (-1)^{n-r} v_r^{a_r-1} \cdots v_1^{a_1-1} \right) 
=c_n\left( \sum_{\substack{a_1+\cdots+a_r=n-r \\ a_i\ge 0}} v_1^{a_1} \cdots v_r^{a_r} \right) \nonumber
\end{align}
Here, we used \eqref{eq: mi} in the first equality.
By combining \eqref{eqn:cal. of left-hand side of middle condition of circ-constant} and \eqref{eqn:cal. of right-hand side of middle condition of circ-constant}, we know that \eqref{eqn:middle condition of circ-constant for b^r_n} is equivalent to the following equation:
$$
\sum_{i=0}^{r-1} circ^i(B_n)(v_1,\dots,v_r)
=c_n\left( \sum_{\substack{a_1+\cdots+a_r=n-r \\ a_i\ge 0}} v_1^{a_1} \cdots v_r^{a_r} \right).
$$
This is nothing but the equation \eqref{eq: circ-neutral}, that is, the mould $B=mi(b)$ is circ-constant for $\{c_n\}_{n>1}$.
Hence, condition $(i)$ is equivalent to condition $(ii)$.\\

\smallskip
Next, we prove the equivalence between $(iii)$ and $(i)$ (or $(ii)$).
Firstly, we prove the equivalence between \eqref{eqn:middle condition of circ-constant for b^r_n} for $r=n$ and (II) in Definition \ref{def: circ-const}.(2).
We calculate both sides of \eqref{eqn:middle condition of circ-constant for b^r_n} for $r=n$.
Because $b_n^n=\langle b\,|\,y^n \rangle y^n$, the left-hand side of \eqref{eqn:middle condition of circ-constant for b^r_n} is calculated as follows:
\begin{align*}
\sum_{i=0}^{n-1} circ^i(b^n_n)
=\sum_{i=0}^{n-1} \langle b\,|\,y^n \rangle y^n
=n\langle b\,|\,y^n \rangle y^n.
\end{align*}
On the other hand, the right-hand side of \eqref{eqn:middle condition of circ-constant for b^r_n} is calculated as follows:
\begin{align*}
(-1)^{n-n}c_n\left( \sum_{\substack{a_1+\cdots+a_n=n \\ a_i\ge 1}} C_{a_n,a_{n-1},\dots,a_1} \right)
=c_nC_{\scriptsize \underbrace{1,1,\dots,1}_{n}}
=c_ny^n.
\end{align*}
Comparing these coefficients, we get $\langle b\,|\,y^n \rangle=\frac{c_n}{n}$.
This means condition (II) in Definition \ref{def: circ-const}.(2).

Secondly, we prove the equivalence between \eqref{eqn:middle condition of circ-constant for b^r_n} for $1<r<n$ and (I) in Definition \ref{def: circ-const}.(2).
We calculate both sides of \eqref{eqn:middle condition of circ-constant for b^r_n} for $1<r< n$.
Because we have $b_n^r=x(b_n^r)^x +y(b_n^r)^y$ and $circ(yf)=y\cdot push(f)$ for $f\in Ass_C$, the left hand side of
\eqref{eqn:middle condition of circ-constant for b^r_n}
is given below:
\begin{align*}
\sum_{i=0}^{r-1} circ^i(b^r_n)
&=\sum_{i=0}^{r-1} circ^i\bigl(x(b^r_n)^x\bigr) + \sum_{i=0}^{r-1} circ^i\bigl(y(b^r_n)^y\bigr) \\
&=y\cdot\sum_{i=0}^{r-1} push^i\bigl((b^r_n)^y\bigr) + \sum_{i=0}^{r-1} circ^i\bigl(x(b^r_n)^x\bigr).
\end{align*}
On the other hand, 
since the only term in $C_{a_r,a_{r-1},\dots,a_1}$ that starts with $y$ is
$(-1)^{a_1+\cdots+a_r-r}yx^{a_r-1}\cdots yx^{a_1-1}$,
the right-hand side of
\eqref{eqn:middle condition of circ-constant for b^r_n}
is given as follows:
\begin{align*}
&(-1)^{n-r}c_n\left( \sum_{\substack{a_1+\cdots+a_r=n \\ a_i\ge 1}} C_{a_r,a_{r-1},\dots,a_1} \right) \\
&=(-1)^{n-r}c_n\left( \sum_{\substack{a_1+\cdots+a_r=n \\ a_i\ge 1}} (-1)^{a_1+\cdots+a_r-r}yx^{a_r-1}\cdots yx^{a_1-1} \right)
+\varphi \\
&=y\cdot \left( c_n\sum_{\substack{a_1+\cdots+a_r=n-r \\ a_i\ge 0}} x^{a_r}yx^{a_{r-1}}\cdots yx^{a_1} \right)
+\varphi.
\end{align*}
Here  $\varphi$ is its terms which belong to  $x\cdot Ass_2$.
By comparing the terms which start with  $y$, we get
$$
\sum_{i=0}^{r-1} push^i\bigl((b^r_{n,0})^y\bigr)
{=}\sum_{i=0}^{r-1} push^i\bigl((b^r_n)^y\bigr)
=c_n\sum_{\substack{a_1+\cdots+a_r=n-r \\ a_i\ge 0}} x^{a_1}yx^{a_{2}}\cdots yx^{a_r},
$$
for $1<r<n$
where we put $b_{n,0}:=\sum_{s=2}^{n-1}b_n^s$.
Because $(b^r)^y=(b^y)^{r-1}$ for $b\in Ass_{C,\ge1}$ and $r\ge1$, we get
$$
\sum_{i=0}^{r-1} push^i\bigl((b^y_{n,0})^{r-1}\bigr)
=c_n\sum_{\substack{a_1+\cdots+a_r=n-r \\ a_i\ge 0}} x^{a_1}yx^{a_{2}}\cdots yx^{a_r},
$$
for $1<r<n$.
By substituting $r$ with $r+1$, we obtain the following:
$$
\sum_{i=0}^{r} push^i\bigl((b^y_{n,0})^{r}\bigr)
=c_n\sum_{\substack{a_0+\cdots+a_r=(n-1)-r \\ a_i\ge 0}} x^{a_0}yx^{a_{1}}\cdots yx^{a_r},
$$
for {$0<r<n-1$}.
Since the weight of $b_{n,0}^y$ is $n-1$, this equation means that $b_{n,0}^y$ is push-constant for $c_n$.
Hence, we obtain (I) in Definition \ref{def: circ-const}.(2).

On the other hand, we assume (I) in Definition \ref{def: circ-const}.(2), that is, we assume
\begin{equation}\label{eqn:push-constant for 1<r<n}
\sum_{i=0}^{r-1} push^i\bigl((b^y_{n,0})^{r-1}\bigr)
=c_n\sum_{\substack{a_1+\cdots+a_r=n-r \\ a_i\ge 0}} x^{a_1}yx^{a_{2}}\cdots yx^{a_r},
\end{equation}
for $1<r<n$.
By $b_{n,0}^{r}\in Ass_{C,\ge1}^{(r,n)}$, there exists $\langle b_{n,0}^{r} \ |\ a_r,\dots,a_1 \rangle\in\Q$ such that
$$
b_{n,0}^{r}
=\sum_{\substack{a_1+\cdots+a_r=n-r \\ a_k\ge0}}
\langle b_{n,0}^{r} \ |\ a_r,\dots,a_1 \rangle C_{a_r+1,a_{r-1}+1,\dots,a_1+1}
$$
for $1<r<n$.
We have $(C_{a_r+1,\dots,a_1+1})^y=(-1)^{n-r}x^{a_r}yx^{a_{r-1}}\cdots yx^{a_1}$ with $a_1+\cdots+a_r=n-r$, so for $1<r<n$, we get
\begin{align*}
&\sum_{i=0}^{r-1} push^i\bigl((b^y_{n,0})^{r-1}\bigr)
=\sum_{i=0}^{r-1} push^i\bigl((b^{r}_{n,0})^y\bigr) \\
&=\sum_{i=0}^{r-1} push^i\left( \sum_{\substack{a_1+\cdots+a_r=n-r \\ a_k\ge0}}
\langle b^{r}_{n,0} \ |\ a_r,\dots,a_1 \rangle
(-1)^{n-r}x^{a_r}yx^{a_{r-1}}\cdots yx^{a_1} \right) \\
&=\sum_{i=0}^{r-1} \sum_{\substack{a_1+\cdots+a_r=n-r \\ a_k\ge0}}
\langle b^{r}_{n,0} \ |\ a_r,\dots,a_1 \rangle
(-1)^{n-r}x^{a_{r+i}}yx^{a_{r-1+i}}\cdots yx^{a_{1+i}} \\
&= \sum_{\substack{a_1+\cdots+a_r=n-r \\ a_k\ge0}}
\left(
(-1)^{n-r}\sum_{i=0}^{r-1}
\langle b^{r}_{n,0} \ |\ a_{r+i},\dots,a_{1+i} \rangle
\right)
x^{a_{1}}yx^{a_{2}}\cdots yx^{a_{r}}.
\end{align*}
By comparing this with \eqref{eqn:push-constant for 1<r<n}, we get
\begin{equation}\label{eqn:cyclic sum of coeff. of bnr=cn}
(-1)^{n-r}\sum_{i=0}^{r-1}
\langle b^{r}_{n,0} \ |\ a_{r+i},\dots,a_{1+i} \rangle
=c_n
\end{equation}
for $1<r<n$.
By using \eqref{eq: mi} and $b_{n,0}^r=b_n^r$ for $1<r<n$, we have
\begin{align*}
&\sum_{i=0}^{r-1} circ^i\circ mi(b^r_n)(v_1,\dots,v_r) \\
&=\sum_{i=0}^{r-1} circ^i
\left(
\sum_{\substack{a_1+\cdots+a_r=n \\ a_k\ge1}}
\langle b_{n,0}^r \ |\ a_r,\dots,a_1 \rangle
(-1)^{n-r} v_1^{a_1-1}v_2^{a_2-1}\cdots v_r^{a_r-1}
\right) \\
&=\sum_{i=0}^{r-1} 
\sum_{\substack{a_1+\cdots+a_r=n \\ a_k\ge1}}
(-1)^{n-r}\langle b_{n,0}^r \ |\ a_r,\dots,a_1 \rangle
v_{1+i}^{a_1-1}v_{2+i}^{a_2-1}\cdots v_{r+i}^{a_r-1} \\
&=\sum_{\substack{a_1+\cdots+a_r=n \\ a_k\ge1}}
(-1)^{n-r}\left(
\sum_{i=0}^{r-1} \langle b_{n,0}^r \ |\ a_{r+i},\dots,a_{1+i} \rangle
\right)
v_{1}^{a_1-1}v_{2}^{a_2-1}\cdots v_{r}^{a_r-1} \\
&=c_n\cdot\sum_{\substack{a_1+\cdots+a_r=n \\ a_k\ge1}}
v_{1}^{a_1-1}v_{2}^{a_2-1}\cdots v_{r}^{a_r-1}.
\end{align*}
Here, we used \eqref{eqn:cyclic sum of coeff. of bnr=cn} in the last equality.
This equation means {\eqref{eq: circ-neutral}} for $1<r< n$.
Since  \eqref{eq: circ-neutral} 
has been shown to be equivalent to \eqref{eqn:middle condition of circ-constant for b^r_n} in the proof of the equivalence between (i) and (ii), the proof is now complete.
\end{proof}

\begin{cor}
Let $b\in \lie{lie}_C$ be of homogeneous weight $n\ge 3$.
Then $b$ is circ-neutral if and only if $mi(b)$ is a circ-neutral mould.
\end{cor}

\begin{proof}
By putting $c_n=0$ for $n>1$ in Proposition \ref{circneutprop}, we obtain this claim.
\end{proof}

\vspace{.2cm}
The notion of circ-constance will play a role later in \S \ref{subsec:the map from krv to krvell}. 
The concept of circ-neutrality  also plays a fundamental role in mould theory.

\begin{defn}
Let $\mathcal F$ be a family of functions.
    We define
    \footnote{In the case when $\mathcal F=\mathcal F_\rat$, it 
    is  $ARI(\{e\})_{push/pusnu}$  in the notation of \cite{FK}.}
    $$
    {ARI}(\mathcal F)_{push/circneut}\subset {ARI}(\mathcal F)
    $$
    as the linear subspace of  $ARI(\mathcal F)$
    consisting of moulds which are push-invariant and  whose swaps are circ-neutral.
\end{defn}

An important property of ${ARI}(\mathcal F)_{circneut}$ is the compatibility with  the $ari$-bracket.

\begin{prop}
\label{circneutclosed}
The space $\overline{ARI}(\mathcal F)_{circneut}$
of circ-neutral moulds in $\overline{ARI}$
forms a Lie subalgebra  of
 $\overline{ARI}(\mathcal F)$
 under the $\overline{ari}$-bracket.  
Further, the space ${ARI}(\mathcal F)_{push/circneut}$
forms a Lie subalgebra  of
 ${ARI}(\mathcal F)$
 under the ${ari}$-bracket.  
\end{prop}

\begin{proof}
 A direct complete proof of this result, along with its generalization, is presented in \cite[Proposition 1.30 and Theorem 1.32]{FK}.
\end{proof}

\begin{prop}\label{prop: ARI al push circneut Lie algebra}
Let $\mathcal F$ be a family of functions.
    We define
$$
 ARI(\mathcal F)_{al+push/circneut}:= ARI(\mathcal F)_{al}\cap ARI(\mathcal F)_{push/circneut}.
$$
Then it forms a Lie algebra under the $ari$-bracket.
\end{prop}

\begin{proof}
It follows from Propositions \ref{prop: al closed under ari} and \ref{circneutclosed}.
\end{proof}

\begin{prop}\label{prop:isom lkrvARI}
The map $ma$ gives a linear space isomorphism
\begin{equation*}
ma:\lie{lkrv}\buildrel\sim\over\rightarrow ARI(\mathcal F_\ser)_{al+push/circneut},
\end{equation*}

\end{prop}

\begin{proof}
The first claim  is because
we showed that a polynomial $b$ lies in $\lie{lkrv}$, i.e.~$b$ is a  Lie polynomial that is push-invariant and circ-neutral, if and only if
the associated mould $ma(b)$ is alternal (by Lemma \ref{Liestuff} (ii)), 
push-invariant (by Proposition \ref{pushinvprop}) and its swap is circ-neutral (by Proposition \ref{circneutprop}).
The second claim follows from Lemma \ref{lem: ma Lie alg isom} and 
Proposition \ref{prop: ARI al push circneut Lie algebra}.
\end{proof}

\bigskip
{\bf Proof of Proposition \ref{lkrvLie}.}
By the above proposition  with  Lemma \ref{lem: ma Lie alg isom} and 
Proposition \ref{prop: ARI al push circneut Lie algebra},
we see that  $\lie{lkrv}$ forms a Lie algebra under the bracket $\{,\}$.
It is immediate  to see that it is bigraded by weight and depth. 
\qed
\bigskip

We note that there is  a Lie algebra isomorphism 
\begin{equation}\label{lkrvARI}
    ma:(\lie{lkrv},\{,\})\overset{\sim}{\to}
    (ARI(\mathcal F_\ser)_{al+push/circneut},ari).
\end{equation}


	%

\subsection{The inclusion $\lie{ls}\hookrightarrow \lie{lkrv}$ (Theorem \ref{firstthm})}\label{sec: proof of firstthm}
    


In order to prove the theorem, we first reformulate the statement in terms of moulds and give its proof.

 \begin{lem}\label{lem:ma ls ARIal/al}        
The map $ma$ gives a {Lie algebra} isomorphism
\begin{equation*}
    ma:\lie{ls}\buildrel\sim\over\rightarrow  ARI(\mathcal F_\ser)_{\underline{al}/\underline{al}}.
\end{equation*}
For the left hand side,  see Definition \ref{def: ls} 
and the right hand side, see Definition \ref{def: ARI al/al}.
 \end{lem} 
 
\begin{proof}
It has been demonstrated in \cite{S2} that the result holds, 
and this  claim has been further extended in \cite{FK}
\end{proof}

				
\begin{thm}\label{mouldfirstthm} 
Let $\mathcal F$ be a family of functions.
There is an inclusion of mould subspaces
\begin{equation}\label{eq:mouldfirstthm}
ARI({\mathcal F})_{\underline{al}/\underline{al}}\subset ARI({\mathcal F})_{al+push/circneut}.
\end{equation}
Moreover in depths $r\le 3$, the following holds for any family $\mathcal F$
$$ARI(\mathcal F)^r\cap ARI(\mathcal F)_{\underline{al}/\underline{al}}
=ARI(\mathcal F)^r\cap ARI(\mathcal F)_{al+push/circneut}.$$
\end{thm}

\begin{proof}
The result that $A\in ARI_{\underline{al}/\underline{al}}(\mathcal F)$ 
is push-invariant can be found in \cite[Lemma 2.5.5]{S2}.
and that $swap(A)$ is circ-neutral 
can be found in \cite[Proposition 3.12]{FK}.



Let us now prove the isomorphism in the cases $r=1,2,3$.:
The case $r=1$ is trivial.
Indeed, elements  in depth one components of elements in both
$ARI(\mathcal F)_{\underline{al}/\underline{al}}$ and
$ARI({\mathcal F})_{al+push/circneut}$ 
must satisfy only one condition: they must be even functions.


Now consider the case $r=2$. 
Let $A\in ARI({\mathcal F})_{al+push/circneut}$ be concentrated in depth 2. 
The circ-neutral property of the swap is explicitly given in depth 2 by $swap(A)(v_1,v_2)+swap(A)(v_2,v_1)=0$.
But this is also the alternality condition on $swap(A)$, so $A\in ARI({\mathcal F})_{al/al}$. 
The isomorphism in depth 2 is thus trivial. 

Finally, we consider the case $r=3$. Let $A\in ARI({\mathcal F})_{al+push/circneut}$
be concentrated in depth 3, and let $B=swap(A)$.  
Again, we only need to show that $B$ is alternal, which in depth 3 means that $B$ must satisfy
he single equation
\begin{equation}\label{alternalB}
B(v_1,v_2,v_3)+B(v_2,v_1,v_3)+B(v_2,v_3,v_1)=0.
\end{equation}

The circ-neutrality condition on $B$ is given by
					\begin{equation}\label{circneutB}
					B(v_1,v_2,v_3)+B(v_3,v_1,v_2)+B(v_2,v_3,v_1)=0.
					\end{equation}
					It is enough to show that $B$ satisfies the equality
					\begin{equation}\label{mantarB3}
					B(v_1,v_2,v_3)=B(v_3,v_2,v_1),
					\end{equation}
					since applying this to the middle term of (\ref{circneutB}) immediately yields
					the alternality property (\ref{alternalB}) in depth $3$. So let us show
					how to prove (\ref{mantarB3}).

					We rewrite the push-invariance condition in the $v_i$, which gives
					\begin{align}\label{(1)} B(v_1,v_2,v_3)&=B(v_2-v_1,v_3-v_1,-v_1)\\
						\label{(2)}&=B(v_3-v_2,-v_2,v_1-v_2)\\
						\label{(3)}&=B(-v_3,v_1-v_3,v_2-v_3).
						\end{align}

						Making the variable change exchanging $v_1$ and $v_3$, this gives
						\begin{align}\label{(4)}
						B(v_3,v_2,v_1)&=B(v_2-v_3,v_1-v_3,-v_3)\\
								\label{(5)}&=B(v_1-v_2,-v_2,v_3-v_2)\\
								\label{(6)}&=B(-v_1,v_3-v_1,v_2-v_1).
								\end{align}

								By (\ref{(1)}), the term $B(v_2-v_1,v_3-v_1,-v_1)$ is circ-neutral with respect
								to the cyclic permutation of $v_1,v_2,v_3$, so we have
								\begin{equation}\label{(star)}
								B(v_2-v_1,v_3-v_1,-v_1)=-B(v_3-v_2,v_1-v_2,-v_2)-B(v_1-v_3,v_2-v_3,-v_3).
								\end{equation}
								But the circ-neutrality of $B$ also lets us cyclically permute the three arguments
								of $B$, so we also have
								\begin{equation*}
								-B(v_3-v_2,v_1-v_2,-v_2)=B(-v_2,v_3-v_2,v_1-v_2)+B(v_1-v_2,-v_2,v_3-v_2).
								\end{equation*}
Using (\ref{(1)}) and substituting this into the right-hand side of (\ref{(star)}) 
	yields
	\begin{align}\label{(starstar)}
	B(v_1,v_2,v_3)&=B(-v_2,v_3-v_2,v_1-v_2)\notag\\
			&\ +B(v_1-v_2,-v_2,v_3-v_2)-B(v_1-v_3,v_2-v_3,-v_3).
			\end{align}
			Now, exchanging $v_1$ and $v_2$ in (\ref{(6)}) gives
			\begin{equation*}
			B(v_3,v_1,v_2)=B(-v_2,v_3-v_2,v_1-v_2),
			\end{equation*}
			and doing the same with (\ref{(4)}) gives
			\begin{equation*}
			B(v_3,v_1,v_2)=B(v_1-v_3,v_2-v_3,-v_3).
			\end{equation*}
			Substituting these two expressions as well as (\ref{(5)}) into the 
			right-hand side of (\ref{(starstar)}), we obtain the desired equality 
			(\ref{mantarB3}).  This concludes the proof of Theorem \ref{mouldfirstthm}. 
 \end{proof}           

The result of Theorem \ref{mouldfirstthm} is related to the injective map 
\begin{equation}\label{dsellinkrvell}
\lie{ds}_{ell}\hookrightarrow\lie{krv}_{ell}
\end{equation}discussed in the introduction (cf.~\eqref{commdiag0}).

\begin{cor}\label{mouldfirstcor} 
Let $\mathcal F$ be a family of functions.
The inclusion of Theorem \ref{mouldfirstthm} extends to
an inclusion of spaces
$$ARI(\mathcal F)_{\underline{al}*\underline{al}}\subset ARI(\mathcal F)_{al+push*circneut},$$
where 
$ARI(\mathcal F)_{\underline{al}*\underline{al}}$ is given in 
Definition \ref{def: ARI al/al}
and $ARI(\mathcal F)_{al+push*circneut}$ is the set of moulds 
which are alternal and push-invariant and with swap being circ-neutral
up to addition of a constant-valued mould.
\end{cor}

\begin{proof}
Let $A\in ARI(\mathcal F)_{\underline{al}*\underline{al}}$. We know
by Lemma 2.5.5 of \cite{S2} that $A$ is push-invariant. Let $B:=swap(A)$ and
let $B_0$ be the constant mould such that $B+B_0$ is alternal. Then by 
Theorem \ref{mouldfirstthm}, $B+B_0$ is also circ-neutral, and therefore
by definition $B=swap(A)$ is $*$circ-neutral, i.e.~circ-neutral up to 
addition of a constant-valued mould.
\end{proof}
The above corollary will be used in the next section.

\bigskip
\noindent
{\bf Proof of Theorem \ref{firstthm}.}
By Theorem \ref{mouldfirstthm}, there is an inclusion
				in terms of moulds as 
\begin{equation*}
ARI({\mathcal F}_\ser)_{\underline{al}/\underline{al}}\subset ARI({\mathcal F}_\ser)_{al+push/circneut}.
\end{equation*}
By  \eqref{eq:ma ls ARIal/al} and \eqref{lkrvARI},
the theorem is proven.
\qed

\vspace{.3cm}
\begin{rem}
We speculate that the inclusions of Theorem 
\ref{mouldfirstthm} and Corollary \ref{mouldfirstcor} 
may extend to isomorphisms. 
But even the proof of the simple equality (\ref{mantarB3}) is
surprisingly complicated in depth $3$, let alone in higher depth.  
Computational experiments suggest the following general pattern:
If $A\in ARI(\mathcal F)_{al+push/circneut}$ and 
$B=swap(A)$, then $B$ is mantar-invariant (cf. Definition \ref{def: various mould operators}),
that is, 
for all $r>0$, we  have
\begin{equation}\label{mantarB}
B(v_1,\ldots,v_r)=(-1)^{r-1}B(v_r,\ldots,v_1).
\end{equation}
\end{rem}

                    The mantar invariance (\ref{mantarB}) would also yield the following
					useful partial result, which is the mould analog for $\lie{lkrv}$ of a
					result that is well-known for $\lie{ls}$, namely that the bigraded
					part $\lie{ls}_n^r=0$ when $n\not\equiv r$ mod 2.

\begin{lem}\label{lem: parity al+push/circneut}
Fix $1\le r\le n$. Let 
$A\in ARI(\mathcal F_\ser)^{(r,n)}\cap
ARI({\mathcal F}_\ser)_{al+push/circneut}$ and let $B=swap(A)$.  
Assume that $B$ satisfies (\ref{mantarB}). Then if $n-r$ is odd, $A=0$.
\end{lem}

For $ARI(\mathcal F_\ser)^{(r,n)}$, see Lemma \ref{Liestuff}.

\begin{proof}
We denote by \textit{mantar} the operator on $\overline{ARI}(\mathcal F_\ser)$, which is defined in the same manner as on $ARI$, except that the variables $v_i$ are used instead of $u_i$
Then it is easy to check the following identity of operators noted by 
\'Ecalle: 
	$$neg\circ push=mantar\circ swap\circ mantar\circ swap,$$
Let $A\in ARI({\mathcal F}_\ser)_{al+push/circneut}$; then $A$ is push-invariant, so applying the
	left-hand operator to $A$ gives $neg(A)$.  Assuming (\ref{mantarB}) for
	$B=swap(A)$, i.e.~assuming that $B=mantar(B)$, we see that applying the 
	right-hand operator to $A$ fixes $A$ since on the one hand $swap\circ swap=id$
	and on the other, $mantar(A)=A$ for all alternal moulds (cf.~\cite{S2},
			Lemma 2.5.3).  Thus $A$ must satisfy $neg(A)=A$, i.e.~if $A\ne 0$ then
	the degree $d=n-r$ of $A$ must be even.
\end{proof}

The following parity result, which is the analogy for $\lie{lkrv}$ of
	the similar well-known result on $\lie{ls}$.

	\begin{cor} 
    If the swaps of all elements of $ARI({\mathcal F}_\ser)_{al+push/circneut}$
	are $mantar$-invariant, then $ARI(\mathcal F_\ser)^{(r,n)}\cap
	ARI({\mathcal F}_{\ser})_{al+push/circneut}=0$ whenever $r-n$ is odd, i.e. 
	$$\lie{lkrv}_n^r=0\ \ \ {\rm when}\ n\not\equiv r\ {\rm mod}\ 2$$
	\end{cor}
    
\begin{proof}
By Proposition \ref{prop:isom lkrvARI}, there is  an isomorphism
\begin{equation*}\label{bigradediso}
\lie{lkrv}_n^r\simeq 
ARI(\mathcal F_\ser)^{(r,n)}\cap ARI(\mathcal F_\ser)_{al+push/circneut},
\end{equation*}
of each bigraded piece.
Thus our claim follows from Lemma \ref{lem: parity al+push/circneut}.
\end{proof}

	\section{The elliptic Kashiwara-Vergne Lie algebra $\lie{krv}_{ell}$}\label{sec:The elliptic Kashiwara-Vergne Lie algebra}

	\vspace{.2cm}
	In this section we follow the procedure of \cite{S3} for the double shuffle
Lie algebra to define a natural candidate for the elliptic Kashiwara-Vergne Lie algebra, closely related to the linearized
Kashiwara-Vergne Lie algebra, and give some of its properties (Theorem \ref{krvellisLie}).

	\subsection{Definition of the elliptic Kashiwara-Vergne Lie algebra $\lie{krv}_{ell}$}\label{subsec:Definition of the elliptic Kashiwara-Vergne Lie algebra}
    We first introduce the $Dari$-bracket, a secondary Lie bracket structure on $ARI(\mathcal F_\Lau)$. We then demonstrate that $\lie{krv}_{ell}$  is closed under this bracket, thereby establishing its Lie algebraic structure.

	\subsubsection{Dari bracket}\label{411}

\begin{defn}
Take $\mathcal F=\mathcal F_\Lau$.
We define the linear map
$$\Delta:ARI(\mathcal F_\Lau)\to ARI(\mathcal F_\Lau) $$ 
given by
	\begin{equation}\label{eq: map Delta}
\Delta(A)(u_1,\ldots,u_r)=u_1\cdots u_r(u_1+\cdots +u_r)A(u_1,\ldots,u_r)
	\end{equation}
	for $r\ge 1$.
    We define 

\begin{equation}\label{Daribrack}
Dari(A,B) = \Delta\Bigl(ari\bigl(\Delta^{-1}(A),\Delta^{-1}(B)\bigr)
\Bigr).
\end{equation}
\end{defn}

This definition implies that $\Delta$ gives an isomorphism of Lie algebras
\begin{equation}\label{Deltaiso}
\Delta:ARI(\mathcal F_\Lau)_{ari}\buildrel\sim\over\rightarrow ARI(\mathcal F_\Lau)_{Dari}.
\end{equation}
Here $ARI(\mathcal F_\Lau)_{Dari}$ is the Lie algebra 
$ARI(\mathcal F_\Lau)$ equipped with ${Dari}$-bracket.

\begin{rem}
It is shown in \cite{S3}, Prop. 3.2.1 that we have a second definition
for the $Dari$-bracket, which is more complicated but sometimes very
useful in certain proofs.
Let $$dar: ARI(\mathcal F_\Lau)\to ARI(\mathcal F_\Lau)$$
denote the mould operator
defined by 
$$dar(A)(u_1,\ldots,u_r)=u_1\cdots u_r\,A(u_1,\ldots,u_r).$$
We begin by introducing,
for each $A\in ARI(\mathcal F)$, an associated derivation $Darit(A)$ of $ARI(\mathcal F)_{lu}$
by the following formula:
\begin{equation}\label{defDarit}
Darit(A)=dar\circ\Bigl(-arit\bigl(\Delta^{-1}(A)\bigr)+ad\bigl(\Delta^{-1}(A)\bigr)\Bigr)\circ dar^{-1},
\end{equation}
where $ad(A)\cdot B=lu(A,B)$.  Then $Dari$ corresponds to
the bracket of derivations, in the sense that
\begin{equation}\label{Daridef2}
Dari(A,B)=Darit(A)\cdot B-Darit(B)\cdot A.
\end{equation}
\end{rem}

\begin{defn}
     Let $ARI(\mathcal F_\Lau)^\Delta$ denote the space of 
	moulds $A$ such that $\Delta(A)$  (for $\Delta$, see \eqref{eq: map Delta})
    is in $ARI(\mathcal F_\ser)$, i.e.~the denominator
			of  $A$ is ``at worst'' $u_1\cdots u_r(u_1+\cdots+u_r)$.
	We write $ARI(\mathcal F_\Lau)^\Delta_{\mathcal P}$ for the space of moulds in $ARI(\mathcal F_\Lau)^\Delta\cap ARI(\mathcal F_\Lau)_{\mathcal P}$,
	where ${\mathcal P}$ may represent any (or no) properties on moulds in $ARI(\mathcal F_\Lau)$; we will
	consider properties ${\mathcal P}$ such as for example $al$, $push$, combinations of these
	etc.
\end{defn}

\begin{lem}\label{lem: ARI Delta under aribracket}
The space $ARI(\mathcal F_\Lau)^\Delta$ forms a Lie algebra under the ari-bracket \eqref{Poisson}.    
\end{lem}

\begin{proof}
    It follows from  \cite[Proposition 4.2]{E2}.
\end{proof}

\begin{lem}\label{lem: ARI Fser Dari}
The space $ARI(\mathcal F_\ser)$ forms a Lie algebra under the Dari-bracket
\eqref{Daribrack}.
\end{lem}

\begin{proof}
    Let $A,B\in ARI(\mathcal F_\ser)$.
    Then we have $\Delta^{-1}(A), \Delta^{-1}(B)\in ARI(\mathcal F_\Lau)^\Delta$.
    By Lemma \ref{lem: ARI Delta under aribracket}, 
    $ari(\Delta^{-1}(A), \Delta^{-1}(B))\in ARI(\mathcal F_\Lau)^\Delta$.
    Therefore  $\Delta(ari(\Delta^{-1}(A), \Delta^{-1}(B)))\in ARI(\mathcal F_\ser)$.
    By \eqref{Daribrack}, we get $Dari(A,B)\in ARI(\mathcal F_\ser)$.
\end{proof}

For convenience, when considering the Lie algebra structure on
$ARI(\mathcal F_\ser)$ with respect to the Dari-bracket,
we denote it by $ARI(\mathcal F_\ser)_{Dari}$,

    

\begin{defn}\label{defnkrvell} 
{\it The mould-version elliptic Kashiwara-Vergne linear space}  is the subspace of $ARI(\mathcal F_\ser)$ given by
$$\Delta\bigl(ARI(\mathcal F_\Lau)^{\Delta}_{al+push*circneut}\bigr).$$
\end{defn}

The elliptic Kashiwara-Vergne linear space 
$\lie{krv}_{ell} (\subset \lie{lie}_C)$
introduced in Definition \ref{mouldkrvell}
agrees with the above mould version
under the map $ma$ given in \eqref{eq: ma}.

\begin{lem}
The following equality holds:
		\begin{equation}\label{defkrvell}
ma\bigl(\lie{krv}_{ell}\bigr)=\Delta\bigl(ARI(\mathcal F_\Lau)^{\Delta}_{al+push*circneut}\bigr).
	\end{equation}
\end{lem}

\begin{proof}
By Definition \ref{mouldkrvell},
$b\in \lie{lie}_C$ is in $\lie{krv}_{ell}$ 
if and only if
$B_\ast=\Delta^{-1}\circ ma(b)$ is in 
$ARI(\mathcal F_\Lau)_{al+push*circneut}$.
Our claim follows because we have
$ma(\lie{lie}_C)=ARI(\mathcal F_\ser)_{al}$
and
$\Delta (ARI(\mathcal F_\ser)_{al})=ARI(\mathcal F_\Lau)^\Delta_{al}$.
\end{proof}

The operator $\Delta$ trivially respects push-invariance of moulds,
so the space $\lie{krv}_{ell}$ lies in the space 
		$\lie{lie}_C^{push}$ of push-invariant elements of $\lie{lie}_C$
        (cf. Definition \ref{defn:push}).  
		We will now show that the subspace $\lie{krv}_{ell}$ is actually a Lie
		subalgebra of $\lie{lie}_C^{push}$, which is itself a Lie algebra under a new Lie bracket in Corollary \ref{liepush} below,
	 of which a more explicit version (with a formula for 
				the partner) is proved in \cite{S3} (Lemma 2.1.1).

		\begin{lem}\label{partner} 
		Let $b\in \lie{lie}_C$.  Then $b\in\lie{lie}_C^{push}$ if and only if
		there exists a unique element $a\in \lie{lie}_C$ (the partner of $b$), 
		such that if $D_{b,a}$ is the derivation of $\lie{lie}_2$ defined by 
		$x\mapsto b$, $y\mapsto a$ (as introduced in \S \ref{subsec:Special types of derivations}), then $D_{b,a}$ annihilates $[x,y]$.
		\end{lem}

		\vspace{.1cm}
		By identifying $\lie{lie}_C^{push}$ with the space 
        $\lie{oder}_2$ (see \S \ref{subsec:Special types of derivations}) of derivations that
		annihilate $[x,y]$, we see that $\lie{lie}_C^{push}$ is
		a Lie algebra under the bracket of derivations.  We state this as a corollary.

		\begin{cor}\label{liepush} The map $b\mapsto D_{b,a}$ gives an isomorphism of linear spaces
		\begin{equation}\label{partial}
		i_o:\lie{lie}_C^{push}\overset{\sim}\rightarrow \lie{oder}_2
		\end{equation}
		whose inverse is $D_{b,a}\mapsto D_{b,a}(x)=b$, and this becomes
		a Lie algebra isomorphism when $\lie{lie}_C^{push}$ is equipped with the
		Lie bracket
		\begin{equation}\label{pushbracket}
		\langle b,b'\rangle = [D_{b,a},D_{b',a'}](x)=D_{b,a}(b')-D_{b',a'}(b).
		\end{equation}
		\end{cor}

 We note that the above derivation $\langle ,\rangle$ is different from the derivation
$\{,\}$ in \eqref{Poisson}.
We compare the $Dari$-bracket to the bracket
$\langle\,,\,\rangle$ on $\lie{lie}_C^{push}$.

	\begin{prop}\label{compat}
    The isomorphism $\mathrm{ma}$ in \eqref{eq: ma lie C push} is also compatible with the Lie brackets $\langle \cdot, \cdot \rangle$ and ${Dari}$ in the sense that
    \begin{equation}\label{eq:bracket-compat}
    ma\bigl(\langle b,b'\rangle\bigr)=Dari\bigl(ma(b),ma(b')\bigr).
     \end{equation}
    Consequently $ARI(\mathcal F_\ser)_{al+push}$ forms a Lie algebra under $Dari$-bracket.
    Moreover, we have a Lie algebra isomorphism
	$$ma:(\lie{lie}_C^{push},\langle \ ,\ \rangle)\rightarrow (ARI(\mathcal F_\ser)_{al+push},{Dari}).$$
 	\end{prop}

\begin{proof}
Since we know that the map $ma$ in \eqref{eq: ma lie C push} is bijective,
it is enough to prove that it is a Lie algebra homomorphism \eqref{eq:bracket-compat}.
The key point is the following non-trivial result,
	which is one of the main results of \cite{BS}:
	if $D_1$ and $D_2$ lie in $\lie{oder}_2$, then
	the map
	\begin{align}
	\Psi:\lie{oder}_2&\rightarrow ARI(\mathcal F_\Lau)_{ari}\label{psimap}  \\
		D&\mapsto \Delta^{-1}\bigl(ma\bigl(D(x)\bigr)\bigr), \notag
        \end{align}
	is an injective Lie algebra homomorphism, i.e.
	\begin{equation*}
\Delta^{-1}\Bigl(ma\bigl([D_1,D_2](x)\bigr)\Bigr)=ari\Bigl(\Delta^{-1}\bigl(ma(D_1(x))\bigr),\Delta^{-1}\bigl(ma(D_2(x))\bigr)\Bigr)
	\end{equation*}
	(see Theorem 3.5 of \cite{BS}).
	Applying $\Delta$ to both sides of this and using \eqref{Daribrack}, this
	is equivalent to
	\begin{equation}
ma\bigl([D_1,D_2](x)\bigr)=Dari\Bigl(ma\bigl(D_1(x)\bigr),ma\bigl(D_2(x)\bigr)
		\Bigr).
	\end{equation}
	By Lemma \ref{lem: ARI Fser Dari}, we see that 
	\begin{equation}\label{BS}
\Delta\circ \Psi:\lie{oder}_2\rightarrow ARI(\mathcal F_\ser)_{Dari}
\end{equation}
is a Lie algebra homomorphism.
We saw in Corollary \ref{liepush} that we have a Lie isomorphism
$i_o: (\lie{lie}_C^{push},\langle,\rangle)\buildrel\sim\over\rightarrow\lie{oder}_2$
when $\lie{lie}_C^{push}$ is equipped with the Lie bracket
(\ref{pushbracket}), so by composition, we have an injective Lie algebra homomorphism
$$b\buildrel i_o\over\mapsto D_{b,a}\buildrel\Psi\over\mapsto \Delta^{-1}\bigl(ma(D_{b,a}(x))\bigr)
	\buildrel\Delta\over\mapsto ma(b)$$
	(where $i_o$ is as in Corollary \ref{liepush} and $\Psi$ is as in 
	 \eqref{psimap}) is an injective Lie algebra homomorphism $\lie{lie}_C^{push}\rightarrow 
	ARI(\mathcal F_\ser)_{Dari}$, which proves the result.
\end{proof} 

It looks remarkable that   $\lie{lie}_C^{push}$ 
(and whence $ARI(\mathcal F_\ser)_{al+push}$)
is encoded with
two Lie algebraic structures:

\begin{rem}
(i).
  From  Corollary \ref{cor: Lie C push closed under the bracket} and
  Proposition \ref{compat},
  we deduce that 
  $\lie{lie}_C^{push}$ forms a Lie algebra with respect to both  the brackets $\{,\}$ and $\langle,\rangle$
  and $ARI(\mathcal F_\ser)_{al+push}$ inherits a Lie algebraic structure under both the
 $ari$ and $Dari$-brackets. 
 Furthermore, there exist sequences of Lie algebra maps:
 \begin{align*}
& (ARI({\mathcal F}_\ser)_{al+push},ari)\overset{ma}{\simeq}
 (\lie{lie}_C^{push},\{\ ,\ \})\overset{i_y}{\hookrightarrow}
\lie{tder}_2^{(y)}, \\
& (ARI({\mathcal F}_\ser)_{al+push},Dari)\overset{ma}{\simeq}
 (\lie{lie}_C^{push},\langle \ ,\ \rangle)\overset{i_o}{\simeq}
\lie{oder}_2.
 \end{align*}

 (ii).
{By combining the $\Q$-linear (but not Lie algebraic) isomorphism $i_{o,z}$
in \eqref{eq: oder and sder}, we obtain another identification of $\Q$-linear spaces
\begin{equation*}\label{eq: lie push and sder}
i_z 
: \quad  \lie{lie}_C^{push} \simeq\lie{sder}_2^{(z)}
\end{equation*}
which is  also obtained from Proposition \ref{specialprop}.
We note that the partner $a$ in Proposition \ref{specialprop} corresponds to $-a$ in Lemma \ref{partner}. 
Since one can show that 
$$i_z(b)\circ \nu=i_y(b)+\mathrm{ad}(-b)$$
where $\nu$ is the involution defined in \eqref{nu},
it follows that the Lie bracket on $\lie{sder}_2^{(z)}$ given by \eqref{eabbracket},
which also induces a new Lie algebra structure 
$\lie{lie}_C^{push}$,
is nevertheless related to the bracket
$\{\ ,\ \}$.
}
 \end{rem}

\subsubsection{Lie algebra structure on $\lie{krv}_{ell}$ (Theorem \ref{krvellisLie})}\label{krvellLie}         
The Lie algebra $(\lie{lie}_C^{push},\langle,\rangle)$  contains
the elliptic Kashiwara-Vergne space $\lie{krv}_{ell}$ as a linear subspace. 
This leads to our first main theorem concerning $\lie{krv}_{ell}$ (Theorem \ref{krvellliealg}), from which Theorem \ref{krvellisLie} is deduced.
As a corollary (Corollary \ref{lastone}), we obtain the inclusion
$\lie{lkrv}\rightarrow \lie{krv}_{ell}$
which establishes Proposition \ref{liemapprop}.

Before stating them, let us recall here the similar definition from of the elliptic double shuffle Lie algebra $\lie{ds}_{ell}$ from \cite{S3}, whose construction was the inspiration for the definition of the elliptic Kashiwara-Vergne Lie algebra.

\begin{defn}\label{defndsell} 
{\it The mould-version elliptic double shuffle Lie algebra}
was defined in \cite{S3} as the linear subspace 
$$
\Delta\bigl(ARI(\mathcal F_\Lau)^\Delta_{\underline{al}*\underline{al}}\bigr)
\subset ARI(\mathcal F_\ser)
.$$
{\it The elliptic double shuffle Lie algebra}
{\rm is the subspace $\lie{ds}_{ell}\subset \lie{lie}_C$ such that
$$ma(\lie{ds}_{ell})=\Delta\bigl(ARI(\mathcal F_\Lau)^\Delta_{\underline{al}*\underline{al}}\bigr).$$
It was also shown in \cite{S3} that 
$$\lie{ds}_{ell}\subset \lie{lie}^{push}_C,$$
and that $\lie{ds}_{ell}$ is closed under the Lie bracket 
$\langle \ , \ \rangle$ of \eqref{pushbracket}.
}
\end{defn}

\vspace{.3cm}
By Corollary \ref{mouldfirstcor}, we know that $ARI(\mathcal F_\Lau)_{\underline{al}*\underline{al}}\subset ARI(\mathcal F_\Lau)_{al+push*circneut}$, so we obtain the inclusions
$$ARI(\mathcal F_\Lau)^\Delta_{\underline{al}*\underline{al}}\subset ARI(\mathcal F_\Lau)^\Delta_{al+push*circneut}$$
and thus
$$\Delta\bigl(ARI(\mathcal F_\Lau)^\Delta_{\underline{al}*\underline{al}}\bigr)\subset \Delta\bigl(ARI(\mathcal F_\Lau)^\Delta_{al+push*circneut}\bigr).$$
In view of the above definitions of the elliptic double shuffle 
(in Definition \ref{defndsell})
and Kashiwara-Vergne Lie algebras
(in Definitions \ref{mouldkrvell} and \ref{defnkrvell}),
this gives us the linear space inclusion
\begin{equation}\label{dsellinkrvell}
\lie{ds}_{ell}\hookrightarrow \lie{krv}_{ell}: \qquad b\mapsto b
\end{equation}
announced in the introduction (see \eqref{commdiag0}).
We will show in the following theorem that $\lie{krv}_{ell}$ is also a Lie algebra. 

\begin{thm}\label{krvellliealg} 
The linear space  $\lie{krv}_{ell}$ forms a Lie algebra
under the Lie bracket $\langle\ ,\ \rangle$ in \eqref{pushbracket}.
And the following is a sequence of Lie algebra inclusions:
$$\lie{ds}_{ell}\subset \lie{krv}_{ell}\subset
		\lie{lie}_C^{push}.$$
		\end{thm}

\begin{proof}
    In view of the fact that $\lie{ds}_{ell}$ is known to be a Lie subalgebra of the Lie algebra $(\lie{lie}_C^{push}, \langle \ , \ \rangle)$, 
it remains only to show that 
    the subspace
	$\lie{krv}_{ell}\subset \lie{lie}_C^{push}$ is closed under the
	bracket $\langle \ ,\ \rangle$.

                \vspace{.2cm}
\noindent {\bf Step 1.}
By Proposition \ref{compat}, the map $ma$ gives an injective Lie algebra 
	morphism
	$$(\lie{lie}_C^{push},\langle\ ,\ \rangle)\rightarrow ARI(\mathcal F_\ser)_{Dari}.$$
	Thus proving  that $\lie{krv}_{ell}$ is closed under $\langle,\rangle$
    is equivalent to proving that its image $\Delta\bigl(ARI(\mathcal F_\Lau)^{\Delta}_{al+push*circneut}\bigr)$ under the map $ma$
    is closed under the
	$Dari$-bracket.  Since we saw above that
	$$\Delta^{-1}:ARI(\mathcal F_\Lau)_{Dari}\rightarrow ARI(\mathcal F_\Lau)_{ari}$$
    is a Lie algebra homomorphism
	it is equivalent to show that $ARI(\mathcal F_\Lau)^\Delta_{al+push*circneut}$
	is a Lie subalgebra of $ARI(\mathcal F_\Lau)_{ari}$.

				\vspace{.2cm}
				\noindent {\bf Step 2.} The space $ARI(\mathcal F_\Lau)^\Delta_{al+push}$ is a Lie algebra
				under $ari$.  Indeed, the definition of $\Delta$ shows that this operator 
				does not change the properties of push-invariance
				or alternality, i.e.~$\Delta^{-1}(ARI(\mathcal F_\Lau)_{al+push})=ARI(\mathcal F_\Lau)_{al+push}$. Restricted
				to $\mathcal F_\ser$-valued moulds, we have $\Delta^{-1}(ARI(\mathcal F_\ser)_{al+push})=
				ARI(\mathcal F_\Lau)^\Delta_{al+push}$.  Since $\Delta$ is an isomorphism from $ARI(\mathcal F_\Lau)_{ari}$ to 
				$ARI(\mathcal F_\Lau)_{Dari}$ by virtue of (\ref{Deltaiso}) and $ARI(\mathcal F_\ser)_{al+push}$ is
				a Lie subalgebra of $ARI(\mathcal F_\ser)_{Dari}$ by Proposition \ref{compat}, its image 
				$ARI(\mathcal F_\Lau)^\Delta_{al+push}$ under $\Delta^{-1}$ is thus a Lie subalgebra of
				$ARI(\mathcal F_\Lau)_{ari}$.

				\vspace{.3cm}
				\noindent {\bf Step 3.} We now complete the proof of Theorem 
				\ref{krvellliealg} by showing that the space $ARI(\mathcal F_\Lau)^\Delta_{al+push*circneut}$ is
				a Lie algebra under $ari$.  

            Let $A,B$ lie in 
				$ARI(\mathcal F_\Lau)^\Delta_{al+push*circneut}$, and let us show that $ari(A,B)$ lies
				in the same space.  By Step 2, we know that 
				$ari(A,B)\in ARI(\mathcal F_\Lau)^\Delta_{al+push}$, so we only need to show that
				$swap\bigl(ari(A,B)\bigr)$ is $*$circ-neutral\footnote{We recall that $*$circ-neutral means circ-neutral up to addition of a constant mould, 
                (see footnote\ \ref{foot: conventionARIsymbol}).}.  But we will show that in fact 
				this mould is actually circ-neutral.  To see this, let $A_0$ and $B_0$ be
				the constant-valued moulds such that $swap(A)+A_0$ and $swap(B)+B_0$ are
				circ-neutral.  
                Since $A+A_0$ and $B+B_0$ are push-invariant,
                by Proposition \ref{circneutclosed}, we have
				$$\overline{ari}\bigl(swap(A)+A_0,swap(B)+B_0\bigr)\in \overline{ARI}(\mathcal F_\Lau)_{circneut}.$$
				Using the identity $$swap\bigl(ari(M,N)\bigr)=\overline{ari}\bigl(swap(M),
						swap(N)\bigr),$$ 
                        valid whenever $M$ and $N$ are push-invariant moulds
				(cf.~[S], (2.5.6)), as well as the fact that constant-valued moulds are
				both push and swap invariant, we have
				\begin{align*}
				\overline{ari}\bigl(&swap(A)+A_0,swap(B)+B_0\bigr)=
				\overline{ari}\bigl(swap(A+A_0),swap(B+B_0)\bigr)\notag\\
					&=swap\cdot ari(A+A_0,B+B_0)\notag\\
					&=swap\cdot ari(A,B)+swap\cdot ari(A,B_0)+swap\cdot ari(A_0,B)+swap\cdot
					ari(A_0,B_0)\\
						&=swap\cdot ari(A,B)\notag
						\end{align*}
						since the definition of the $ari$-bracket shows that $ari(C,M)=0$ whenever
						$C$ is a constant-valued mould.  So $swap\cdot ari(A,B)$ is circ-neutral, which follows
                        $ari(A,B)\in ARI(\mathcal F_\Lau)^{\Delta}_{al+push*circneut}$.
						The proof of Theorem \ref{krvellliealg} is completed.  
\end{proof}

						\vspace{.2cm}

{\bf Proof of Theorem \ref{krvellisLie}}.
The claim (i) is immediate, since the conditions of push-invariance and circ-neutrality up to addition of constants are both homogeneous.

The claims (ii) and (iii) follow directly from Corollary \ref{liepush} and Theorem \ref{krvellliealg}.
\qed
       						\vspace{.2cm}

\begin{cor}\label{lastone}
The linear map $\Delta$ defined in \eqref{eq: map Delta} gives a Lie algebra morphism
\begin{equation}\label{Deltamap}
\Delta:(ARI({\mathcal F}_{\ser})_{al+push/circneut},ari)\rightarrow (\Delta(ARI(\mathcal F_\Lau)^\Delta_{al+push*circneut}),Dari),
\end{equation}
which induces a Lie algebra morphism
\begin{align}\label{liemap} 
(\lie{lkrv},\{\ , \ \})&\hookrightarrow (\lie{krv}_{ell},\langle \ , \ \rangle)\notag\\
b(x,y)&\mapsto [x,b(x,[x,y])].
\end{align}
\end{cor}

\begin{proof}
For the first statement, composing the inclusion map 
\begin{align}\label{incl2}
ARI({\mathcal F}_{\ser})_{al+push/circneut}
&=ARI({\mathcal F}_{\ser})\cap ARI({\mathcal F}_{\Lau})^\Delta_{al+push/circneut} \\ \notag
&\qquad\qquad \subset ARI({\mathcal F}_{\Lau})^\Delta_{al+push*circneut}
\end{align}
with the operator $\Delta$, considered as an injective 
linear map on moulds gives an injective linear map 
$$
\Delta:
ARI({\mathcal F}_{\ser})_{al+push/circneut}\hookrightarrow \Delta(ARI({\mathcal F}_{\Lau})^\Delta_{al+push*circneut}).$$
It is shown in Step 3 of the proof of Theorem \ref{krvellliealg}  that 
$ARI({\mathcal F}_{\Lau})^\Delta_{al+push*circneut}$ is a Lie algebra under the $ari$-bracket, and in 
Proposition \ref{prop: ARI al push circneut Lie algebra}
that $ARI({\mathcal F}_{\ser})_{al+push/circneut}$ is a Lie subalgebra
of it. 
A basic property of the linear map $\Delta$ is that it transforms 
the $ari$-bracket into the $Dari$-bracket 
(cf. \eqref{Daribrack}), 
so the space 
$\Delta(ARI({\mathcal F}_{\Lau})^\Delta_{al+push*circneut})$ 
is a Lie algebra under the $Dari$-bracket.
Thus the map in \eqref{Deltamap} is a Lie algebra morphism from a Lie subalgebra of $ARI({\mathcal F}_{\ser})_{ari}$ to a Lie subalgebra of $ARI({\mathcal F}_{\Lau})_{Dari}$.

  Finally, by \eqref{lkrvARI} we have 
$$ma(\lie{lkrv})=ARI({\mathcal F}_{\ser})_{al+push/circneut}\subset ARI({\mathcal F}_{\ser})_{ari}$$
and by \eqref{defkrvell} we have
$$ma(\lie{krv}_{ell})=\Delta(ARI({\mathcal F}_{\Lau})^{\Delta}_{al+push*circneut})\subset ARI({\mathcal F}_{\ser})_{Dari},$$
so \eqref{Deltamap} translates directly under $ma^{-1}$ to a Lie algebra homomorphism $\lie{lkrv}\rightarrow \lie{krv}_{ell}$. 
The presentation \eqref{liemap} follows from \eqref{eq: map Delta}   and \eqref{eq: ma}. 
\end{proof}

Whence Proposition \ref{liemapprop} is proved. \qed
\vspace{.2cm}

\subsection{The map
$\lie{krv}\hookrightarrow\lie{krv}_{ell}$
(Theorem \ref{krvsection})}
 \label{subsec:the map from krv to krvell}

 In this subsection we prove our next main result on the elliptic
Kashiwara-Vergne Lie algebra, which is analogous to known results on 
the elliptic Grothendieck-Teichm\"uller Lie algebra of \cite{E1} and the
elliptic double shuffle Lie algebra of \cite{S3}. 
\S \ref{doubleshuf} below is devoted to connections between these three situations.

To state the assumption of our result, we prepare  some definitions on moulds.
\begin{defn}\label{defsenary} 
\noindent (i) Let 
$$teru:ARI(\mathcal F)\to ARI(\mathcal F)$$ be the operator defined for $A\in ARI(\mathcal F)$ as follows:
$teru(A)$ is equal to $A$ in depths $0$ and $1$, and for depths $r>1$, we have
						\begin{equation}\label{teru}
						teru(A)(u_1,\ldots,u_r)=
						\qquad \qquad \qquad \qquad \qquad \qquad \qquad \qquad
						\qquad \qquad \qquad 
						\end{equation}
						\begin{equation*}
						\quad \qquad A(u_1,\ldots,u_r)+
						\frac{1}{u_r}
						\Bigl(A(u_1,\ldots,u_{r-2},u_{r-1}+u_r)-A(u_1,\ldots,u_{r-2},u_{r-1})\Bigr).
						\end{equation*}
                        
(ii) A mould $A\in ARI(\mathcal F)$ is said to satisfy {\it the senary relation} (cf.~(3.64)
		in \S 3.5 of \cite{Ec}) if
	\begin{equation}\label{senaryrelation}
	teru(A)=push\circ mantar\circ teru\circ mantar(A),
	\end{equation}
	and the {\it twisted senary relation} if
    \begin{equation}\label{senarybis}
	teru\circ pari(B)=push\circ mantar\circ teru\circ pari(B).
	\end{equation}
    (for $pari$  see Definition \ref{def: various mould operators}),
    that means that  $pari(A)$ satisfies the senary  relation.
    
(iii) We define the mould subspace
	\begin{equation}\label{eq: ARIal+tsen/circconst}
	        ARI(\mathcal F)_{al+tsen/circconst} \qquad
	(\text{resp.} \ ARI(\mathcal F)_{al+tsen*circconst})
    	\end{equation}
    to be the subspace of 
	 moulds $A\in ARI(\mathcal F)$ such that $swap(A)$ is circ-constant  (see Definition \ref{circconstance})
	(resp.~up to adding a constant mould) and $A$ satisfies the twisted senary relation \eqref{senarybis}.
    \end{defn}

\begin{rem}\label{rem:constant c/n}
Observe that if $swap(A)\in\overline{ARI}(\mathcal F)$ 
for ${\mathcal F}={\mathcal F}_\ser$ or ${\mathcal F}_\pol$
is a polynomial of homogeneous degree $n$ which is circ-constant up to addition of
	a constant-valued mould, then the constant-valued mould is uniquely 
	determined as being the mould whose only non-zero value is the constant
	value $\frac{c}{n}$ in depth $n$, where $c$ is given by 
	$$swap(A)(v_1)=cv_1^{n-1}.$$
    \end{rem}

\vspace{.2cm}
\noindent {\bf Theorem \ref{krvsection}.} {\it 
Assume 
\footnote{{This assertion is announced to hold in \cite{Ka}.}}
that the adjoint action $Ad_{ari}(pal):ARI(\mathcal F_\Lau)_{al}\to ARI(\mathcal F_\Lau)_{al}$
of the mould $pal\in GARI(\mathcal F_\Lau)$ (cf. Definition \ref{defn:pal and pil})
restricts to a bijection between
 \begin{equation}\label{Ecstatement}
 \xymatrix{Ad_{ari}(pal):ARI(\mathcal F_\Lau)_{push}\ar@{->}[r]^{\quad\quad \sim}&ARI(\mathcal F_\Lau)_{sen},}
 \end{equation}
 where $ARI(\mathcal F_\Lau)_{push}$ denotes the set of  $push$-invariant moulds \eqref{pushinv}, and $ARI(\mathcal F_\Lau)_{sen}$
denotes the set of moulds satisfying the senary relation \eqref{senaryrelation}.
 Then there exists an injective Lie algebra morphism
\begin{equation}\label{krvtokrvell}
\lie{krv}\hookrightarrow\lie{krv}_{ell}.
\end{equation}
}

The proof constructs the morphism from $\lie{krv}$ to $\lie{krv}_{ell}$
in four main steps as follows. 

\vspace{.2cm}
\noindent {\bf Step 1.} We first consider a twisted version
of the Kashiwara-Vergne Lie algebra, or rather of the associated
polynomial space $V_{\lie{krv}}$ of Definition \ref{Vkrv}, via the map
\begin{align}\label{VtoW}
\nu:V_{\lie{krv}}&\buildrel\sim\over\rightarrow W_{\lie{krv}}\\
f\ &\mapsto \nu(f), \nonumber
\end{align}
where $\nu$ is the automorphism of $Ass_2$ defined by \eqref{nu}.
In paragraph \ref{Wkrv}, we prove that $W_{\lie{krv}}$ is a Lie algebra under
the Poisson or Ihara bracket, and give a description of 
$W_{\lie{krv}}$ via two properties, the ``twisted'' versions of the two
defining properties of $V_{\lie{krv}}$ given in Definition \ref{Vkrv}.

\vspace{.2cm}
\noindent {\bf Step 2.} 
	In paragraph \ref{mouldWkrv}, we study the mould space 
	$ma\bigl(W_{\lie{krv}}\bigr)$.  Thanks to the compatibility
	of the $ari$-bracket with the Poisson bracket (see \eqref{ariPoisson}),
	this space is a Lie subalgebra of $ARI({\mathcal F}_\ser)_{ari}$.  Just as we reformulated
	the defining properties of $\lie{lkrv}$ in mould terms in \S \ref{sec: mould theory}, proving that
	$ma(\lie{lkrv})=ARI({\mathcal F}_\ser)_{al+push/circneut}$, in \eqref{lkrvARI}, we will reformulate the
	defining properties of $W_{\lie{krv}}$ in mould terms and show that
	\begin{equation}\label{Wkrvsen}
	ma\bigl(W_{\lie{krv}}\bigr)=ARI({\mathcal F}_\ser)_{al+tsen*circconst}.
	\end{equation}

	\vspace{.2cm}
	 \noindent {\bf Step 3.} 
In paragraph \ref{subsubsec: Step 3: Construction}, we consider the map 
$$\Xi:=Ad_{ari}(invpal)\circ pari:ARI({\mathcal F}_\Lau)_{ari}\rightarrow ARI({\mathcal F}_\Lau)_{ari},$$ 
and 
we show that it gives an injective Lie algebra homomorphism 
\begin{equation}\label{Xi}
\xymatrix{
ARI({\mathcal F}_\ser)_{al+tsen*circconst}\ar[r]^{pari}
&ARI({\mathcal F}_\ser)_{al+sen*circconst}\ar@{^{(}->}[rr]^{Ad_{ari}(invpal)}&&
ARI({\mathcal F}_\Lau)^\Delta_{al+push*circneut}}
\end{equation}
of subalgebras of $ARI({\mathcal F}_\ser)_{ari}$.

\vspace{.2cm}
\noindent {\bf Step 4.} The final step in the paragraph
\ref{subsubsec: Step 4: Comparison} is to compose (\ref{Xi}) with 
the Lie algebra homomorphism 
$\Delta:ARI({\mathcal F}_\Lau)_{ari}\rightarrow ARI({\mathcal F}_\Lau)_{Dari}$, obtaining
an injective Lie algebra homomorphism 
$$ARI({\mathcal F}_\ser)_{al+tsen*circconst} \rightarrow \Delta\bigl(ARI({\mathcal F}_\Lau)^\Delta_{al+push*circneut}\bigr),$$
where the left-hand space is a subalgebra of $ARI({\mathcal F}_\ser)_{ari}$ and the right-hand one of $ARI({\mathcal F}_\ser)_{Dari}$.
Since the right-hand space is equal to $ma(\lie{krv}_{ell})$,
the desired injective Lie algebra homomorphism from $\lie{krv}$ to $\lie{krv}_{ell}$
is obtained by composing all the maps described above, as shown in the 
following diagram: 
\begin{equation}\xymatrix{
\lie{krv}\ar[d]_{{\rm by}\ \eqref{krviso}}&\\
 V_{\lie{krv}}\ar[d]_{{\rm by}\ \eqref{VtoW}}^\nu &
 \lie{krv}_{ell}\\
W_{\lie{krv}}\ar[d]_{{\rm by}\ \eqref{Wkrvsen}}^{ma}& \Delta\bigl(ARI({\mathcal F}_\Lau)^\Delta_{al+push*circneut}\bigr)\ar[u]^{ma^{-1}}_{{\rm by}\
\eqref{defkrvell}}\\
ARI({\mathcal F}_\ser)_{al+tsen*circconst}\ar@{->}[r]^{\Xi}_{{\rm by}\ \eqref{Xi}}&
ARI({\mathcal F}_\Lau)^\Delta_{al+push*circneut}\ar[u]^{\Delta}
 }
\end{equation}

\vspace{.2cm}
\subsubsection{Step 1: The twisted space $W_{\lie{krv}}$}\label{Wkrv}

\begin{prop}\label{WkrvLie}
Let 
$$W_{\lie{krv}}=\nu(V_{\lie{krv}})$$
(for $V_{\lie{krv}}$, see Definition \ref{Vkrv}).
Then 
$W_{\lie{krv}}$ is a Lie algebra under the Poisson bracket
$\{ \ , \ \}$ in \eqref{Poisson}.
\end{prop}

\begin{proof}
We use Lemma \ref{isLie} to complete the proof of Proposition \ref{WkrvLie}.
Write 
		$$\lie{krv}^\nu=\{\nu\circ E\circ \nu\,|\,E\in \lie{krv}\}\subset
		\lie{sder}_{2}^{(x)}.$$
		By restricting the isomorphism (\ref{conjbynu}) to
		the subspace $\lie{krv}\subset \lie{sder}_2^{(z)}$, we obtain a commutative
		diagram of isomorphisms of linear spaces
		$$\xymatrix{
			\lie{krv}\ar[r]\ar[d]&\lie{krv}^\nu\ar[d]\\
				V_{\lie{krv}}\ar[r]^\nu&W_{\lie{krv}}
		}$$
where the left-hand vertical arrow is the isomorphism
(\ref{krviso}) mapping $E_{a,b}\mapsto b$, 
	and the right-hand vertical map sends an 
	Ihara derivation $d_f$ to $f$.  Equipping $W_{\lie{krv}}$ with the
	Lie bracket inherited from $\lie{krv}^\nu$ makes this
	into a commutative diagram of Lie isomorphisms. But this bracket is
	nothing other than the Poisson bracket $\{ \ , \ \}$ in \eqref{Poisson}
    since $\lie{krv}^\nu
	\subset\lie{sder}_{2}^{(x)}$.
\end{proof}

	We now give a characterization of $W_{\lie{krv}}$
	by two defining properties which are the twists by $\nu$  defined in \eqref{nu} of those defining
	$V_{\lie{krv}}$.  Recall that $\beta$ is the the backwards operator given in
	Definition \ref{backwards}.  

	\begin{prop}\label{Wtwoprops} The space $W_{\lie{krv}}$ is the space spanned by polynomials
	$b\in \lie{lie}_C$, of homogeneous degree $n\ge 3$, such that 

	\vspace{.2cm}
\noindent (i)\ $b_y-b_x$ is {\rm anti-palindromic}, i.e.~$\beta(b_y-b_x)
	=(-1)^{n-1}(b_y-b_x)$, and

	\vspace{.2cm}
	\noindent (ii)\ {$b+\frac{c}{n}y^n$  with $c=(b|x^{n-1}y)$
    is {\rm  circ-constant}
     for $\{c_k\}_{k>1}$ with $c_k=\delta_{k,n}\cdot (b|x^{n-1}y)$.}
	\end{prop}

	\vspace{.1cm}
	\textit{Proof.} {Let $b\in W_{\lie{krv}}$. Put $f=\nu(b)$, so we have $f\in V_{\lie{krv}}$.}  Then the 
	property that $b_y-b_x$ is anti-palindromic is precisely equivalent to 
	the push-invariance of $f$ (this is proved as the equivalence of properties 
			(iv) and (v) of Theorem 2.1 of \cite{S1}). This proves (i).

	For (ii), we note that since $f\in V_{\lie{krv}}$, $f^y-f^x$ is
	push-constant for the value $c=(f|x^{n-1}y)=(-1)^{n-1}(b|x^{n-1}y)$.
	We have
	$$b(x,y)=xb^x(x,y)+yb^y(x,y),$$
	so
	$$f(x,y)=b(z,y)=zb^x(z,y)+yb^y(z,y)=-xb^x(z,y)-yb^x(z,y)+yb^y(z,y).$$
	Thus since $f(x,y)=xf^x(x,y)+yf^y(x,y)$, this gives 
	$$f^x=-b^x(z,y)\ \ {\rm and}\ \ 
	f^y=-b^x(z,y)+b^y(z,y),$$
	so
	$$f^y-f^x=b^y(z,y)=\nu(b^y).$$
	Thus to prove the result, it suffices to prove that the following
	statement: if $g\in Ass_C$ is a polynomial of homogeneous degree $n$
	that is push-constant for $(-1)^{n-1}c$, then 
	$\nu(g)$ is circ-constant for $c$, since taking $g=f^y-f^x$ then shows that 
	$\nu(g)=b^y$ is circ-constant for $c$. The proof of this statement
	is straightforward using the 
	substitution $z=-x-y$ (but see the proof of Lemma 3.3 in \cite{S1} for details).
	To complete the proof of (ii), we note that when 
	$f\in V_{\lie{krv}}$ is of even degree $n$ we have $c=0$. In fact this 
	follows from Corollary \ref{firstcor}, which states that 
	$\lie{lkrv}_n^1=0$ when $n$ is even; this means that
	there are no elements in $\lie{krv}$ of even weight $n$ and depth 1,
	so there are no such elements in $V_{\lie{krv}}$. Since $c$ is the coefficient
	of the depth 1 term $x^{n-1}y$, we have $c=0$ when $n$ is even. 
	This completes the proof of (ii).  \ \qed

	\vspace{.2cm}
	\subsubsection{Step 2: The mould version $ma(W_{\lie{krv}})$}\label{mouldWkrv}
	The space $ma(W_{\lie{krv}})$ is closed under the $ari$-bracket by 
	(\ref{ariPoisson}), since $W_{\lie{krv}}$ is closed under the Poisson bracket by Proposition \ref{WkrvLie}.

	Let $b\in W_{\lie{krv}}$ and let $B=ma(b)$.  Then since $b$ is
	a Lie polynomial, $B$ is an alternal polynomial mould.  Let us give the
	mould reformulations of properties (i) and (ii) of Proposition \ref{Wtwoprops}.
	The second property is easy since we already showed, in Proposition 
	\ref{circneutprop}, that a polynomial $b$ is circ-constant if and only 
	if $swap(B)$ is circ-constant.

	Expressing the first property in terms of moulds is more complicated and
	calls for an identity discovered by \'Ecalle.  We need to use the
	mould operators $mantar$  and $pari$ in Definition \ref{def: various mould operators}.

	The operator $pari$ extends the operator $y\mapsto -y$ on polynomials to 
	all moulds, and $mantar$ extends the operator $f\mapsto (-1)^{n-1}\beta(f)$.

	\begin{lem}\label{terulem} Let $b\in\lie{lie}_C$. 
    {Assume $\deg b>2$.}
    Then the following are
	equivalent:

	\vspace{.2cm}
	(i) $b_y-b_x$ is anti-palindromic;

	\vspace{.1cm}
	(ii) if $B=ma(b)$, then $B$ satisfies the twisted senary relation
    \eqref{senarybis}.
	\end{lem}
	(Note that since $b$ is a Lie element, $B$ and $pari(B)$ are alternal
	 and thus $mantar$-invariant, so we can drop the right-hand $mantar$ from the 
	 senary relation \eqref{senaryrelation}.)

	\begin{proof} It suffices to prove the statement for an element
	$b$ of homogeneous degree $n$.
	The statement is a consequence of the following result, proved 
	in Proposition A.3 of the Appendix of \cite{S1}
    {(see also \cite[Lemma 4.2]{FK2}).}
	Let $\tilde b\in\lie{lie}_C$ and let $\tilde B=ma(\tilde b)$
    {with $n>2$.}
	Write $\tilde b=\tilde b_xx+\tilde b_yy$ as usual.
	Then for each depth $r$ part $(\tilde b_x+\tilde b_y)^r$ of the polynomial
	$\tilde b_x+\tilde b_y$ ($1\le r\le n-1$),
	the anti-palindromic property 
	\begin{equation}\label{kalinda}
	(\tilde b_x+\tilde b_y)^r=(-1)^{n-1}\beta(\tilde b_x+\tilde b_y)^r
	\end{equation}
	translates directly to the following relation on $\tilde B$: 
	\begin{equation}\label{senary}
	teru(\tilde B)(u_1,\ldots,u_r)=
	push\circ mantar\circ teru(\tilde B)(u_1,\ldots,u_r).
	\end{equation}
Let us deduce the equivalence of the claims (i) and (ii) from that of (\ref{kalinda})
	and (\ref{senary}).
	Let $\tilde b$ be defined by $\tilde b(x,y)=b(x,-y)$. 
	This implies that $(b_x)^r=(-1)^r(\tilde b_x)^r$, $(b_y)^r
	=(-1)^{r-1}(\tilde b_y)^r$,
	and $\tilde B=pari(B)$.
	Thus $b_y-b_x$ is anti-palindromic if and only if $\tilde b_y+\tilde b_x$ is, 
	i.e.~if and only if (\ref{kalinda}) holds for $\tilde b$, which is the case
	if and only if (\ref{senary}) holds for $\tilde B$, which is equivalent
	to (\ref{senarybis}) for $B$. 
    This proves the lemma.
    \end{proof}

	\begin{cor}\label{alsen} 
	We have the 
	isomorphism of Lie algebras 
    \begin{equation}\label{eq:ma Wkrv to ARIal+tsen*circconst}
	\xymatrix{ma:W_{\lie{krv}}\ar[r]^{\!\!\!\!\!\!\!\!\!\!\!\!\!\!\!\!\!\!\!\!\!\!\!\!\!\!\!\!\!\!\!\!\!\!\!\!\!\!\!\!\!\!\!\!\!\!\!\!\!\!\!\!\sim}&ARI(\mathcal F_\ser)_{al+tsen*circconst}
	\subset ARI(\mathcal F_\ser)_{ari}}.
    \end{equation}
	\end{cor}

\begin{proof}
        By Proposition \ref{Wtwoprops}, the space 
	$W_{\lie{krv}}$ is the space of
	Lie polynomials $b$ satisfying (i) $b_y-b_x$ is antipalindromic and
	(ii) $b+\frac{c}{n}y^n$ for $c=(b|x^{n-1}y)$
    is 
    {circ-constant for $\{c_k\}_{k>1}$ with $c_k=\delta_{k,n}\cdot (b|x^{n-1}y)$}.
	By Lemma \ref{terulem}, property (i) is equivalent to the fact that
	$pari(B)$ satisfies the senary relation \eqref{senaryrelation}. 
    By Proposition \ref{circneutprop} the fact that $b$ is circ-constant is
	equivalent to $swap(B)$ being circ-constant 
    (and Remark \ref{rem:constant c/n} 
            shows that the constant is necessarily unique
			and the same). But by Definition \ref{defsenary},
	$ARI(\mathcal F_\ser)_{al+tsen*circconst}$ is precisely the space of alternal polynomial
	moulds satisfying precisely these two mould properties.
\end{proof}

\vspace{.2cm}
\subsubsection{Step 3: Construction of the map $\Xi$}\label{subsubsec: Step 3: Construction}
In this section we finally arrive at the main step of the construction
of our map $\lie{krv}\rightarrow \lie{krv}_{ell}$, namely the construction
of the map $\Xi$ given in the following proposition.

					\begin{prop}\label{keypropb}
				The operator $\Xi=Ad_{ari}(invpal)\circ pari$ gives an injective Lie algebra homomorphism 
					of Lie subalgebras of $ARI({\mathcal F}_\Lau)_{ari}$:
\begin{equation}\label{themap}
\Xi:ARI({\mathcal F}_\ser)_{al+tsen*circconst}\hookrightarrow ARI({\mathcal F}_\Lau)^\Delta_{al+push*circneut}.
\end{equation}
				\end{prop}

\textit{Proof.}
We have already shown that
both spaces are Lie subalgebras of $ARI({\mathcal F}_\Lau)_{ari}$, the first in 
Corollary \ref{alsen} and the second in \S \ref{krvellLie}.  Furthermore, since $pari$ and $Ad_{ari}(invpal)$ are both invertible and respect the $ari$-bracket, the proposed map is a Lie algebra map invertible on its image, 
and therefore injective.
Thus it remains only to show that the image of $ARI({\mathcal F}_\ser)_{al+tsen*circconst}$
under $\Xi$ really lies in the target space
$ARI({\mathcal F}_\Lau)^\Delta_{al+push*circneut}$.
We will show separately that if
$B\in ARI({\mathcal F}_\ser)_{al+tsen*circconst}$ and $A=\Xi(B)$, then

\vspace{.2cm}
(i) $A$ is push-invariant,

\vspace{.1cm}
(ii) $A$ is alternal,

\vspace{.1cm}
(iii) $swap(A)$ is $\ast$circ-neutral (i.e.~circ-neutral up to addition of a constant-valued mould),

\vspace{.1cm}
(iv) $A\in ARI({\mathcal F}_\Lau)^\Delta$.

\vspace{.3cm}
\noindent \underline{\it Proof of (i): $A$ is push-invariant.} This statement follows directly on \'Ecalle's
senary property \eqref{senaryrelation}. 
Indeed, since $B$ satisfies (\ref{senarybis}), 
$\tilde B:=pari(B)$ satisfies (\ref{senary}), so \eqref{Ecstatement} implies that
$Ad_{ari}(invpal)(\tilde B)=\Xi(B)=A$ is push-invariant.
\qed

\vspace{.3cm}
\noindent \underline{\it Proof of (ii): $A$ is alternal.} The subspace of alternal moulds 
$ARI({\mathcal F}_\Lau)_{al}$ is closed under $ari$ 
(cf.~\cite{SS}), so $exp_{ari}(ARI({\mathcal F}_\Lau)_{al})$ forms a subgroup of $GARI({\mathcal F}_\Lau)_{gari}$,
which we denote by $GARI({\mathcal F}_\Lau)^{as}_{gari}$ (the superscript $as$ stands
for {\it symmetral}).  
The mould $pal$ is known to be symmetral 
(cf.~\cite{Ec2}, or in more detail \cite{S2}, Theorem 4.3.4).
Thus, since $GARI({\mathcal F}_\Lau)^{as}_{gari}$ is a group, the $gari$-inverse mould $invpal$ is also symmetral.
Therefore the adjoint action $Ad_{ari}(invpal)$ on $ARI({\mathcal F}_\Lau)$ 
restricts to an adjoint action on the Lie subalgebra $ARI({\mathcal F}_\Lau)_{al}$ of alternal moulds. 
If $B$ is alternal, then $pari(B)$ is alternal, and so 
$A=\Xi(B)$ is alternal.
This completes the proof of (ii).
\qed

\vspace{.3cm}
For the assertions (iii) and (iv), we will make use of \'Ecalle's 
fundamental identity in the version (\ref{goodfund}) given in \S \ref{adariinvpal},
with $N=pari(B)$ (recall that (\ref{goodfund}) is valid whenever
$Ad_{ari}(invpal)\cdot N$ is push-invariant, which is the case for $pari(B)$ thanks to (i) above).  
The key point is that the operators $\overline{ganit}(poc)$ and $Ad_{\overline{ari}}(invpil)$
on the left-hand side of (\ref{goodfund}) are better adapted to tracking the circ-neutrality and the denominators than 
the right-hand operator $Ad_{ari}(invpal)$ considered directly.

\vspace{.3cm}
\noindent \underline{\it Proof of (iii): $swap(A)$ is $\ast$circ-neutral.} 
Let $b\in W_{\lie{krv}}$, and assume that $b$ is of homogeneous degree $n$.
Let $B=ma(b)$.  Then by Corollary \ref{alsen},
$swap(B)$ is circ-constant, and even circ-neutral
if $n$ is even.  

We need to show that 
$swap\cdot \Xi(B)=swap\cdot Ad_{ari}(invpal)\cdot pari(B)$ is 
$\ast$circ-neutral.
To do this, we use \eqref{goodfund} with $N=pari(B)$, 
and in fact show the result on the left-hand side, which is equal to
$$Ad_{\overline{ari}}(invpil)\cdot \overline{ganit}(poc)\cdot pari\cdot swap(B)$$
(noting that $pari$ commutes with $swap$).
We prove that this mould is
$\ast$circ-neutral in three steps:
\begin{itemize}
\item First we show that the operator
$\overline{ganit}(poc)\cdot pari$ changes a circ-constant mould 
into one that is circ-neutral (Proposition \ref{ganitcirc}). 
\item Secondly, we show that the operator $Ad_{\overline{ari}}(invpil)$ preserves 
the property of circ-neutrality (Proposition \ref{adaripilcirc}).
\item Finally, we show that if $M$ is a mould that is not circ-constant but only $\ast$circ-constant, and if $M_0$ is the (unique) constant-valued mould such that $M+M_0$ is circ-constant, then
$Ad_{\overline{ari}}(invpil)\cdot \overline{ganit}(poc)\cdot pari(M)+M_0$ is circ-neutral, 
which says that
$Ad_{\overline{ari}}(invpil)\cdot \overline{ganit}(poc)\cdot pari(M)$ is $\ast$circ-neutral
(a paragraph after Proposition \ref{adaripilcirc}).
\end{itemize}

\begin{prop}\label{ganitcirc} Fix $n\ge 1$, and let $M\in \overline{ARI}({\mathcal F}_\pol)$ 
be a circ-constant (cf.~Definition \ref{circconstance})
polynomial-valued mould of homogeneous degree $n$.  
Then $\overline{ganit}(poc)\cdot pari(M)$ is circ-neutral.  
\end{prop}

\vspace{.2cm} When $n=1$, the homogeneous weight
$n=1$ means that $M$ is concentrated in depth 1 where it has constant value
$M(v_1)=c$.  Thus the condition of being circ-constant is automatically 
satisfied. Direct computation shows that 
$$\overline{ganit}(poc)\cdot pari(M)(v_1,\ldots,v_r)=
\frac{-c}{(v_1-v_2)(v_2-v_3)\cdots (v_{r-1}-v_r)}.$$
Summing the images of this mould under powers of the circ-operator and
putting them over the common denominator $(v_r-v_1)(v_1-v_2)\cdots (v_{r-1}-v_r)$ shows that this mould is circ-neutral. For the remainder of the proof, we
assume that $n>1$.

\vspace{.3cm}
\textit{Notation for the proof of Proposition \ref{ganitcirc}.} 
Let ${\bf v}=(v_1,\ldots,v_r)$, and let ${\bf W}_{\bf v}$ be the set of 
decompositions $d_{\bf v}$ of ${\bf v}$ into chunks $d_{\bf v}=
{\bf a}_1{\bf b}_1\cdots {\bf a}_s{\bf b}_s$ as in (\ref{chunks}).  
We always denote consecutive chunks ${\bf b}_i$ in the form of tuples
$(v_k,v_{k+1},\ldots,v_{k+l})$, but they can also be considered simply
as subsets of $\{v_1,\ldots,v_r\}$.

For any decomposition $d_{\bf v}$, we let its ${\bf b}$-part be the subset 
${\bf b}_1\cup \cdots\cup {\bf b}_s$ of $\{v_1,\ldots,v_r\}$ and 
its ${\bf a}$-part be the subset ${\bf a}_1\cup\cdots\cup{\bf a}_s$; the 
${\bf a}$-part is equal to the complement of the ${\bf b}$-part 
in $\{v_1,\ldots,v_r\}$. We write $|{\bf a}|$ for the number of letters 
in the ${\bf a}$-part, i.e. $|{\bf a}|=|{\bf a}_1|+\cdots+|{\bf a}_s|$. 

Set 
\begin{equation}\label{thesetW} {\bf W}=\coprod_{1\le i\le r} 
{\bf W}_{\sigma_r^i({\bf v})}, 
\end{equation} 
where the $\sigma_r^i({\bf v})$ are the cyclic permutations 
of ${\bf v}=(v_1,\ldots,v_r)$. For a fixed subset ${\bf b}\subsetneq 
\{v_1,\ldots,v_r\}$, let ${\bf W}^{\bf b}$ denote the subset of decompositions 
in ${\bf W}$ having ${\bf b}$-part equal to ${\bf b}$; thus we have 
\begin{equation}\label{Wb} 
{\bf W}=\coprod_{{\bf b}\subsetneq\{v_1,\ldots,v_r\}} {\bf W}^{\bf b}.  
\end{equation} 
Let ${\bf w}=(v_{i+1},\ldots,v_r,v_1,\ldots,v_i)$ be any cyclic permutation 
of ${\bf v}=(v_1,\ldots,v_r)$; let $d_{\bf w}={\bf a_1b_1\cdots a_sb_s}$ 
be a decomposition of ${\bf w}$, and let 
${\bf b}={\bf b_1}\cup \cdots\cup {\bf b_s}$ be its ${\bf b}$-part. We will 
list all the elements of ${\bf W}^{\bf b}$, i.e.~all decompositions of all 
cyclic permutations of ${\bf v}$ having ${\bf b}$-part equal to ${\bf b}$. 
Let ${\bf a}={\bf a_1}\cup \cdots\cup {\bf a_s}$ be the ${\bf a}$-part of 
$d_{\bf w}$. Then there exists a decomposition of a cyclic permutation of 
${\bf v}$ having ${\bf b}$-part equal to ${\bf b}$ if and only if the cyclic 
permutation begins with a letter $v_k\in {\bf a}$; for such a cyclic 
permutation, there is exactly one decomposition with ${\bf b}$-part 
${\bf b}$, obtained by cyclically shifting the pieces of the decomposition 
$d_{\bf w}$.

\begin{exm}
Let ${\bf w}=(v_3,v_4,v_5,v_6,v_7,v_1,v_2)$ and 
consider the decomposition 
$${\bf w}={\bf a_1b_1a_2b_2a_3}=(v_3,v_4)(v_5)(v_6)(v_7,v_1)(v_2).$$
Then ${\bf b}=\{v_1,v_5,v_7\}$ and ${\bf a}=\{v_2,v_3,v_4,v_6\}$. The 
only cyclic permutations of ${\bf v}=(v_1,\ldots,v_7)$ admitting the 
${\bf b}$-part $\{v_1,v_5,v_7\}$ are the ones starting with $v_k\in {\bf a}$, 
and for each one, there is a unique decomposition determined by ${\bf b}$:
\begin{equation*}
\begin{cases}(v_2,v_3,v_4)(v_5)(v_6)(v_7,v_1)\\
(v_3,v_4)(v_5)(v_6)(v_7,v_1)(v_2)\\
(v_4)(v_5)(v_6)(v_7,v_1)(v_2,v_3)\\
(v_6)(v_7,v_1)(v_2,v_3,v_4)(v_5).
\end{cases}
\end{equation*}
The set ${\bf W}^{\bf b}$ consists of these four decompositions.
\end{exm}
					Let the {\it ordered ${\bf a}$-part} of a decomposition $d_{\bf w}$ of
					a cyclic permutation ${\bf w}$ of ${\bf v}=(v_1,\ldots,v_r)$ be the word 
					${\bf a_1\cdots a_s}$ of the decomposition $d_{\bf w}$. Then by the
					above, there are exactly $|{\bf a}|$ decompositions in ${\bf W}^{\bf b}$,
				and their ordered ${\bf a}$-parts are given by
					\begin{equation}\label{aparts}
			\{\sigma^j_{|{\bf a}|}({\bf a_1\cdots a_s})\,|\,j=0,\ldots,|{\bf a}|-1\}
			\end{equation}
			i.e. the cyclic permutations of the letters of ${\bf a_1\cdots a_s}$.

				\vspace{.2cm}
			\textit{Proof of Proposition \ref{ganitcirc}.} 
			Let $c=\bigl(M(v_1)\,|\,v_1^{n-1}\bigr)$, and
				let $N=pari(M)$ (where $pari$ is defined in 
                Definition \ref{def: various mould operators}), so that $N$ is a polynomial
				mould of fixed homogeneous degree
				$n$, with $N(v_1)=-cv_1^{n-1}$.  Since $M$ is circ-constant for $c$, 
                we have
				\begin{equation}\label{littlec}
			N(v_1,\ldots,v_r)+\cdots +N(v_r,v_1,...,v_{r-1})=(-1)^rc
				\sum_{{{e_1+\cdots+e_r=n-r}\atop{e_i\ge 0}}} v_1^{e_1}\cdots v_r^{e_r}.
				\end{equation}
			By the explicit formula (\ref{ganitQ}), we have
				\begin{equation}\label{ganitpoc}
			\bigl(\overline{ganit}(poc)\cdot N\bigr)(v_1,\ldots,v_r)=\sum_{{\bf W}_{\bf v}}
			poc(\lfloor {\bf b}_1)\cdots poc(\lfloor {\bf b}_s)
				N({\bf a}_1\cdots {\bf a}_s),\\
				\end{equation}
			where we recall that the sum over ${\bf W}_{\bf v}$ is the sum over all
				decompositions $d_{\bf v}={\bf a}_1{\bf b}_1\cdots {\bf a}_s{\bf b}_s$ 
				of ${\bf v}=(v_1,\ldots,v_r)$ as in \eqref{chunks}.
				Now, summing over the set ${\bf W}$ which 
				is the union of such decompositions
				not only for ${\bf v}$ but also for all the cyclic permutations of
				${\bf v}$ (cf.~\eqref{thesetW}), we have
				\begin{align}
			\sum_{i=0}^{r-1} \bigl(\overline{ganit}&(poc)\cdot N\bigr)
				\bigl(\sigma_r^i({\bf v})\bigr) =
                \sum_{{\bf a}_1{\bf b}_1\cdots{\bf a}_s{\bf b}_s\in {\bf W}}
			poc(\lfloor {\bf b}_1)\cdots poc(\lfloor {\bf b}_s)
				N({\bf a}_1\cdots {\bf a}_s)\notag\\
				&=
               \sum_{\substack{{\bf b}\subsetneq\{v_1,\dots,v_r\} }}\sum_{{\bf a}_1{\bf b}_1\cdots{\bf a}_s{\bf b}_s\in {\bf W}^{\bf b}}
				poc(\lfloor {\bf b}_1)\cdots poc(\lfloor {\bf b}_s)
				N({\bf a}_1\cdots {\bf a}_s)\notag\\
				&=
			\sum_{\substack{{\bf b}\subsetneq\{v_1,\dots,v_r\} \\ {\bf b}={\bf b}_1\cup\cdots\cup{\bf b}_s}}
            poc(\lfloor {\bf b}_1)\cdots poc(\lfloor{\bf b}_s) \sum_{j=0}^{|{\bf a}|-1} 
			N\bigl(\sigma_{|{\bf a}|}^j({\bf a}_1\cdots {\bf a}_s)\bigr)\notag\\
				&= c
			\sum_{\substack{{\bf b}\subsetneq\{v_1,\dots,v_r\} \\ {\bf b}={\bf b}_1\cup\cdots\cup{\bf b}_s}}
            (-1)^{|{\bf a}|}
			poc(\lfloor {\bf b}_1)\cdots poc(\lfloor{\bf b}_s) 
				\sum_{{{e_1+\cdots+e_{|{\bf a}|}=n-|{\bf a}|}\atop{e_j\ge 0}}}
			v_{i_1}^{e_1}\cdots v_{i_{|{\bf a}|}}^{e_{|{\bf a}|}},\label{temp1}
			\end{align}
			where the disjoint consecutive chunks ${\bf b}_1,\dots,{\bf b}_s$ in the
				sum run
				over the possible ${\bf b}$-parts of decompositions as in \eqref{chunks}.
				Here, the first equality is the definition of $\overline{ganit}(poc)$, the
				second equality follows directly from \eqref{Wb}, the third follows directly 
				from \eqref{aparts}, and the last equality from \eqref{littlec} using
				the notation ${\bf a}_1\cup\cdots \cup {\bf a}_s=\{v_{i_1},\ldots,v_{i_{|{\bf a}|}}\}$ for the subset ${\bf a}$ of $\{v_1,\ldots,v_r\}$.

\vspace{.2cm}
If $c=0$, i.e.~if $M$ is a circ-neutral mould, the expression \eqref{temp1} 
is trivially equal to zero in all depths $r>1$, proving Proposition
\ref{ganitcirc} in the case where $M$ is circ-neutral.
In order to deal with the case where $M$ is circ-constant for a value $c\ne 0$,
we use a trick and subtract off a known mould that is also circ-constant for $c$.  

\vspace{.2cm}
For $A\subset \{1,\ldots,n\}$, let $S^A_d$ denote the sum of 
all monomials of degree $d$ in the letters $v_i$, $i\in A$. In particular, if $d=0$ we set $S^A_d=1$ and if $d<0$ we set $S^A_d=0$.

\begin{lem}\label{Tnc} For $n>1$ and any constant $c$, let $T^n_c$ be the 
homogeneous polynomial mould of degree $n$ (cf.~Definition \ref{circconstance})
defined by $T^n_c(\emptyset)=0$ and 
$$T^n_c(v_1,\ldots,v_r)=\frac{c}{r}S^{\{1,\ldots,r\}}_{n-r}$$ 
for $r>1$; in particular we have $T^n_c(v_1,\ldots,v_n)=\frac{1}{n}$ 
and $T^n_c(v_1,\ldots,v_r)=0$ when $r>n$.  Then $T^n_c$ is circ-constant and 
$\overline{ganit}(poc)\cdot pari(T^n_c)$ is circ-neutral.  
\end{lem}

			The proof of this lemma is surprisingly long and technical, so we have 
				relegated it to Appendix \ref{Appendix B}.  Using the result, we can now finish the proof of
				Proposition \ref{ganitcirc}. Indeed  we have $M(v_1)=T^n_c(v_1)=cv_1^{n-1}$,
					    so the mould $M-T^n_c$ is circ-neutral and thus
						    $\overline{ganit}(poc)\cdot pari(M-T^n_c)$ is also circ-neutral.  But 
						    Lemma \ref{Tnc} shows that $\overline{ganit}(poc)\cdot pari(T^n_c)$ is itself
						    circ-neutral, so we have
						    $$\sum_{i=1}^r circ^i\bigl(\overline{ganit}(poc)\cdot pari(M)\bigr)=
						    \sum_{i=1}^r circ^i\bigl(\overline{ganit}(poc)\cdot pari(T^n_c)\bigr)=0$$
						    and thus $\overline{ganit}(poc)\cdot pari(M)$ is also circ-neutral, completing
						    the proof of Proposition \ref{ganitcirc}.
						    \ \qed

						    \vspace{.3cm}
			We now proceed to the second step, showing that the operator 
				$Ad_{\overline{ari}}(invpil)$ preserves circ-neutrality.

				\begin{prop}\label{adaripilcirc} If $M\in\overline{ARI}(\mathcal F_\Lau)$ is circ-neutral
				then $Ad_{\overline{ari}}(invpil)\cdot M$ is also circ-neutral.  
				\end{prop}

			\textit{Proof.} By (\ref{adariinvpil}), we have
				\begin{equation}\label{expsumpil}
			Ad_{\overline{ari}}(invpil) = exp\bigl(ad_{\overline{ari}}(-lopil)\bigr)
				= \sum_{n\ge 0}\frac{(-1)^n}{n!}\bigl(ad_{\overline{ari}}(lopil)\bigr)^n.
				\end{equation}
			The definition of $lopil$ in (\ref{deflopil}) shows that
				$lopil$ is trivially circ-neutral.  Thus, since $M$ is circ-neutral  (as explained in Example \ref{example: lopil is circ-neutral}),
				$ad_{\overline{ari}}(lopil)\cdot M=\overline{ari}(lopil,M)$ is also
					circ-neutral by Proposition \ref{circneutclosed}, and successively so are
					all the terms $ad_{\overline{ari}}(lopil)^n(M)$.  Thus 
					$Ad_{\overline{ari}}(invpil)\cdot M$ is circ-neutral.   \qed

\vspace{.4cm}
\noindent \underline{\it Proof of (iii) (continued):} 
Finally, we now assume that $B\in ARI({\mathcal F_{\ser}})_{al+tsen*circconst}$.
In particular, $swap(B)$ is a $\ast$circ-constant polynomial-valued 
mould in $\overline{ARI}({\mathcal F_{\ser}})$ of homogeneous degree $n$.
Let $B_0$ be the 
(unique) constant-valued mould such that $swap(B)+B_0$ is circ-constant.  
Then by Propositions \ref{ganitcirc} and \ref{adaripilcirc}, 
the mould 
$$Ad_{\overline{ari}}(invpil)\cdot \overline{ganit}(poc)\cdot pari 
\bigl(swap(B)+B_0\bigr)$$ 
is circ-neutral. This mould breaks up as the sum 
$$Ad_{\overline{ari}}(invpil)\cdot \overline{ganit}(poc)\cdot pari(swap(B))+ 
Ad_{\overline{ari}}(invpil)\cdot \overline{ganit}(poc)\cdot pari(B_0),$$ 
but the operator $Ad_{\overline{ari}}(invpil)\cdot \overline{ganit}(poc)$ 
preserves constant-valued moulds (cf.~\cite{S2}, Lemma 4.6.2 for the proof).  
Thus the mould 
$$Ad_{\overline{ari}}(invpil)\cdot \overline{ganit}(poc)\cdot pari\bigl(swap(B)+B_0\bigr)=\qquad\qquad$$ 
$$\qquad\qquad\qquad\qquad Ad_{\overline{ari}}(invpil)\cdot \overline{ganit}(poc)\cdot pari(swap(B))+B_0$$
is circ-neutral, or equivalently, 
$$Ad_{\overline{ari}}(invpil)\cdot ganit(poc)\cdot pari(swap(B))$$
is $\ast$circ-neutral. However, using the fact that $pari$ trivially commutes 
with $swap$ and also the fact that by Proposition \ref{keypropb} (which we 
recall relies on \'Ecalle's assertion \eqref{Ecstatement}) 
$Ad_{ari}(invpal)\cdot pari(B)$ is push-invariant, we can apply 
\eqref{goodfund} to find that 
\begin{align*}
Ad_{\overline{ari}}(invpil)\cdot ganit(poc)\cdot swap\bigl(pari(B)\bigr)
&=swap\cdot Ad_{ari}(invpal)\cdot pari(B)\\
&=swap\cdot \Xi(B).
\end{align*}
Thus $swap\cdot \Xi(B)$ is $\ast$circ-neutral, which concludes the proof of (iii).
\qed

\vspace{.3cm}
\noindent \underline{\it Proof of (iv): $A\in ARI(\mathcal F_\Lau)^\Delta$.}  
				\vspace{.3cm}
For any mould $A\in \overline{ARI}({\mathcal F_\Lau})$, let us use the notation 
$$\overline{\Delta}(A)(v_1,\ldots,v_r):=v_1(v_1-v_2)\cdots (v_{r-1}-v_r)v_r\,A(v_1,\ldots,v_r),$$ 
so that we have
\begin{equation}\label{Deltabar}
swap\bigl(\Delta(A)\bigr)=\overline{\Delta}\bigl(swap(A)\bigr).
\end{equation}
We write $\overline{ARI}({\mathcal F_\Lau})^{\overline{\Delta}}$ for the subspace of moulds 
$M\in \overline{ARI}({\mathcal F_\Lau})$ such that $\overline{\Delta}(M)$ is polynomial.  
By \eqref{Deltabar}, we see that
$M\in \overline{ARI}({\mathcal F_\Lau})^{\overline{\Delta}}$ means that 
$swap(M)\in ARI({\mathcal F_\Lau})^\Delta$.

We will again use \'Ecalle's assertion \eqref{Ecstatement} and the equality
\eqref{goodfund}; this time we will study the left-hand side of $\eqref{goodfund}$ to track the denominators that appear 
in the right-hand side.  
By (\ref{goodfund}), if $B$ is a polynomial-valued mould satisfying the 
twisted senary relation, and if $A=\Xi(B)=Ad_{ari}(invpal)\cdot pari(B)$ is
push-invariant, then 
$A$ lies in $ARI({\mathcal F_\Lau})^\Delta$ if and only if 
\begin{equation}\label{wanted1} 
Ad_{\overline{ari}}(invpil)\cdot \overline{ganit}(poc)\cdot swap\cdot pari(B)
 \in \overline{ARI}({\mathcal F_\Lau})^{\overline\Delta}.
\end{equation}

The identity \eqref{goodfund} applies in the situation where $B\in ARI_{al+tsen*circconst}({\mathcal F_\ser})$ and $A=\Xi(B)$, since we proved in (i) that $A$ is indeed push-invariant. 
We will prove that $A\in ARI({\mathcal F_\Lau})^\Delta$ by studying 
the denominators that arise in \eqref{wanted1},
produced first by applying $\overline{ganit}(poc)$ 
and then by applying $Ad_{\overline{ari}}(invpil)$.  
Let us first show that the denominators introduced by applying 
$\overline{ganit}(poc)$ to a polynomial-valued mould are at worst of the form 
$(v_1-v_2)\cdots (v_{r-1}-v_r)$.

 \begin{lem}[$\!\!${\cite[Proposition 4.38]{B}}]\label{ganitdem}  
 Let $M\in{\overline{ARI}}({\mathcal F}_\ser)$.
 Then
$$\overline{ganit}(poc)\cdot M\in \overline{ARI}({\mathcal F_\Lau})^{\overline\Delta}$$ 
and in fact $(v_1-v_2)\cdots (v_{r-1}-v_r)$ is 
{a common denominator}
for all terms arising in $\bigl(\overline{ganit}(poc)\cdot M\bigr)(v_1,\ldots,v_r)$
{for all $M$.}
\end{lem}

\textit{Proof.} The explicit expression for $\overline{ganit}(Q)$ given in 
(\ref{ganitQ}) shows that the only denominators that can occur 
in $\overline{ganit}(poc)\cdot M$ come from the factors 
\begin{equation}\label{dens} 
poc(\lfloor {\bf b}_1)\cdots poc(\lfloor {\bf b}_s) 
\end{equation} 
for all decompositions $d_{\bf v}={\bf a}_1{\bf b}_1\cdots 
{\bf a}_s{\bf b}_2$ of ${\bf v}=(v_1,\ldots,v_r)$ into chunks as in 
(\ref{chunks}), and 
$$\lfloor {\bf b}_i=(v_k-v_{k-1},v_{k+1}-v_{k-1},\ldots,v_{k+l}-v_{k-1})$$ 
(for $k>1$) as in (\ref{flex}).  Since $poc$ is defined as in (\ref{poc}), 
the only factors that can appear in (\ref{dens}) are 
$(v_l-v_{l-1})$ where $v_l$ is a letter in one of ${\bf b}_i$, and these 
factors appear in each term with multiplicity one.  Since 
the sum ranges over all possible decompositions, the only letter of 
${\bf v}$ that never belongs to any ${\bf b}_i$ is $v_1$; all the other 
factors $(v_i-v_{i-1})$ appear.  Thus $(v_1-v_2)(v_2-v_3)\cdots (v_{r-1}-v_r)$ 
is a common denominator for all the terms in the sum defining 
$\overline{ganit}(poc)\cdot M$, which proves the lemma.\ \qed

\vspace{.2cm} 
\begin{lem}\label{circden} 
Let $M,N\in \overline{ARI}({\mathcal F_\Lau})_{*circneut}^{\overline{\Delta}}$ be 
two moulds. 
Then $\overline{ari}(M,N)$ also lies in $ARI({\mathcal F_\Lau})^{\overline{\Delta}}$.  
\end{lem}

\textit{Proof.} We begin by showing that if $M\in\overline{ARI}({\mathcal F_\Lau})_{*circneut}^{\overline{\Delta}}$
then for all $r>1$, $\overline{\Delta}(M)$ satisfies the identity
\begin{equation}\label{id0}
\overline{\Delta}(M)(0,v_2,\ldots,v_r)=\overline{\Delta}(M)(v_2,\ldots,v_r,0).
\end{equation}
In fixed depth $r>1$, the $\ast$circ-neutrality of $M$ means that there exists a 
constant $c_r$ such that $M(v_1,\ldots,v_r)+c_r$ is circ-neutral, i.e.
\begin{equation}\label{id1}
M(v_1,\ldots,v_r)+M(v_2,\ldots,v_r,v_1)+\cdots+M(v_r,v_1,\ldots,v_{r-1})+rc_r=0.
\end{equation}
Writing
$$M(v_1,\ldots,v_r)=\frac{\overline{\Delta}(M)(v_1,\ldots,v_r)}{v_1(v_1-v_2)\cdots(v_{r-1}-v_r)v_r}$$
allows us to rewrite \eqref{id1} in terms of $\overline{\Delta}(M)$, as
\begin{align}\label{id2}
&\frac{\overline{\Delta}(M)(v_1,\ldots,v_r)}{v_1(v_1-v_2)\cdots(v_{r-1}-v_r)v_r}+\frac{\overline{\Delta}(M)(v_2,\ldots,v_r,v_1)}{v_2(v_2-v_3)\cdots(v_r-v_1)v_1}\notag\\
&\qquad +\cdots+\frac{\overline{\Delta}(M)(v_r,v_1,\ldots,v_{r-1})}{v_r(v_r-v_1)\cdots(v_{r-2}-v_{r-1})v_{r-1}}+rc_r=0.
\end{align}
Only the first two terms of this sum have a pole at $v_1=0$, so multiplying
\eqref{id2} by $v_1$ and then setting $v_1=0$ leaves only the first two
terms, which become
\begin{equation}\label{id3}
\frac{\overline{\Delta}(M)(0,v_2,\ldots,v_r)}{-v_2(v_2-v_3)\cdots(v_{r-1}-v_r)v_r}
+\frac{\overline{\Delta}(M)(v_2,\ldots,v_r,0)}{v_2(v_2-v_3)\cdots(v_r)}=0.
\end{equation}
Since the two terms in \eqref{id3} have the same denominator with opposite
signs, this is equivalent to the desired identity \eqref{id0}.

\vspace{.1cm}
We can now proceed to the proof of Lemma \ref{circden}. By additivity, 
it is enough to consider the case where 
$M$ and $N$ are moulds concentrated in depths $r\ge 1$ and $s\ge 1$ 
respectively.  

In \S 4.7 of [S2] (proof of Theorem 4.7.1, see also \S 4.3.4 of Baumard's
thesis or Prop. 5.1 of [BS])
the explicit expressions for $\overline{\Delta}\bigl(\overline{arit}(M)\cdot N\bigr)$,
$\overline{\Delta}\bigl(\overline{arit}(N)\cdot M\bigr)$ and $\overline{\Delta}\bigl(lu(M,N)\bigr)$ are calculated, so
as to examine the poles of each term. It is shown there that 

\vspace{.1cm}\noindent 
(i) $\overline{\Delta}\bigl(\overline{arit}(M)\cdot N\bigr)$ has potential poles (linear
factors in the denominator) only at $v_i-v_{i+r+1}$ 
($1\le i\le s-1$)
and $v_i-v_{i+r}$ 
($1\le i\le s$), and at $v_{r+1}=0$ and $v_s=0$;

\vspace{.1cm}\noindent 
(ii) The residues at the poles $v_i=v_{i+r+1}$ are equal to zero with
only the hypothesis $M,N\in \overline{ARI}({\mathcal F_\Lau})^{\overline{\Delta}}$;

\vspace{.1cm}\noindent 
(iii) The residues at the poles $v_i=v_{i+r}$ are equal to zero thanks
to \eqref{id0};

\vspace{.1cm}\noindent
(iv) $\overline{\Delta}\bigl(\overline{arit}(M)\cdot N\bigr)$ does have simple poles 
at $v_{r+1}=0$ and $v_s=0$, with residues
\begin{equation}\label{residue1}
\overline{\Delta}(N)(0,v_{r+2},\ldots,v_{r+s})\overline{\Delta}(M)(v_1,\ldots,v_r)\ \ \ \hbox{at}\ v_{r+1}=0
\end{equation}
and
\begin{equation}\label{residue2}
\overline{\Delta}(N)(v_1,\ldots,v_{s-1},0)\overline{\Delta}(M)(v_{s+1},
\ldots,v_{r+s})\ \ \ \hbox{at}\ v_s=0.
\end{equation}

\vspace{.1cm}\noindent
The symmetric statements hold for $\overline{\Delta}\bigl(arit(N)\cdot M\bigr)$,
which thus has simple poles only at $v_r=0$ and $v_{s+1}=0$ with residues
\begin{equation}\label{residue3}
\overline{\Delta}(M)(0,v_{s+2},\ldots,v_{r+s})\overline{\Delta}(N)(v_1,\ldots,v_s)\ \ \ \hbox{at} \ v_{s+1}=0
\end{equation}
and
\begin{equation}\label{residue4}
\overline{\Delta}(M)(v_1,\ldots,v_{r-1},0)\overline{\Delta}(N)(v_{r+1},
\ldots,v_{r+s}) \ \ \ \hbox{at}\ v_r=0.
\end{equation}

\vspace{.2cm}\noindent
Since we have 
\begin{equation}\label{arireminder}
\overline{\Delta}\bigl(\overline{ari}(M,N)\bigr)=
\overline{\Delta}\bigl(\overline{arit}(N)\cdot M\bigr)
-\overline{\Delta}\bigl(\overline{arit}(M)\cdot N\bigr)+\overline{\Delta}\bigl(lu(M,N)\bigr),
\end{equation}
we can show that $\overline{ari}(M,N)\in \overline{ARI}({\mathcal F_\Lau})^{\bar\Delta}$ by showing that
$\overline{\Delta}\bigl(\overline{ari}(M,N)\bigr)$ is polynomial, i.e.~that it has
no poles, by showing firstly that $\overline{\Delta}\bigl(lu(M,N)\bigr)$ has no poles outside of $v_r=0$, $v_s=0$, $v_{r+1}=0$, $v_{s+1}=0$, and
secondly that the residues at these poles cancel out exactly with those
of $\overline{\Delta}\bigl(\overline{arit}(N)\cdot M-\overline{arit}(M)\cdot N\bigr)$.
We write out $\overline{\Delta}\bigl(lu(M,N)\bigr)$ explicitly as
\begin{align*}\label{Deltalu}
&\overline{\Delta}\bigl(lu(M,N)\bigr)(v_1,\dots,v_{r+s}) \\
&=v_1(v_1-v_2)\cdots(v_{r+s-1}-v_{r+s})v_{r+s} \\
&\ \ \ \ \cdot\bigl(M(v_1,\ldots,v_r)N(v_{r+1},\ldots,v_{r+s})
-N(v_1,\ldots,v_s)M(v_{s+1},\ldots,v_{r+s})\bigr)\notag\\
&=\frac{v_r-v_{r+1}}{v_rv_{r+1}}\overline{\Delta}(M)(v_1,\ldots,v_r)\overline{\Delta}(N)(v_{r+1},\ldots,v_{r+s})\notag\\
&\ \ \ -\frac{v_s-v_{s+1}}{v_sv_{s+1}}\overline{\Delta}(N)(v_1,\ldots,v_s)\overline{\Delta}(M)(v_{s+1},\ldots,v_{r+s}).
\end{align*}
This shows that the only poles are at $v_r=0$, $v_s=0$, $v_{r+1}=0$, 
$v_{s+1}=0$. Multiplying by $v_{r+1}$ and setting $v_{r+1}=0$, we
find that the residue at $v_{r+1}=0$ comes from the first line
only and is equal to \eqref{residue1} which is the residue at $v_{r+1}=0$
in $\overline{\Delta}\bigl(\overline{arit}(M)\cdot N\bigr)$. Thus these poles
cancel out in \eqref{arireminder}. Similarly, the pole at $v_s=0$
has residue equal to \eqref{residue2}, the pole at $v_{s+1}=0$ has
residue equal to the negative of \eqref{residue3}, and the pole at $v_r=0$ 
has residue equal to the negative of \eqref{residue4}. So all these
poles cancel out in \eqref{arireminder}, and we conclude that
$\overline{\Delta}\bigl(\overline{ari}(M,N)\bigr)\in \overline{ARI}({\mathcal F_\ser})$,
i.e.~$\overline{ari}(M,N)\in \overline{ARI}({\mathcal F_\Lau})^{\overline{\Delta}}$, as desired.\ \qed

\begin{cor}\label{usefulcor} 
If $P\in \overline{ARI}({\mathcal F_\Lau})$ is a circ-neutral mould 
lying in $\overline{ARI}({\mathcal F_\Lau})^{\overline{\Delta}}$, 
then we also have 
\begin{equation}\label{Ad} 
Ad_{\overline{ari}}(invpil)\cdot P\in \overline{ARI}(\mathcal F_\Lau)^{\overline{\Delta}}.  
\end{equation}
\end{cor}

\textit{Proof.}
The mould $lopil\in \overline{ARI}(\mathcal F_\Lau)$  given in \eqref{deflopil} is circ-neutral 
(as explained in Example \ref{example: lopil is circ-neutral})
and 
by its defining equation \eqref{deflopil}, $lopil$ lies in 
$\overline{ARI}(\mathcal F_\Lau)^{\overline{\Delta}}$. 
Therefore by 
Lemma \ref{circden}, we have $\overline{ari}(lopil,P)\in 
\overline{ARI}({\mathcal F_\Lau})^{\overline{\Delta}}$. Furthermore, by Proposition \ref{circneutclosed},
$\overline{ari}(lopil,P)$ is also circ-neutral. 
Then, applying Lemma \ref{circden} successively shows that 
$ad_{\overline{ari}}(lopil)^n(P)$ is a circ-neutral mould lying in 
$\overline{ARI}({\mathcal F_\Lau})^{\overline{\Delta}}$ for all $n\ge 1$.
Since $Ad_{\overline{ari}}(invpil)\cdot P$ is obtained by summing these
terms by (\ref{expsumpil}), we obtain (\ref{Ad}).\ \qed

\vspace{.4cm}
\noindent \underline{\it Proof of (iv) (continued):} 
We can now complete the proof of (iv) of Proposition \ref{keypropb}. Recall
that $B\in ARI({{\mathcal F}_\ser})_{al+tsen*circconst}$ and $A=\Xi(B)$. We may assume
that $B$ is homogeneous of fixed degree $n\ge 3$. 
Then there exists a constant mould $B_0$ such that $swap(B+B_0)$ is circ-constant.
By Lemma \ref{ganitdem},
we have
\begin{equation}\label{thismould}
\overline{ganit}(poc)\cdot swap\cdot pari(B+B_0)\in \overline{ARI}({\mathcal F_\Lau})^{\overline\Delta}.
\end{equation}
Since $swap$ and $pari$ commute, the mould in \eqref{thismould} is equal to
\begin{equation}\label{thismould2}
\overline{ganit}(poc)\cdot pari\cdot swap(B+B_0)\in \overline{ARI}({\mathcal F_\Lau})^{\overline\Delta}.
\end{equation}
By Proposition \ref{ganitcirc}, this mould 
is circ-neutral, and by \eqref{thismould2} it
lies in $\overline{ARI}({\mathcal F_\Lau})^{\overline\Delta}$.
Therefore we can apply
Corollary \ref{usefulcor} with $P=\overline{ganit}(poc)\cdot swap\cdot pari(B+B_0)$
to conclude that
$$Ad_{\overline{ari}}(invpil)\cdot \overline{ganit}(poc)\cdot
swap\cdot pari(B+B_0)\in \overline{ARI}({\mathcal F_\Lau})^{\overline\Delta}.$$
Observe that this mould is equal to
$$Ad_{\overline{ari}}(invpil)\cdot \overline{ganit}(poc)\cdot
swap\cdot pari(B)({u_1,\ldots,u_n})+(-1)^nB_0({u_1,\ldots,u_n})$$
since $pari\cdot swap(B_0)({u_1,\ldots,u_n})=(-1)^nB_0({u_1,\ldots,u_n})$ and $Ad_{ari}(invpil)\cdot \overline{ganit}(poc)$
preserves constant moulds (cf.~\cite{S2}, Lemma 4.6.2)
which follows from $\overline{ganit}(pic)^{-1}=\overline{ganit}(poc)$
in \eqref{eq: ganit pic and poc}.
Thus
$$Ad_{\overline{ari}}(invpil)\cdot \overline{ganit}(poc)\cdot
swap\cdot pari(B)\in \overline{ARI}({\mathcal F_\Lau})^{\overline\Delta}.$$
Swapping this mould and applying \eqref{goodfund} with 
$N=pari(B)$, we finally find that
$$Ad_{ari}(invpal)\cdot pari(B)=\Xi(B)\in ARI({\mathcal F_\Lau})^\Delta,$$
which completes the proof of (iv), and therefore the proof of
Proposition \ref{keypropb}, and thus Step 3 of the proof of Theorem
\ref{krvsection}.

\subsubsection{Step 4: Composing with $\Delta$.}
\label{subsubsec: Step 4: Comparison}
The final step in the proof is very easy; it consists simply of composing
the injective map \eqref{themap}
$$
\Xi:ARI({\mathcal F_\ser})_{al+tsen*circconst}\rightarrow
ARI({\mathcal F_\Lau})^\Delta_{al+push*circneut}
$$
with the injective map $\Delta$ to finally obtain a map
$$
\Delta\circ\Xi:ARI({\mathcal F_\ser})_{al+tsen*circconst}\rightarrow
\Delta(ARI({\mathcal F_\Lau})^\Delta_{al+push*circneut})
$$
As explained just before paragraph \ref{Wkrv},
since by definition the left-hand space is the image of the space $W_{\lie{krv}}$
under the mould map $ma$, where $W_{\lie{krv}}$ is isomorphic to $\lie{krv}$, and
the right-hand space is just the image of $\lie{krv}_{ell}$ under
the mould map $ma$ (see diagram \eqref{defkrvell}),
we thus obtain an injective map
$$\lie{krv}\rightarrow\lie{krv}_{ell}$$
as shown by the sequence of vertical maps in the 
right-hand side of the diagram \eqref{bigone}
below. This completes the proof of Theorem \ref{krvsection}.\ \qed

\vspace{.3cm}
\subsection{Relations with elliptic Grothendieck-Teichm\"uller and 
double shuffle Lie algebras (Theorem \ref{bigdiagram})}\label{doubleshuf}
The final result in this paper is the proof of 
Theorem \ref{bigdiagram}.  In fact, this result is simply a consequence of
putting together the results of the previous sections with known results.
Indeed, the commutativity of the diagram
\begin{equation}\label{littlediagram}\xymatrix{
\lie{grt}\ar@{^{(}->}[rr]\ar@{^{(}->}[d]&&\lie{ds}\ar@{^{(}->}[d]\\
\widetilde{\lie{grt}}_{ell}\ar@{^{(}->}[rr]\ar@{^{(}->}[dr]&&\lie{ds}_{ell}
\ar@{^{(}->}[dl]\\
&\lie{oder}_2&
}\end{equation}
where $Ad_{ari}(invpal):\lie{ds}\rightarrow \lie{ds}_{ell}$ is the right-hand
vertical map is shown in \cite{S3}.

Let $b=b(x,y)\in \lie{ds}$. By (10), the injective map 
$\lie{ds}\hookrightarrow\lie{krv}$
sends $b$ to the derivation 
of ${\rm Lie}[x,y]$ given by $y\mapsto \hat b(x,y):=b(-x-y,-y)$ 
and $[x,y]\mapsto 0$
(which determines the value of the derivation on $x$ uniquely).
If $b(x,y)\in \lie{ds}$, then $b(x,-y)$ lies in $W_{\lie{krv}}$ and
$b(z,-y)$ lies in $V_{\lie{krv}}$, so this map unpacks to
$$\xymatrix{\lie{ds}\ar@{->}[r]^{y\mapsto -y}&W_{\lie{krv}}
\ar[r]^{\!\!x\mapsto z}&V_{\lie{krv}}\ar[r]& \lie{krv},}$$ 
where the last map comes from (\ref{krviso}).
We can thus construct a commutative square
\begin{equation}\label{arrow}\xymatrix{
\lie{ds}\ar[d]\ar@{->}[r]&\lie{krv}\ar@{->}[d]\\
\lie{ds}_{ell}\incl[r]&\lie{krv}_{ell}
}\end{equation}
given in detail by
\begin{equation}\label{bigone}\xymatrix{
\lie{ds}\ar@{^{(}->}[rr]^{y\mapsto -y}\ar[d]_{ma}&&W_{\lie{krv}}\simeq \lie{krv}\ar[d]^{ma}\\
ARI({\mathcal F_\ser})_{\underline{al}*\underline{il}}\ar@{^{(}->}[rr]^{pari}\ar[d]_{Ad_{ari}(invpal)}&&ARI({\mathcal F_\ser})_{al+tsen*circconst}\ar@{->}[d]^{Ad_{ari}(invpal)\circ pari}\\
ARI({\mathcal F_\Lau})^\Delta_{\underline{al}*\underline{al}}\incl[rr]\ar[d]_{\Delta}&&ARI({\mathcal F_\Lau})^\Delta_{al+push*circneut}\ar[d]^{\Delta}\\
\Delta(ARI({\mathcal F_\Lau})^\Delta_{\underline{al}*\underline{al}})\ar[d]_{ma^{-1}}\incl[rr]&& \Delta(ARI({\mathcal F_\Lau})^\Delta_{al+push*circneut})\ar[d]^{ma^{-1}}\\
\lie{ds}_{ell}\incl[rr]&&\lie{krv}_{ell}.
}\end{equation}
The second line of this diagram 
is the direct mould translation of the top line, as the left-hand space
is exactly $ma(\lie{ds})$, the right-hand space is $ma(W_{\lie{krv}})$ by
(\ref{Wkrvsen}), and the map $pari$ restricted to Lie series is
nothing other than $y\mapsto -y$. The proof of the vertical morphism
$$
Ad_{ari}(invpal):ARI({\mathcal F_\ser})_{\underline{al}*\underline{il}}\rightarrow ARI({\mathcal F_\Lau})^\Delta_{\underline{al}*\underline{al}}
$$
has two parts: the fact that $Ad_{ari}(invpal)$ maps $ARI({\mathcal F_\Lau})_{\underline{al}*\underline{il}}$
to $ARI({\mathcal F_\Lau})_{\underline{al}*\underline{al}}$ is one of the fundamental results of \'Ecalle's mould theory,
and follows directly from \'Ecalle's fundamental identity \eqref{Ecallefund} (see \cite{S2}, 
Theorem 4.6.1), while the fact that restricted to 
$ARI(\mathcal F_\ser)_{\underline{al}*\underline{il}}$, the operator $Ad_{ari}(invpal)$ produces denominators
at worst $\Delta$ was proved in \cite{B}, Thm. 4.35.
The vertical morphism
$${
Ad_{ari}(invpal)\circ pari:ARI(\mathcal F_\ser)_{al+tsen*circconst}
\hookrightarrow
ARI(\mathcal F_\Lau)^\Delta_{al+push*circneut}
}
$$
{follows directly from our assumption \eqref{Ecstatement}.}
Since $pari$ is an involution,
the square formed by the second and third lines of the diagram commutes,
where the horizontal inclusion of the third line comes from 
Corollary \ref{mouldfirstcor}.
The fourth horizontal inclusion is obtained 
simply by applying $\Delta$ to the third line.
Finally, the last line of the diagram
comes from the definitions $ma(\lie{ds}_{ell})=\Delta(ARI({\mathcal F_\Lau})^\Delta_{\underline{al}*\underline{al}})$
($\!\!$\cite{S3}) and $ma(\lie{krv}_{ell})=\Delta(ARI({\mathcal F_\Lau})^\Delta_{al+push*circneut})$
by Definition \ref{defnkrvell}.

This diagram shows that the diagram (\ref{littlediagram}) above can be
completed by the diagram (\ref{arrow}) to the commutative diagram of Theorem 
\ref{bigdiagram}.
\qed

\vspace{.3cm}
\appendix

\vspace{.3cm}
\section{Proof of Lemma \ref{Tnc}}\label{Appendix B}

Let us recall the statement of the technical lemma \ref{Tnc}. Recall that
for $A\subset \{v_1,\ldots,v_r\}$, we let $M^A_d$ denote the set of
all monomials of degree $d$ in the letters of $A$, and $S^A_d$ the 
sum of all monomials in $M^A_d$. Recall from the notation given
between the statement of Proposition \ref{ganitcirc} and its proof
that for any cyclic permutation
permutation $(v_i,\ldots,v_r,v_1,\ldots,v_{i-1})$ of
$(v_1,\ldots,v_r)$ and any decomposition of it as
${\bf a}_1{\bf b}_1\cdots {\bf a}_s{\bf b}_s$ for some $s\ge 1$ (where 
none of these subsets is allowed to be empty except for possibly ${\bf b}_s$),
we write ${\bf b}={\bf b}_1\cup\cdots\cup{\bf b}_s$ and call this the
${\bf b}$-part of the decomposition. 
Recall also the notation ${\bf W}$, ${\bf W}^{\bf b}$ etc.~given just
after the statement of Proposition \ref{ganitcirc}. 
We will also use the following further notation:
for each $0\le i\le r$, let ${\mathcal{B}}_i$ denote the set of all 
${\bf b}$-parts (occurring in the sum in \eqref{temp3}) that contain $v_i$ 
but not $v_{i+1}, \ldots,v_r$; in other words, a given ${\bf b}$-part 
${\bf b}={\bf b}_1\cup\cdots\cup{\bf b}_s$ 
lies in ${\mathcal{B}}_i$ if and only if $i$ is the largest index such that 
$v_i$ occurs in ${\bf b}$.  

The following are examples for $i=0,1,2$. 
\begin{exm}
We have ${\mathcal{B}}_0=\{\emptyset\}$,
since a ${\bf b}$-part that contains no $v_i$ can only be the empty set; 
empty ${\bf b}$-parts arise in the trivial decompositions
$$\sigma^j_r({\bf v})=(v_{j+1},\ldots,v_r,v_1,\ldots,v_j)={\bf a}_1$$
for $0\le j\le r-1$.  The set ${\mathcal{B}}_1$ contains only the single 
element ${\bf b}=(v_1)$, and corresponds to the decompositions
$$\sigma^j_r({\bf v})=(v_{j+1},\ldots,v_r,v_1,\ldots,v_j)=(v_{j+1},\ldots,v_r)
(v_1)(v_2,\ldots,v_j)={\bf a}_1{\bf b}_1{\bf a}_2$$
for $1\le j\le r-1$ (in fact just ${\bf a}_1{\bf b}_1$ for $j=1$).  
The set ${\mathcal{B}}_2$ contains two different 
${\bf b}$-parts, namely $(v_2)$ and $(v_1,v_2)$. The ${\bf b}$-part
$(v_2)$ occurs in the decompositions
\begin{equation*}
\begin{cases}
(v_1,\ldots,v_r)=(v_1)(v_2)(v_3,\ldots,v_r)={\bf a}_1{\bf b}_1{\bf a}_2\\
(v_3,\ldots,v_r,v_1,v_2)=(v_3,\ldots,v_r,v_1)(v_2)={\bf a}_1{\bf b}_1\\
(v_{j+1},\ldots,v_r,v_1,\ldots,v_{j})=(v_{j+1},\ldots,v_r,v_1)(v_2)(v_3,\ldots,v_j)={\bf a}_1{\bf b}_1{\bf a}_2\ \ {\rm for}\ 3\le j\le r-1.
\end{cases}
\end{equation*}
The ${\bf b}$-part $(v_1,v_2)$ occurs in the decompositions
\begin{equation*}
\begin{cases}
(v_3,\ldots,v_r)(v_1,v_2)={\bf a}_1{\bf b}_1\\
(v_j,\ldots,v_r)(v_1,v_2)(v_3,\ldots,v_{j-1})={\bf a}_1{\bf b}_1{\bf a}_2\ \ {\rm for}\ 4\le j\le r\\
\end{cases}
\end{equation*}
Indeed, for $1\le i\le r-1$, the set ${\mathcal{B}}_i$ is simply in bijection 
with the set of all subsets $B\subset \{1,\ldots,i-1\}$, by associating
$B$ to the ${\bf b}$-part $\{v_j|j\in B\}\cup \{v_i\}$); when $i=r$,
${\mathcal{ B}}_r$ is in bijection with the set of all strict subsets of 
$\{1,\ldots,r-1\}$.
\end{exm}

\noindent {\bf Lemma \ref{Tnc}.} {\it For $n>1$ and any constant $c\ne 0$, 
let $T^n_c$ be the homogeneous polynomial mould of degree $n$ defined by 
$$T^n_c(v_1,\ldots,v_r)=\frac{c}{r}\,S^{\{v_1,\ldots,v_r\}}_{n-r}.$$
Then $T^n_c$ is circ-constant and 
$\overline{ganit}(poc)\cdot pari(T^n_c)$ is circ-neutral.}

\vspace{.3cm}
\textit{Proof.} The mould $T^n_c$ is trivially circ-constant for the 
value $c$. For the rest of this proof we set $c=1$ and $T^n=T^n_1$;
it suffices to multiply all identities in the proof below
by the constant $c$ to prove the general case.

Let $N=pari(T^n)$. In order to show that $\overline{ganit}(poc)\cdot N$
is circ-neutral, we start by recalling from the beginning of the
proof of Proposition \ref{ganitcirc} that for each $r>1$, the cyclic sum 
\begin{equation}\label{ToProve}
\overline{ganit}(poc)\cdot N(v_1,\ldots,v_r)+\cdots+
\overline{ganit}(poc)\cdot N(v_r,v_1,\ldots,v_{r-1})
\end{equation}
is equal to the expression \eqref{temp1}
\begin{equation}\label{temp3}
\sum_{\substack{{\bf b}\subsetneq\{v_1,\dots,v_r\} \\ {\bf b}={\bf b}_1\cup\cdots\cup{\bf b}_s}}
(-1)^{|{\bf a}|}
poc(\lfloor{\bf b}_1)\cdots poc(\lfloor{\bf b}_s)S_{n-|{\bf a}|}^{\bf a},
\end{equation}
where the sum runs over all the distinct ${\bf b}$-parts that can arise from
decomposing the cyclic permutations 
$\sigma^i_r({\bf v})=(v_{i+1},\ldots,v_r,v_1,
\ldots,v_i)$ into chunks ${\bf a}_1{\bf b}_1\cdots {\bf a}_s{\bf b}_s$
in which only ${\bf b}_s$ can be empty (so in particular ${\bf b}={\bf b}_1\cup\cdots\cup{\bf b}_s$ cannot be 
the full set $\{v_1,\ldots,v_r\}$). For each term of the sum, ${\bf a}$ 
denotes the subset of $\{v_1,\ldots,v_r\}$ which is the complement of the
${\bf b}$-part ${\bf b}_1\cup\cdots\cup{\bf b}_s$.

To prove the Lemma, we will show that \eqref{temp3} is equal to zero for 
all $r>1$ by breaking up the sum into simpler parts that can be 
expressed explicitly. 
To do this, we observe that every 
${\bf b}$-part corresponding to a decomposition of a cyclic permutation 
$\sigma^j_r({\bf v})$ lies in a unique ${\mathcal{B}}_i$.
Therefore if we set
\begin{equation}\label{Ri}
R^r_i:=\sum_{{\bf b}\in {\mathcal{B}}_i} (-1)^{|{\bf a}|}
poc(\lfloor{\bf b}_1)\cdots poc(\lfloor{\bf b}_s)S_{n-|{\bf a}|}^{\bf a}
\end{equation}
for $0\le i\le r$, we can write the sum \eqref{temp3} as 
\begin{equation}\label{sumTis}
\sum_{\substack{{\bf b}\subsetneq\{v_1,\dots,v_r\} \\ {\bf b}={\bf b}_1\cup\cdots\cup{\bf b}_s}}
(-1)^{|{\bf a}|}
poc(\lfloor{\bf b}_1)\cdots poc(\lfloor{\bf b}_s)S_{n-|{\bf a}|}^{\bf a}=
R^r_0+\cdots+R^r_r.
\end{equation}
\vspace{.2cm}
We have
\begin{equation}\label{R0}
R^r_0 =(-1)^rS_{n-r}^{\{v_1,\ldots,v_r\}}
=\sum_{j=0}^{r-1} N\bigl(\sigma^j_r({\bf v})\bigr),
\end{equation}
where the first equality comes from \eqref{Ri} and the second from the fact 
that $N=pari(T^n)$ and $T^n$ is circ-constant. For $R^r_1$, the only
possible ${\bf b}$-part is $(v_1)$ and we have
\begin{equation}\label{R1r}
R^r_1={{(-1)^r}\over{v_1-v_r}}S_{n-r+1}^{\{v_2,\ldots,v_r\}}.
\end{equation}
Let us now consider $R_i^r$ for $i>1$. Consider any decomposition
\begin{equation}\label{decomp}
{\bf a}_1 {\bf b}_1\cdots {\bf a}_s {\bf b}_s
\end{equation}
of any cyclic permutation of $(v_1,\ldots,v_r)$. For any non-empty 
chunk ${\bf b}_j$ of this decomposition, write
${\bf b}_j=(v_{k+1},v_{k+2},\ldots,v_l)$, with 
indices $k$ and $l$ considered mod $r$ from $1$ to $r$ (for example 
${\bf b}_j=(v_{r-1},v_r,v_1)$ with $k=r-2$ and $l=1$).
Then by the definition of $poc$ (cf.~\eqref{poc}), we have
\begin{align*}
poc(\lfloor {\bf b}_j)&=poc(v_{k+1}-v_k,v_{k+2}-v_k,\ldots,v_l-v_k)\\
&={{-1}\over{(v_{k+1}-v_k)(v_{k+1}-v_{k+2})\cdots (v_{l-1}-v_l)}}\\
&=\prod_{v_m\in {\bf b}_j} {{1}\over{(v_{m-1}-v_m)}},
\end{align*}
again with indices $m$ mod $r$ with values from $1$ to $r$.
Thus, writing as usual ${\bf a}$ for the ${\bf a}$-part of a decomposition 
as in \eqref{decomp}, \eqref{Ri} can be written
\begin{equation}\label{si}
R^r_i=\sum_{{\bf b'}\subseteq \{v_1,\ldots,v_{i-1}\}} {{(-1)^{|{\bf a}|}S_{n-|{\bf a}|}^{\bf a}}\over{\prod_{v_j\in {\bf b}}(v_{j-1}-v_j)}}
\end{equation}
for $1\le i\le r-1$, 
where ${\bf b}'$ runs over all subsets of $\{v_1,\ldots,v_{i-1}\}$ so
${\bf b}={\bf b}'\cup\{v_i\}$ runs over the elements of ${\mathcal{B}}_i$,
and for $i=r$ we have
\begin{equation}\label{sir}
R^r_r=\sum_{{\bf b'}\subsetneq \{v_1,\ldots,v_{r-1}\}} {{(-1)^{|{\bf a}|}S_{n-|{\bf a}|}^{\bf a}}\over{\prod_{v_j\in {\bf b}}(v_{j-1}-v_j)}}.
\end{equation}

We will use these explicit expressions in the proofs of the following
two claims, which will allow us to easily conclude the proof of Lemma \ref{Tnc}.
The proofs of the claims are given farther below.

\vspace{.3cm}
\noindent {\bf Claim 1.} {\it (i) For $i=1$, we have
\begin{equation}\label{Claim1i}
R_1^r={{(-1)^{r-1}S_{n-r+1}^{\{v_2,\ldots,v_r\}}}\over{(v_r-v_1)}}.
\end{equation}
\noindent (ii) For $2\le i\le r-1$, we have
\begin{equation}\label{Claim1ii}
R_i^r={{(-1)^{r-i}S^{\{v_{i-1},v_{i+1},\ldots,v_{r-1}\}}_{n-r+i}}\over{(v_r-v_1)(v_1-v_2)\cdots (v_{i-1}-v_i)}}.
\end{equation}
\noindent (iii) For $i=r$, we have
\begin{equation}\label{Claim1iii}
R_r^r={{v_{r-1}^{n-1}}\over{(v_r-v_1)(v_1-v_2)\cdots (v_{r-2}-v_{r-1})}}.
\end{equation}}

\vspace{.3cm}
\noindent {\bf Claim 2.} {\it For $0\le i\le r-1$ we have 
\begin{equation}\label{Claim2}
R_0^r+\cdots+R_i^r={{(-1)^{r-i}S^{\{v_i,\ldots,v_{r-1}\}}_{n-r+i}}\over{(v_r-v_1)(v_1-v_2)\cdots (v_{i-1}-v_i)}}.
\end{equation}
}

\vspace{.3cm}
Accepting these two claims, we observe that by \eqref{Claim2} when $i=r-1$, we have
\begin{equation}
R_0^r+\cdots+R_{r-1}^r={{-v_{r-1}^{n-1}}\over{(v_r-v_1)(v_1-v_2)\cdots
(v_{r-2}-v_{r-1})}},
\end{equation}
but this equal to $-R_r^r$ by Claim 1 (iii). 
This shows that we have $R_0^r+\cdots+R_r^r=0$. Thus the left-hand
side of \eqref{sumTis} is equal to zero. But this left-hand side is
equal to \eqref{temp3}, which in turn is equal to \eqref{ToProve},
showing that $\overline{ganit}(poc)\cdot N$ is circ-neutral. 
This completes the proof of Lemma \ref{Tnc}. \qed

\vspace{.3cm}
It remains only to prove Claims 1 and 2. We begin by proving Claim 2 using
Claim 1, then give an elementary (but surprisingly complicated) proof of
Claim 1.

\vspace{.3cm}
\noindent {\bf Proof of Claim 2.} 
We first note the following trivial but useful identity. 
Recall the notation $V_m=\{v_1,\ldots,v_m\}$ for $1\le m\le r$.
Let $A'\subsetneq V_r$, let $v_j\not\in A'$, and
let $A=A' \cup \{v_j\}$: then we have the useful identity
\begin{equation}\label{cute1}
S^{A'}_{d+1}+v_jS^A_d=S^A_{d+1}.
\end{equation}
Indeed, the first term is the sum of all monomials of degree $d+1$ in the
elements of $A'$, and the second is the sum of all monomials of degree $d+1$
in the letters of $A$ that contain $v_j$, so their sum forms the sum of all 
monomials of degree $d+1$ in the letters of $A$.

We will prove Claim 2 by induction on $i$.  The base case
$i=0$ is given by \eqref{R0}. 
Now let $1\le i\le r-1$ and assume \eqref{Claim2} up to $i-1$.  Then by the 
induction hypothesis and Claim 1, we have
\begin{align*}\label{temp2}R_0^r+&\cdots+R_i^r=(R_0^r+\cdots+R_{i-1}^r)+R_i^r\\ 
\\ 
&={{(-1)^{r-i+1}S^{\{v_{i-1},\ldots,v_{r-1}\}}_{n-r+i-1}}\over{
(v_r-v_1)(v_1-v_2)\cdots (v_{i-2}-v_{i-1})}}+{{(-1)^{r-i}S_{n-r+i}^{\{v_{i-1},v_{i+1},
\ldots,v_{r-1}\}}}\over{(v_r-v_1)(v_1-v_2)\cdots(v_{i-1}-v_i)}}\\ 
\\ 
&={{(-1)^{r-i+1}\Bigl(v_{i-1}S_{n-r+i-1}^{\{v_{i-1},\ldots,v_{r-1}\}}-v_iS_{n-r+i-1}^{\{v_{i-1}, \ldots,v_{r-1}\}}-S_{n-r+i}^{\{v_{i-1},v_{i+1},\ldots,v_{r-1}\}}\Bigr)
}\over{(v_r-v_1)(v_1-v_2)\cdots(v_{i-1}-v_i)}}.
\end{align*}
By \eqref{cute1}, the second term and third terms in the numerator sum to
$-S_{n-r+i}^{\{v_{i-1},\ldots,v_{r-1}\}}$, so the numerator becomes
$(-1)^{r-i+1}\Bigl(v_{i-1}S_{n-r+i-1}^{\{v_{i-1},\ldots,v_{r-1}\}}-S_{n-r+i}^{\{v_{i-1},\ldots,v_{r-1}\}}\Bigr)$ which again by \eqref{cute1} sums to
$(-1)^{r-i}\Bigl(S_{n-r+i}^{\{v_i,\ldots,v_{r-1}\}}\Bigr)$. Thus we have
$$R_0^r+\cdots+R_i^r={{(-1)^{r-i}S_{n-r+i}^{\{v_i,\ldots,v_{r-1}\}}}\over{(v_r-v_1)(v_1-v_2)\cdots(v_{i-1}-v_i)}},$$
which proves Claim 2.  \qed

\vspace{.5cm}
Finally, we proceed to the proof of Claim 1.

\vspace{.3cm}
\noindent {\bf Proof of Claim 1.} (i) When $i=1$ we have ${\mathcal{B}}_1=\{v_1\}$,
so here there is only one term in the sum \eqref{si} corresponding to
${\bf b}'=\emptyset$, ${\bf b}=\{v_1\}$, ${\bf a}=\{v_2,\ldots,v_r\}$,
$|{\bf a}|=r-1$, so that \eqref{si} for $i=1$ comes down to \eqref{Claim1i}.
 
\vspace{.2cm}
\noindent (ii) Let $V_m=\{v_1,\ldots,v_m\}$ for $1\le m\le r$. Fix a value
of $i$ with $2\le i\le r-1$.
Recall that $R_i^r$ was given in \eqref{si} as the sum
\begin{equation}\label{siagain}
R^r_i=\sum_{{\bf b'}\subseteq V_{i-1}} {{(-1)^{|{\bf a}|}S_{n-|{\bf a}|}^{\bf a}}\over{\prod_{v_j\in {\bf b}}(v_{j-1}-v_j)}}.
\end{equation}
where ${\bf b}={\bf b}'\cup \{v_i\}$ and ${\bf a}$ is the complement
of ${\bf b}$ in $V_r$.

Multiplying $R_i^r$ by the common denominator $(v_r-v_1)\cdots (v_{i-1}-v_i)$
and setting $v_0=v_r$ as usual, 
we rewrite \eqref{siagain} as
\begin{equation}\label{cute2}
\prod_{j=1}^i (v_{j-1}-v_j)R_i^r=
\sum_{{\bf b}'\subseteq V_{i-1}} (-1)^{|{\bf a}|}
\prod_{v_j\in V_{i-1}\setminus {\bf b}'}
(v_{j-1}-v_j)S_{n-|{\bf a}|}^{{\bf a}}.
\end{equation}

To conclude the proof, we need one more claim. 

\vspace{.2cm}
\noindent {\bf Claim 3.} {\it For each pair $i,k$ with $1<i<r$ and
$1\le k\le i-1$, define the polynomial $Q_k^i$ by
\begin{equation*}
\sum_{v_1,\ldots,v_k\not\in B'\subset V_{i-1}} (-1)^{r-|B'|+k-1}\biggl(\prod_{v_j\in V_{i-1}\setminus (B'\cup\{v_1,\ldots,v_k\})} \!\!\!\!\!\!\!\!\!\!(v_{j-1}-v_j)\biggr)
S^{V_r\setminus (B'\cup \{v_1,\ldots,v_{k-1}\}\cup \{v_i,v_r\})}_{n-r+|B'|+k+1}
\end{equation*}
Let $Q_0^i$ denote the right-hand side of \eqref{cute2}. Then we have
$$Q_0^i=Q_1^i=Q_2^i=\cdots=Q_{i-1}^i.$$}

\vspace{.1cm}
\noindent {\bf{Proof of claim 3}}. Let us show that for $0\le k\le i-2$ 
we have $Q_k^i=Q_{k+1}^i$.  We break the expression for $Q_k^i$
into the terms with $v_{k+1}\in B'$ and those with $v_{k+1}\not\in B'$, 
writing $Q_k^i$ as
$$\sum_{{{v_1,\ldots,v_k\not\in B'\subseteq V_{i-1}}\atop{v_{k+1}\in B'}}} 
\!\!\!\!\!\!\!\!\!\!(-1)^{r-|B'|+k-1}\Bigl(\prod_{v_j\in V_{i-1}\setminus (B'\cup\{v_1,\ldots,v_k\})} \!\!\!\!\!\!\!\!\!\!\!\!\!\!\!(v_{j-1}-v_j)\Bigr) S_{n-r+|B'|+k+1}^{V_r\setminus (B'\cup\{v_1,\ldots,v_{k-1}\}\cup\{v_i,v_r\})}$$
$$+\!\!\!\sum_{v_1,\ldots,v_{k+1}\not\in B'\subseteq V_{i-1}} \!\!\!\!\!\!\!\!\!\!\!\!(-1)^{r-|B'|+k-1}\Bigl(\prod_{v_j\in V_{i-1}\setminus (B'\cup\{v_1,\ldots,v_k\})} \!\!\!\!\!\!\!\!\!\!\!\!\!\!\!(v_{j-1}-v_j)\Bigr)S_{n-r+|B'|+k+1}^{V_r\setminus (B'\cup\{v_1,\ldots,v_{k-1}\}\cup \{v_i,v_r\})}.$$
We will use an analogous decomposition for $Q_0^i$, breaking it up as
\begin{align*}
\prod_{j=1}^i (v_{j-1}-v_j)R_i^r&=
\sum_{v_1\in B'\subseteq V_{i-1}} (-1)^{|A|}\prod_{v_j\in V_{i-1}\setminus B'} (v_{j-1}-v_j)S_{n-|A|}^{A}\\
&\ \ \ \ \ + \sum_{v_1\not\in B'\subseteq V_{i-1}} (-1)^{|A|}\prod_{v_j\in V_{i-1}\setminus B'} (v_{j-1}-v_j)S_{n-|A|}^{A},
\end{align*} 
where $B'={\bf b}'$ and $A={\bf a}$ (to harmonize the notation with the
$Q_k^i$), and the $v_0$ that occurs in the second line is equal to $v_r$ as
in \eqref{cute2}.

Now fix $k\in \{0,\ldots,i-2\}$, and write $B'':=B'\setminus\{v_{k+1}\}$ in 
the first line of the decomposition of $Q_k^i$,
and simply replace the notation $B'$ by $B''$ in the second line, obtaining
$$\!\!\!\!\!\sum_{v_1,\ldots,v_{k+1}\not\in B''\subseteq V_{i-1}} 
\!\!\!\!\!\!\!\!\!\!(-1)^{r-|B''|+k}\Bigl(\prod_{v_j\in V_{i-1}\setminus (B''\cup\{v_1,\ldots,v_{k+1}\})} \!\!\!\!\!\!\!\!\!\!\!\!\!\!\!(v_{j-1}-v_j)\Bigr) S_{n-r+|B''|+k+2}^{V_r\setminus (B''\cup\{v_1,\ldots,v_{k-1}\}\cup\{v_{k+1},v_i,v_r\})}$$
$$+\!\!\!\!\!\sum_{v_1,\ldots,v_{k+1}\not\in B''\subseteq V_{i-1}} \!\!\!\!\!\!\!\!\!\!(-1)^{r-|B''|+k-1}(v_k-v_{k+1})\times\qquad\qquad\qquad\qquad
\qquad\qquad\qquad\qquad$$
$$\qquad\qquad\qquad\qquad \Bigl(\prod_{v_j\in V_{i-1}\setminus (B''\cup\{v_1,\ldots,v_{k+1}\})} \!\!\!\!\!\!\!\!\!\!\!\!\!\!(v_{j-1}-v_j)\Bigr)S_{n-r+|B''|+k+1}^{V_r\setminus (B''\cup\{v_1,\ldots,v_{k-1}\}\cup \{v_i,v_r\})},$$
except when $k=0$, where (using $|A|=r-|B''|-1$), the $S$-factors are
$S_{n-r+|B''|+1}^{V_r\setminus (B''\cup \{v_1,v_i\}}$ in the first line and 
$S_{n-r+|B''|+1}^{V_r\setminus (B''\cup \{v_i\})}$ in the second.

Now we gather the terms, writing this as
$$Q_k^i=\sum_{v_1,\ldots,v_{k+1}\not\in B''\subseteq V_{i-1}} 
(-1)^{r-|B''|+k}\Bigl(\prod_{v_j\in V_{i-1}\setminus (B''\cup\{v_1,\ldots,v_{k+1}\})} (v_{j-1}-v_j)\Bigr) \times \qquad\qquad\qquad$$
$$\qquad\Biggl(S_{n-r+|B''|+k+2}^{V_r\setminus (B''\cup\{v_1,\ldots,v_{k-1}\}\cup\{v_{k+1},v_i,v_r\}))}
-(v_k-v_{k+1})S_{n-r+|B''|+k+1}^{V_r\setminus (B''\cup\{v_1,\ldots,v_{k-1}\}\cup \{v_i,v_r\})}\Biggr),$$
or in the case $k=0$,
$$Q_0^i=\sum_{B''\subseteq V_{i-1}} 
(-1)^{r-|B''|}\Bigl(\prod_{v_j\in V_{i-1}\setminus (B''\cup\{v_1\})} (v_{j-1}-v_j)\Bigr) \times \qquad\qquad\qquad$$
$$\qquad\Biggl(S_{n-r+|B''|+2}^{V_r\setminus (B''\cup\{v_1,v_i\}))}
-(v_r-v_1)S_{n-r+|B''|+1}^{V_r\setminus (B''\cup\{v_i\})}\Biggr),$$
We will now use \eqref{cute1} twice to simplify the right-hand factor.
For $1\le k\le i-2$, we take $A'=V_r\setminus (B''\cup\{v_1,\ldots,v_{k-1}\}\cup
\{v_{k+1},v_i,v_r\})$ and $A=A'\cup \{v_{k+1}\}$, while for
$k=0$ we take $A'=V_r\setminus (B''\cup \{v_1,v_i\})$ and $A=A'\cup \{v_1\}$.
We can then expand the right-hand factor as
$$S_{n-r+|B''|+k+2}^{A'}
-v_kS_{n-r+|B''|+k+1}^A
+v_{k+1}S_{n-r+|B''|+k+1}^A$$
(recalling that when $k=0$, $v_0=v_r$).
Applying \eqref{cute1}, this simplifies to
$$S_{n-r+|B''|+k+2}^A -v_kS_{n-r+|B''|+k+1}^A.$$
Next, since $v_k\not\in B''$ we must have $v_k\in A$, so by
applying \eqref{cute1} again we see that this simplifies to
$$S_{n-r+|B''|+k+2}^{A\setminus \{v_k\}}
=S_{n-r+|B''|+k+2}^{V_r\setminus (B''\cup\{v_1,\ldots,v_k\}\cup \{v_i,v_r\})}.$$
Thus we can rewrite $Q_k^i$ as
$$\sum_{v_1,\ldots,v_{k+1}\not\in B''\subseteq V_{i-1}} 
(-1)^{r-|B''|+k}\biggl(\prod_{v_j\in V_{i-1}\setminus (B''\cup\{v_1,\ldots,v_{k+1}\})} (v_{j-1}-v_j)\biggr) 
S_{n-r+|B''|+k+2.}^{V_r\setminus (B''\cup\{v_1,\ldots,v_k\}\cup \{v_i,v_r\})}$$
But according to the definition of the polynomials $Q_k^i$ for $1\le k\le i-1$,
this is exactly equal to $Q_{k+1}^i$. Thus we find that
$Q_0^i=Q_1^i=\cdots =Q_{i-1}^i$, completing the proof of Claim 3.\hfill{$\square$}

\vspace{.3cm}
It remains only to prove part (iii) of Claim 1.

\vspace{.1cm}
\noindent (iii) In this final part we have to prove that
\begin{equation}\label{extra}
\prod_{j=1}^{r-1} (v_{j-1}-v_j)R_r^r=v_{r-1}^{n-1}.
\end{equation}
Recall from \eqref{sir} that $R_r^r$ is given by the formula
\begin{equation*}
R^r_r=\sum_{{\bf b'}\subsetneq V_{r-1}} {{(-1)^{|{\bf a}|}S_{n-|{\bf a}|}^{\bf a}}\over{\prod_{v_j\in {\bf b}}(v_{j-1}-v_j)}},
\end{equation*}
where ${\bf b}={\bf b}'\cup \{v_r\}$. Thus the common denominator
of all the terms in the sum is $(v_r-v_1)(v_1-v_2)\cdots (v_{r-1}-v_r)$,
and we have
\begin{equation}\label{cute3}\prod_{j=1}^r (v_{j-1}-v_j)R_r^r=
\sum_{{\bf b}'\subsetneq V_{r-1}} (-1)^{r-|{\bf b}'|-1}\prod_{v_j\in V_{r-1}\setminus {\bf b}'} 
(v_{j-1}-v_j)\ S_{n-r+|{\bf b}'|+1}^{V_{r-1}\setminus {\bf b}'}.
\end{equation}
Let us write ${\bf c}=V_{r-1}\setminus {\bf b}'$, so this equality can be
expressed as
\begin{equation}\label{cute4}\prod_{j=1}^r (v_{j-1}-v_j)R_r^r=
\sum_{\emptyset\ne {\bf c}\subseteq V_{r-1}} (-1)^{|{\bf c}|}\prod_{v_j\in {\bf c}}
(v_{j-1}-v_j)S_{n-|{\bf c}|}^{{\bf c}}.
\end{equation}
For $1\le i\le r-1$ and $n\ge 1$, define the sum $T_i^n$ by
\begin{equation*}
T^n_i:=\sum_{\emptyset\ne {\bf c}\subseteq V_i}(-1)^{|{\bf c}|}\prod_{v_j\in {\bf c}} (v_{j-1}-v_j)S_{n-|{\bf c}|}^{{\bf c}},
\end{equation*}
where we set $S_0^{{\bf c}}=1$ and $S_m^{{\bf c}}=0$ if $m<0$. 
By this definition, the term $T_{r-1}^n$ is equal to the right-hand side
of \eqref{cute4}.  We will prove that 
\begin{equation}\label{cute5}
T_i^n=(v_i-v_r)v_i^{n-1}.
\end{equation}
The equality \eqref{cute5} suffices to prove the desired result \eqref{extra}. 
Indeed, since $T_{r-1}^n$ is equal to the right-hand side of \eqref{cute4},
the left-hand side of \eqref{cute4} is equal to the
right-hand side of \eqref{cute5} with $i=r-1$, i.e.
$$\prod_{j=1}^r(v_{j-1}-v_j)R_r^r=T^n_{r-1}=(v_{r-1}-v_r)v_{r-1}^{n-1}.$$
Canceling out the factor $(v_r-v_{r-1})$ from both sides yields 
the desired identity \eqref{extra}.

It remains only to prove \eqref{cute5}. We proceed by induction on $i$.  
When $i=1$, we have 
${\bf c}=\{v_1\}$ and for all $n\ge 1$, we have
$$T_1^n=-(v_r-v_1)S^{v_1}_{n-1}=(v_1-v_r)v_1^{n-1},$$
proving the base case.

Fix $n\ge 1$ and assume \eqref{cute5} holds for $i-1$.
We break $T^n_i$ into the sum over ${\bf c}$ containing $v_i$ and ${\bf c}$ not
containing $v_i$, as follows:
\begin{align*}
T^n_i&=\sum_{\emptyset\ne {\bf c}\subseteq V_{i-1}}(-1)^{|{\bf c}|}\prod_{v_j\in {\bf c}} (v_{j-1}-v_j)S_{n-|{\bf c}|}^{\bf c}\\
&\qquad \ +\sum_{{{{\bf c}\subseteq V_{i-1}}\atop{{\bf c}'={\bf c}\cup\{v_i\}}}}(-1)^{|{\bf c}'|}(v_{i-1}-v_i)\prod_{v_j\in {\bf c}}
(v_{j-1}-v_j)S_{n-|{\bf c}'|}^{{\bf c}'}\\
&=T^n_{i-1} +(v_{i-1}-v_i)\sum_{{\bf c}\subseteq V_{i-1}}(-1)^{|{\bf c}|+1}\prod_{v_j\in {\bf c}} (v_{j-1}-v_j)S_{n-|{\bf c}|-1}^{{\bf c}\cup \{v_i\}}\\
&=T^n_{i-1}-(v_{i-1}-v_i)v_i^{n-1}-(v_{i-1}-v_i)\sum_{\emptyset\ne{\bf c}\subseteq V_{i-1}}(-1)^{|{\bf c}|}\prod_{v_j\in {\bf c}} (v_{j-1}-v_j)S_{n-|{\bf c}|-1}^{{\bf c}\cup \{v_i\}},
\end{align*}
where the last line comes from separating the sum over ${\bf c}\subseteq V_{i-1}$ into ${\bf c}=\emptyset$ (giving the extra term $(v_{i-1}-v_i)v_i^{n-1}$)
and the sum over ${\bf c}\ne \emptyset$.
Since ${\bf c}$ does not contain $v_i$, we can write
$$S_{n-|{\bf c}|-1}^{{\bf c}\cup\{v_i\}}=S_{n-|{\bf c}|-1}^{{\bf c}}+v_iS_{n-|{\bf c}|-2}^{{\bf c}}+v_i^2S_{n-|{\bf c}|-3}^{{\bf c}}+\cdots
+v_i^{n-|{\bf c}|-2}S_1^{{\bf c}}+v_i^{n-|{\bf c}|-1}.$$
Using this, the above equality becomes
\begin{align*}
&=T^n_{i-1} -(v_{i-1}-v_i)v_i^{n-1}-(v_{i-1}-v_i)\sum_{\emptyset\ne{\bf c}\subseteq V_{i-1}}(-1)^{|{\bf c}|}\prod_{v_j\in {\bf c}} (v_{j-1}-v_j)\times\\
&\qquad\ \Bigl(S_{n-|{\bf c}|-1}^{{\bf c}}
+v_iS_{n-|{\bf c}|-2}^{{\bf c}}+v_i^2S_{n-|{\bf c}|-3}^{{\bf c}}+\cdots
+v_i^{n-|{\bf c}|-2}S_1^{{\bf c}}+v_i^{n-|{\bf c}|-1}\Bigr)\\
&=T^n_{i-1} -(v_{i-1}-v_i)v_i^{n-1}-(v_{i-1}-v_i)\sum_{\emptyset\ne{\bf c}\subseteq V_{i-1}}(-1)^{|{\bf c}|}\sum_{k=0}^{n-|{\bf c}|-1}\prod_{v_j\in {\bf c}} (v_{j-1}-v_j)S^{{\bf c}}_{n-|{\bf c}|-1-k}v_i^k\\
&=T^n_{i-1} -(v_{i-1}-v_i)v_i^{n-1}-(v_{i-1}-v_i)\sum_{k=0}^{n-2}v_i^k\sum_{\emptyset\ne{\bf c}\subseteq V_{i-1}}(-1)^{|{\bf c}|}\prod_{v_j\in {\bf c}} (v_{j-1}-v_j)S^{{\bf c}}_{n-|{\bf c}|-1-k}.
\end{align*}
If $n-|{\bf c}|-1-k<0$ then $S^{\bf c}_{n-|c|-1-k}=0$, so this is
equal to
$$T^n_{i-1} -(v_{i-1}-v_i)v_i^{n-1}-(v_{i-1}-v_i)\sum_{k=0}^{n-2}v_i^kT_{i-1}^{n-k-1}$$
which then by induction is equal to
\begin{align*}
&=(v_{i-1}-v_r)v_{i-1}^{n-1} -(v_{i-1}-v_i)v_i^{n-1}-(v_{i-1}-v_i)\sum_{k=0}^{n-2}v_i^k(v_{i-1}-v_r)v_{i-1}^{n-k-2}\\
&=(v_{i-1}-v_r)v_{i-1}^{n-1} -(v_{i-1}-v_i)v_i^{n-1}-(v_{i-1}-v_i)(v_{i-1}-v_r)\sum_{k=0}^{n-2}v_i^kv_{i-1}^{n-k-2}\\
&=(v_{i-1}-v_r)v_{i-1}^{n-1} -(v_{i-1}-v_i)v_i^{n-1}-(v_{i-1}-v_r)(v_{i-1}^{n-1}-v_i^{n-1})\\
&=(v_i-v_r)v_i^{n-1}.
\end{align*}
This proves \eqref{cute5} and thus completes the proof of Claim 1 (iii).  \ \qed


\begin{thebibliography}{9}
\bibitem[AT]{AT}{A. Alekseev, Ch. Torossian, {\it The Kashiwara-Vergne conjecture and Drinfeld's associators}, Annals of Math., {\bf 175} (2012), no. 2, 415-463.} 
\bibitem[B]{B}{S. Baumard, {\it Aspects modulaires et elliptiques des relations entre multiz\^etas}, Ph.D. thesis (2014).}
\bibitem[BS]{BS}{S. Baumard, L. Schneps, {\it On the derivation representation of the fundamental Lie algebra of mixed elliptic motives}, Ann. Math. Qu\'ebec {\bf 41} (1) (2014), 43-62.}
\bibitem[Bb]{Bbk} N. Bourbaki, {\it \'El\'ements de math\'ematique. Fasc. XXXVII. Groupes et alg\`ebres de Lie. 
Chapitre II: Alg\`ebres de Lie libres. Chapitre III: Groupes de Lie,} Actualit\'es Scientifiques et Industrielles, 
No. 1349. Hermann, Paris, 1972. 
\bibitem[Br]{Br}{F. Brown, {\it Depth-graded motivic multiple zeta values}, Compositio Math. {\bf 157} (2021), no.3, 529-572.}
\bibitem[Ec1]{Ec81}
{J.~\'Ecalle, {\it Les fonctions r\'{e}surgentes. Tome I et II},
Publications Math\'{e}matiques d'Orsay {\bf 81},
6. Universit\'{e} de Paris-Sud, D\'{e}partement de Math\'{e}matique, Orsay, 1981.} 
\bibitem[Ec2]{Ec}{J.~\'Ecalle, {\it The flexion structure and dimorphy: flexion units, singulators, generators, and the enumeration of multizeta irreducibles}, CRM Series, {\bf 12} (2011), 27-211.}
\bibitem[Ec3]{Ec2}{J.\'Ecalle, {\it Eupolars and their bialternality grid. }, Acta Math. Vietnam. {\bf 40} (2015), no.4, 545-636.}
\bibitem[E1]{E1}{B. Enriquez, {\it Elliptic Associators.}, Selecta Math. (N.S.), {\bf 20} (2014), no.2, 491-584.}
\bibitem[E2]{E2}{B. Enriquez, {\it Analogues elliptiques des nombres multiz\^etas}, Bulletin de la SMF {\bf 144} (2016) no.3, 395-427.}
\bibitem[EF1]{EF}{B. Enriquez, H. Furusho, {\it A stabilizer interpretation of double shuffle Lie algebras}, IMRN {\bf 2018} (2017), no.22, 6870-6907.}
\bibitem[EF2]{EF2}{B. Enriquez, H. Furusho, {\it Double shuffle Lie algebra and special derivations}, preprint {\tt arXiv:2505.02265}, 2025}
\bibitem[F1]{F1}{H. Furusho, {\it Pentagon and hexagon equations},  Ann. of Math. (2), {\bf 171} (2010), no.1, 545-556.}
\bibitem[F2]{F2}{H. Furusho, {\it Double shuffle relation for associators},  Ann. of Math. (2), {\bf 174} (2011), no.1, 341-360.}
\bibitem[FHK]{FHK}{H. Furusho, M.  Hirose, N. Komiyama, {\it 
Associators in mould theory},
preprint, {\tt arXiv:2312.15423}}.
\bibitem[FK1]{FK}{H. Furusho, N. Komiyama, {\it Kashiwara-Vergne and dihedral bi-graded Lie algebras in mould theory},
Ann. Fac. Sci. Toulouse Math. (6) {\bf 32} (2023), no. 4, 655-725.}
\bibitem[FK2]{FK2}{H. Furusho, N. Komiyama, {\it Notes on Kashiwara-Vergne and double shuffle Lie algebras}, "Low Dimensional Topology and Number Theory " Springer Proc. Math. Stat., 456, Springer, Singapore, 2025, 63–80.}
\bibitem[G1]{G1}{A. Goncharov, {\it The dihedral Lie algebras and Galois symmetries of $\pi_1^{(l)}({\mathbb P}^1-(\{0,\infty\}\cup\mu_N))$},
Duke Math. J. {\bf 110} (2001), no. 3, 397--487.}
\bibitem[G2]{G}{A. Goncharov, {\it Galois symmetries of fundamental groupoies and noncommutative geometry}, Duke Math. {\bf 128} (2005), no.2, 209-284.}
\bibitem[I]{I}{Y. Ihara, {\it On the stable derivation algebra associated with some braid groups}, Israel J. Math.{\bf 80} (1992), no.1-2, 135-153.}
\bibitem[IKZ]{IKZ}{K. Ihara, M. Kaneko, D. Zagier, {\it Derivation and double shuffle relations for multiple zeta values}, Compos. Math., {\bf 142} (2006), no.2, 307-338.}
\bibitem[Ka]{Ka}{H. Kawamura, 
{\it Ecalle's senary relation and dimorphic structures},
preprint, {\tt arXiv:2509.21252}, 2025.}
\bibitem[Ko]{K}{N. Komiyama, {\it On properties of $adari(pal)$ and $ganit_v(pic)$}, preprint, {\tt arXiv:2110.04834v2}, 2021.}
\bibitem[M]{M}{M. Maassarani, {\it Bigraded Lie algebras related to multiple zeta values}, Publ. Res. Inst. Math. Sci. {\bf 58} (2022), no. 4, 757--791.}
\bibitem[R]{R}{G. Racinet, {\it S\'eries g\'en\'eratrices non-commutatives de polyz\^etas et associateurs de Drinfel'd}, Ph.D. thesis, 2000.}
\bibitem[S1]{S1} {L. Schneps, {\it Double shuffle and Kashiwara-Vergne Lie algebras.},  J. Algebra, {\bf 367} (2012), 54-74.}
\bibitem[S2]{S2}{L. Schneps, {\it ARI, GARI, Zig and Zag: An introduction to Ecalle's theory of multiple zeta values}, {\tt arXiv:1507.01534v3}, 2015.}
\bibitem[S3]{S3}{L. Schneps, {\it Elliptic multiple zeta values, Grothendieck-Teichm\"{u}ller and mould theory}, Annales Math. Qu\'ebec {\bf 44} (2020) no. 2, 261-289.}
\bibitem[S4]{S4}{L. Schneps, {\it The double shuffle Lie algebra injects into the Kashiwara-Vergne Lie algebra}, preprint, {\tt arXiv:2504.14293}, 2025.}
\bibitem[SS]{SS}{A. Salerno, L. Schneps, {\it Mould theory and the double shuffle Lie algebra structure}, in {\it Periods in Quantum Field Theory and Arithmetic}, J.~Burgos Gil, K.~Ebrahimi-Fard, H.~Gangl, eds., Springer PROMS {\bf 314}, 2020,399-430.}
\bibitem[Ts]{Ts}{H. Tsunogai, {\it On some derivations of Lie algebras related to Galois representations. }, Publ. Res. Inst. Math. Sci. ,{\bf 31} (1995), no.1, 113-134.}

\end{thebibliography}
 \end{document}